\documentclass[11pt, letterpaper, leqno]{amsart} 

\usepackage{graphicx} 
\usepackage{amsmath,amsthm,amssymb} 
\usepackage{amsfonts} 
\usepackage{mathrsfs} 
\usepackage[T1]{fontenc} 
\usepackage[utf8]{inputenc} 

\usepackage[dvipsnames, svgnames]{xcolor} 
\usepackage[shortlabels]{enumitem} 
\usepackage[colorlinks,citecolor=citec,hypertexnames=false,urlcolor=urlc,breaklinks=true, pagebackref]{hyperref} 

\usepackage{accents} 

\usepackage{amsbsy} 
\usepackage{amscd} 
\usepackage{latexsym} 
\usepackage{txfonts} 

\usepackage{exscale} 
\usepackage{bbm} 
\usepackage{mathtools} 
\usepackage{bbding}

\usepackage{geometry} 

\makeatletter
\providecommand\@dotsep{4.5}
\makeatother

\usepackage{thmtools} 

\newif\ifsoldark
\newif\ifsollight
\newif\ifclassic
\newif\ifplain

\usepackage{cite} 

\renewcommand*\backref[1]{\ifx#1\relax \else (p. #1) \fi} 

\numberwithin{equation}{section}

\declaretheoremstyle[
headfont=\color{lmcolor}\normalfont\bfseries,
bodyfont=\normalfont\itshape,
sibling={equation},
]{colorlemma}
\declaretheorem[style=colorlemma]{lemma}

\declaretheoremstyle[
headfont=\color{propcolor}\normalfont\bfseries,
bodyfont=\normalfont\itshape,
sibling={equation},
]{colorprop}
\declaretheorem[style=colorprop]{proposition}
\declaretheorem[style=colorprop]{corollary}

\declaretheoremstyle[
headfont=\color{thmcolor}\normalfont\bfseries,
bodyfont=\normalfont\itshape,
sibling={equation},
]{colorthm}
\declaretheorem[style=colorthm]{theorem}

\declaretheoremstyle[
headfont=\color{defcolor}\normalfont\bfseries,
bodyfont=\normalfont,
sibling={equation},
]{colordef}
\declaretheorem[style=colordef]{definition}

\declaretheoremstyle[
headfont=\color{excolor}\normalfont\bfseries,
bodyfont=\normalfont,
sibling={equation},
]{colorexample}

\declaretheoremstyle[
headfont=\color{rmcolor}\normalfont\itshape,
bodyfont=\normalfont,
sibling={equation},
]{colorremark}
\declaretheorem[style=colorremark]{remark}
\declaretheorem[style=colorremark]{notation}

\declaretheorem[style=colorremark]{conclusion}

\makeatletter
\renewcommand{\eqref}[1]{\textup{\eqreftagform@{\ref{#1}}}}
\let\eqreftagform@\tagform@
\def\tagform@#1{%
	\maketag@@@{\color{eqcolor}(\ignorespaces#1\unskip\@@italiccorr)}%
}
\makeatother

\newcommand{\dint}{\int\!\!\!\!\!\int}

\newcommand{\dist}{\operatorname{dist}}
\newcommand{\card}{\operatorname{card}}
\newcommand{\opint}{\operatorname{int}}

\newcommand{\dv}{\operatorname{div}}
\newcommand{\n}[1]{\mathscr{#1}}

\newcommand{\bb}[1]{\mathbb{#1}}

\newcommand{\dH}{\,\text{d}\n H^{n-1}|_{\Sigma}}
\newcommand{\medcup}{\textstyle\bigcup}
\newcommand{\medsum}{\textstyle\sum}
\newcommand{\mres}{\mathbin{\vrule height 1.6ex depth 0pt width 0.2ex\vrule height 0.2ex depth 0pt width 0.7ex}} 

\newcommand{\RNum}[1]{\uppercase\expandafter{\romannumeral #1\relax}}

\DeclareMathOperator{\supp}{supp}

\DeclareMathOperator{\diam}{diam}

\def\XXint#1#2#3{{\setbox0=\hbox{$#1{#2#3}{\int}$}
		\vcenter{\hbox{$#2#3$}}\kern-.5\wd0}}

\def\Yint#1{\mathchoice
	{\YYint\displaystyle\textstyle{#1}}%
	{\YYint\textstyle\scriptstyle{#1}}%
	{\YYint\scriptstyle\scriptscriptstyle{#1}}%
	{\YYint\scriptscriptstyle\scriptscriptstyle{#1}}%
	\!\dint}
\def\YYint#1#2#3{{\setbox0=\hbox{$#1{#2#3}{\iint}$}
		\vcenter{\hbox{$#2#3$}}\kern-.51\wd0}}
\def\longdash{{-}\mkern-3.5mu{-}} 

\def\fiint{\Yint\longdash}

\newcommand{\ep}{\varepsilon}
\newcommand{\ra}{\rightarrow}
\newcommand{\lra}{\longrightarrow}
\newcommand{\m}[1]{\mathcal{#1}}

\newcommand{\f}[1]{\mathfrak{#1}}
\newcommand{\vertiii}[1]{{\left\vert\kern-0.25ex\left\vert\kern-0.25ex\left\vert #1 
		\right\vert\kern-0.25ex\right\vert\kern-0.25ex\right\vert}}
\newcommand{\dADR}{$d$\textup{-ADR }}
\newcommand{\loc}{\operatorname{loc}}

\newcommand{\dyadic}{\textnormal{dyadic}}

\newcommand{\doubling}{\operatorname{doubling}}
\newcommand{\short}{\operatorname{short}}

\allowdisplaybreaks

\setlist{nosep} 

\plaintrue 

 \ifsoldark
 
 \definecolor{solblack}{HTML}{001116} 	
 \definecolor{solnormal}{HTML}{b3cacc}	
 \definecolor{solblue}{HTML}{268bd2}
 \definecolor{solpurple}{HTML}{6c71c4}
 \definecolor{solteal}{HTML}{2aa198}
 \definecolor{solgreen}{HTML}{859900}
 \definecolor{solpink}{HTML}{ffaaff}
 \definecolor{solbrown}{HTML}{b58900}
 \definecolor{solorange}{HTML}{cb4b16}
 \definecolor{solviolet}{HTML}{8d95ff}
 
 \colorlet{citec}{solblue}
 \colorlet{urlc}{solteal}
 \colorlet{toc}{solpink}
 \colorlet{hyperc}{solviolet}
 \colorlet{bpcolor}{solblue}
 \colorlet{smcolor}{solpurple}
 \colorlet{impcolor}{solpink}
 \colorlet{eqcolor}{solblue}
 \colorlet{lmcolor}{solbrown}
 \colorlet{propcolor}{solpurple}
 \colorlet{thmcolor}{solviolet}
 \colorlet{defcolor}{solblue}
 \colorlet{rmcolor}{solorange}
 \colorlet{excolor}{solgreen}
 \pagecolor{solblack}
 {\color{solnormal} 
 	
 	\fi
 	
 	\ifsollight
 	
 	\definecolor{solbackl}{HTML}{fdf6e3}
 	\definecolor{solnormall}{HTML}{3f4d52} 
 	\definecolor{solbluel}{HTML}{268bd2}
 	\definecolor{solblueld}{HTML}{1b679a} 
 	\definecolor{solpurplel}{HTML}{6c71c4}
 	\definecolor{solteall}{HTML}{2aa198}
 	\definecolor{solgreenl}{HTML}{859900}
 	\definecolor{solpinkl}{HTML}{ffaaaa}
 	\definecolor{solbrownl}{HTML}{b58900}
 	\definecolor{solorangel}{HTML}{cb4b16}
 	\definecolor{solvioletl}{HTML}{8d95ff}

 	\colorlet{citec}{solbluel}
 	\colorlet{urlc}{solteall}
 	\colorlet{toc}{solbluel}
 	\colorlet{hyperc}{solbluel} 
 	\colorlet{bpcolor}{solbluel}
 	\colorlet{smcolor}{solpurplel}
 	\colorlet{impcolor}{solpinkl}
 	\colorlet{eqcolor}{solbluel} 
 	\colorlet{lmcolor}{solbrownl}
 	\colorlet{propcolor}{solpurplel}
 	\colorlet{thmcolor}{solblueld}
 	\colorlet{defcolor}{solbluel}
 	\colorlet{rmcolor}{solorangel}
 	\colorlet{excolor}{solgreenl}
 	\pagecolor{solbackl}
 	{\color{solnormall} 
 		
 		\fi
 		
 		\ifclassic
 		
 		\definecolor{niceblue}{HTML}{003399}
 		\colorlet{citec}{blue}
 		\colorlet{urlc}{teal}
 		\colorlet{toc}{NavyBlue}
 		\colorlet{hyperc}{Green}
 		\colorlet{bpcolor}{NavyBlue}
 		\colorlet{smcolor}{purple}
 		\colorlet{impcolor}{ProcessBlue}
 		\colorlet{eqcolor}{ForestGreen}
 		\colorlet{lmcolor}{Mahogany}
 		\colorlet{propcolor}{RoyalPurple}
 		\colorlet{thmcolor}{NavyBlue}
 		\colorlet{defcolor}{Cerulean}
 		\colorlet{rmcolor}{MidnightBlue}
 		\colorlet{excolor}{PineGreen}
 		
 		\fi

 		\ifplain
 		
 		\colorlet{citec}{blue}
 		\colorlet{urlc}{blue}
 		\colorlet{toc}{blue}
 		\colorlet{hyperc}{blue}
 		\colorlet{bpcolor}{NavyBlue}
 		\colorlet{smcolor}{purple}
 		\colorlet{impcolor}{ProcessBlue}
 		\colorlet{eqcolor}{black}
 		\colorlet{lmcolor}{black}
 		\colorlet{propcolor}{black}
 		\colorlet{thmcolor}{black}
 		\colorlet{defcolor}{black}
 		\colorlet{rmcolor}{black}
 		\colorlet{excolor}{black}
 		
 		\fi

\usepackage{everypage} 
\usepackage{tikzpagenodes} 

\begin{document}
 
\author{S. Mayboroda}

\address{Svitlana Mayboroda
	\\
	School of Mathematics
	\\
	University of Minnesota
	\\
	Minneapolis, MN 55455, USA} \email{svitlana@math.umn.edu}

\author{B. Poggi}

\address{Bruno Poggi Cevallos
	\\
	School of Mathematics
	\\
	University of Minnesota
	\\
	Minneapolis, MN 55455, USA} \email{poggi008@umn.edu}

\thanks {
The first author is supported in part by the NSF grants DMS 1344235, DMS 1839077 and the Simons Foundation grant 563916, SM}

\title[Perturbations of operators on domains with  low dimensional boundaries]{Carleson perturbations of elliptic operators on domains with low dimensional boundaries}
\date{\today}

\begin{abstract}We prove an analogue of a perturbation result for the Dirichlet problem of divergence form elliptic operators by Fefferman, Kenig and Pipher, for the degenerate elliptic operators of David, Feneuil and Mayboroda, which were developed to study geometric and analytic properties of sets with boundaries whose co-dimension is higher than $1$. These operators are of the form $-\text{div} A\nabla$, where $A$ is a weighted elliptic matrix   crafted to weigh the distance to the high co-dimension boundary in a way that allows for the nourishment of an elliptic theory. When this boundary is a $d-$Alhfors-David regular set in $\bb R^n$ with $d\in[1,n-1)$ and $n\geq3$, we prove that the membership of the harmonic measure in $A_{\infty}$ is preserved under Carleson measure perturbations of the matrix of coefficients, yielding in turn that the $L^p-$solvability of the Dirichlet problem is also stable under these perturbations (with possibly different $p$).  If the Carleson measure perturbations are suitably small, we establish solvability of the Dirichlet problem in the same $L^p$ space. One of the corollaries of our results together with a previous result of David, Engelstein and Mayboroda, is that, given \emph{any} \dADR boundary $\Gamma$ with $d\in[1,n-2)$, $n\geq3$, there is a family of degenerate operators of the form described above whose harmonic measure is absolutely continuous with respect to the $d-$dimensional Hausdorff measure on $\Gamma$. 
\end{abstract}

\maketitle

{
	\hypersetup{linkcolor=toc}
	\tableofcontents
}
\hypersetup{linkcolor=hyperc}

\section{Introduction}

The study of the Dirichlet problem for the Laplace equation on domains with rough boundaries  has been of great interest in the last fifty years. This investigation has been in part motivated by a broader program of enlightening a connection between the well-posedness of the Dirichlet problem and the underlying geometry, essentially as pioneered by Dahlberg in \cite{dah1}  and Jerison and Kenig \cite{jk1}, with the important precursor work of F. and M. Riesz in \cite{rr}, and others. 

The key machinery allowing the aforementioned link is the construction via the maximum principle of a family of probability measures, known as the harmonic measure, which provides a representation formula for solutions to the Dirichlet problem with continuous data on the boundary. Then, in a robust way, quantifiable well-posedness of the Dirichlet problem is equivalent to quantifiable absolute continuity of the harmonic measure with respect to the boundary surface measure, whenever the latter makes sense \cite{fkp, dkp, djk, zh}. In turn, this property of quantifiable absolute continuity of harmonic measure has been successfully tied to quantifiable geometric and topological properties of the boundary of the domain, when the boundary has co-dimension $1$. We do not attempt to comprehensively review the literature in this area, but let us mention that considerable attention has been devoted to studying the geometric assumptions on co-dimension 1 boundaries for which the harmonic measure is absolutely continuous with respect to the surface measure \cite{dah1, dj, sem89, bl04, hm2}, as well as the converse so-called free boundary problems, where geometric information of the boundary is deduced from a priori solvability properties  \cite{azz1, hmu, ahmnt, ahmmmtv}, culminating in  the recent  results of \cite{ahmmt}, which  gives a complete picture of the relationship between absolute continuity of harmonic measure with respect to surface measure on the one hand, and uniform rectifiability plus a quantitative connectivity property on the other hand, for boundaries of co-dimension $1$.

However, when the boundary has co-dimension larger than $1$, the correspondence between  geometry and the theory of the Dirichlet problem for uniformly elliptic operators is severed, essentially owing to deep (and, as of yet, not completely understood) dimensional constraints on the support of the harmonic measure \cite{bj, bourgain, wolff}. Indeed, if we attempted to solve the Laplace equation outside of a boundary of high enough co-dimension, the equation does not ``see'' the boundary, and thus the uniformly elliptic theory is not a correct lens by which to characterize the geometry of such boundaries.

In \cite{dfm1}, G. David, J. Feneuil and the first author of this paper started a program to characterize the geometry of boundaries of high co-dimension via the theory of certain degenerate elliptic operators, crafted to overcome the fundamental myopia of uniformly elliptic operators vis-\`a-vis the geometry of low dimensional boundaries. Let $d\leq n-1$ be the Hausdorff dimension of the closed set $\Gamma\subset\bb R^n$. Furthermore, suppose that $\Gamma$ is quantifiably $d-$dimensional; more precisely, suppose that $\Gamma$ is $d-$Ahlfors-David regular (cf. Definition \ref{def.dadr}). Formally, consider
\begin{equation}\label{eq.opintro}
(Lu)(X)=-\dv\Big(\frac{\m A(X)}{\dist(X,\Gamma)^{n-d-1}}\nabla u(X)\Big),\qquad X\in\Omega:=\bb R^n\backslash\Gamma,
\end{equation}
where $\m A$ is a uniformly elliptic matrix. Observe that in compactly contained subsets of $\Omega$, the operator $L$ is strongly elliptic, but not uniformly so up to the boundary $\Gamma$, unless $d=n-1$. Instead, the operator must degenerate at a fixed rate which morally forces harmonic functions to respect the high co-dimension sets. It turns out that for these operators, one can recover an elliptic theory \cite{dfm1}, and in particular an analogue of harmonic measure can be devised. Degenerate operators have been considered in many previous works, for instance, \cite{fks} and \cite{fjk2}, but the solvability of boundary value problems was not studied.

Let us review a bit of the theory developed so far for the operators of the form (\ref{eq.opintro}). In \cite{dfm2} the  authors provided an analogue of Dahlberg's result \cite{dah1} which holds for their weighted elliptic operators. More precisely, for a $d-$dimensional Lipschitz graph $\Gamma$ with small Lipschitz constant, David, Feneuil and Mayboroda constructed a weighted elliptic operator of the form (\ref{eq.opintro}) so that the harmonic measure is absolutely continuous with respect to the surface measure on $\Gamma$. The equivalence of quantitative well-posedness of the Dirichlet problem and quantitative absolute continuity of harmonic measure was considered by Mayboroda and Zhao in \cite{mz}, so that from the two works \cite{dfm2}, \cite{mz}  we see the first solvability results of the Dirichlet problem for the operators defined in \cite{dfm1}. More recently, the Dirichlet problem $\operatorname{(D)}_p$ was tackled by Feneuil, Mayboroda, and Zhao in \cite{fmz} under some small Carleson norm assumptions on the coefficient matrix $\m A$, extending results of \cite{dp19} to this setting (see also \cite{dpp}). 

Regarding solvability of the Dirichlet problem on domains with uniformly rectifiable low dimensional boundaries, David and Mayboroda \cite{dm} have shown that for a suitable substitute of the Laplacian in the low dimensional setting, the harmonic measure is absolutely continuous with respect to the surface measure on $d-$dimensional uniformly rectifiable boundaries, with $d\leq n-2$ an integer (see also Feneuil \cite{fen1} for a different proof). On the other hand, in \cite{dem}, David, Engelstein, and Mayboroda manufactured an example which shows that for \emph{any} $d-$Ahlfors-David regular set $\Gamma$ with $d<n-2$,  there exists a special  operator $L_{\operatorname{DEM}}$ formally belonging to the class (\ref{eq.opintro}) whose harmonic measure is absolutely continuous with respect to the surface measure. The latter result lies in sharp contrast to the co-dimension $1$ case with the landmark free-boundary result \cite{ahmmt},  because it implies that for \dADR sets with $d<n-2$ an integer (recall that the complements of these sets always verify the interior Corkscrew  and Harnack Chain properties), uniform rectifiability alone cannot possibly characterize the $A_{\infty}$ property of the harmonic measure for all operators in the class (\ref{eq.opintro}). It is possible that the case of $L_{\operatorname{DEM}}$ is a miraculous arithmetic cancellation, and the free boundary result is still valid in some capacity. We show in this article, however,  that even if so, this arithmetic cancellation produces an entire family of counterexamples - see Corollary \ref{cor.magic} below.
	
We also mention briefly that an axiomatic elliptic theory for domains with mixed-dimensional boundaries is realized in \cite{dfm20},  while the Regularity problem  in the high co-dimension setting is dealt with in an  upcoming paper by Dai, Feneuil, and Mayboroda  \cite{daifm}. 

In this paper, we aim to further develop the theory of these degenerate elliptic operators by showing that quantifiable well-posedness of the Dirichlet problem is an open property.   Given an operator $L$ of the form (\ref{eq.opintro}), we say that the Dirichlet problem for $L$ with $L^p$ data is solvable, or $\operatorname{(D)}_p$ is solvable, if for each $f\in L^p(\Gamma,\sigma)$, there exists a unique solution $u$ to the problem
\[
\left\{\begin{matrix*}[l]Lu=0\qquad&&\text{in }\Omega,\\[2mm] u\lra f\qquad&&\text{non-tangentially},\\[2mm] \Vert Nu\Vert_{L^p(\Gamma,\sigma)}\lesssim\Vert f\Vert_{L^p(\Gamma,\sigma)}.\end{matrix*}\right.
\]
For the definition of the non-tangential maximal function $N$ and non-tangential convergence, see Definition \ref{def.ntmax}.  Moreover, as usual we understand $Lu=0$ in a weak sense, see Definition \ref{def.weaksol}. We are ready to state our main results.
 
\begin{theorem}[Carleson perturbation preserves Dirichlet problem solvability]\label{thm.main} Suppose that $\Gamma$ is \dADR with $d\in[1,n-1)$, $n\geq3$. Let two operators $L_0$ and $L$ be given as in (\ref{eq.opintro}) with associated bounded and uniformly elliptic  real (not necessarily symmetric) matrices $\m A_0$ and $\m A$  (see (\ref{eq.elliptic2}) for the definition of uniform ellipticity). Let $\omega_0,\omega$ denote the respective harmonic measures.  We define the disagreement of the matrices $\m A$, $\m A_0$ as
\begin{equation}\label{eq.disagreement}
\f a(X):=\sup\limits_{Y\in B(X,\delta(X)/2)}|\f E(Y)|,\quad\f E(Y):=\m A(Y)-\m A_0(Y),\qquad X\in\Omega.
\end{equation}
Assume that $\delta(X)^{d-n}\f a^2\,dX$ is a Carleson measure; that is, assume that there exists a constant $C\geq1$ such that for each surface ball $\Delta=B(x,r)\cap\Gamma$, the estimate
\begin{equation}\label{eq.acarlesonintro}
\dint_{T(\Delta)}\frac{\f a(X)^2}{\delta(X)^{n-d}}\,dX\leq C\sigma(\Delta) 
\end{equation}
holds, where $T(\Delta)=B(x,r)\cap\Omega$. Then, if $\operatorname{(D)}_{p'}$ is solvable for  $L_0$ for some $p'\in(1,\infty)$, then there exists $q'\in(1,\infty)$ such that $\operatorname{(D)}_{q'}$ is solvable for $L$. 
\end{theorem}

\begin{theorem}[Small Carleson perturbation preserves $\operatorname{(D)}_p$]\label{thm.main2} Suppose that $\Gamma$ is \dADR with $d\in[1,n-1)$, $n\geq3$.  Let two operators $L_0$ and $L$ be given as in (\ref{eq.opintro}) with associated uniformly elliptic real (not necessarily symmetric) matrices $\m A_0$ and $\m A$. Let $\omega_0,\omega$ denote the respective harmonic measures. Moreover, suppose that there exists $p'\in(1,\infty)$ such that $\operatorname{(D)}_{p'}$ is solvable for $L_0$. Define $\f a$ as in (\ref{eq.disagreement}), and assume that there exists $\ep_0>0$ so that
\begin{equation}\label{eq.acarlesonintro2}
\dint_{T(\Delta)}\frac{\f a(X)^2}{\delta(X)^{n-d}}\,dX\leq\ep_0\sigma(\Delta),\qquad\text{for all }\Delta\subset\Gamma.
\end{equation}
Then for all $\ep_0$ small enough, depending only on $n$, $d$, the \dADR constant of $\Gamma$, the ellipticity of $L_0$ and $L$, and the $RH_p$ characteristic of the Poisson kernel $\frac{d\omega_{L_0}}{d\sigma}$ (see Theorem \ref{thm.equivalence} and Definition \ref{def.ainftyw}), we have that $\operatorname{(D)}_{p'}$ is solvable for $L$ as well, and the $RH_p$ characteristic of $\frac{d\omega_{L_1}}{d\sigma}$ depends only on the same parameters as does $\ep_0$.
\end{theorem}

We note that when $\Gamma=\bb R^n\backslash\bb R^d$, a version of Theorem \ref{thm.main} with a  certain operator $L_{\operatorname{DFM}}$ satisfying some structural assumptions  is already part of the main result in \cite{dfm2}, where Carleson measure perturbations from their operator $L_{\operatorname{DFM}}$ are allowed. Likewise, a version of Theorem \ref{thm.main2} with a specific operator $L_{\operatorname{FMZ}}$ satisfying some structural assumptions is already in \cite{fmz} for $\Gamma=\bb R^n\backslash\bb R^d$. The novelty of our main results is that they hold for any  $d\in[1,n-1)$, any $d-$Ahlfors-David regular set $\Gamma$ (with $d$ not necessarily an integer), and  for any real operator $L_0$ verifying the well-posedness of $\operatorname{(D)}_p$, for some $p\in(1,\infty)$. We also briefly remark that Theorem \ref{thm.main} may not be deduced from Theorem \ref{thm.main2} due to the dependence of the latter on the $RH_p$ characteristic of   $\frac{d\omega_{L_0}}{d\sigma}$. 

Let us discuss a couple of immediate corollaries of our result. First, in the work  \cite{dm} (also \cite{fen1})   the authors obtain the following result. Suppose that $\Gamma\subset\bb R^n$  is a $d-$dimensional uniformly rectifiable set in $\bb R^n$ with $d\leq n-2$, $d\in\bb N$, that $\mu$ is a uniformly rectifiable measure on $\Gamma$, and that 
\[
L_{\mu,\alpha}=-\dv\Big(\frac1{D_{\mu,\alpha}^{n-d-1}}\nabla\Big),\qquad\text{in }\bb R^n\backslash\Gamma,
\]
where $D_{\mu,\alpha}$ is the regularized distance
\begin{equation}\label{eq.regdist}
D_{\mu,\alpha}(X):=\Big(\int_{\Gamma}|X-y|^{-d-\alpha}\,d\mu(y)\Big)^{1/\alpha},\qquad\alpha>0,~ X\in\bb R^n\backslash\Gamma.
\end{equation}
Then $\omega_{\mu,\alpha}\in A_{\infty}(\mu)$. We remark that in this setting, $\mu$ is equivalent to $\sigma=\n H^d|_{\Gamma}$, and for any $\alpha>0$, $D_{\mu,\alpha}\approx_{\alpha}\delta(X)$. Using our Theorem \ref{thm.main}, we are able to extend their class of solvable problems.

\begin{corollary}[Uniform rectifiability implies solvability of many operators] Let $\Gamma\subset\bb R^n$ be a closed $d-$dimensional uniformly rectifiable set with $d\leq n-2$ an integer. Suppose that $\mu$ and $L_{\mu,\alpha}$ are as described above, with $\alpha>0$. Let $L$ be an operator of the form (\ref{eq.opintro}) with matrix $\m A$ for which the disagreement of $\m A$ with $\m A_{\mu,\alpha}$ satisfies the Carleson measure perturbation (\ref{eq.acarlesonintro}). Then $\omega_L\in A_{\infty}(\n H^d|_{\Gamma})$.
\end{corollary}

Next, we recall the ``magic'' $\alpha$ example in \cite{dem}. Let $\Gamma$ be a (possibly purely unrectifiable)  closed unbounded $d-$Ahlfors-David regular set in $\bb R^n$ with $d\in(0,n-2)$ not necessarily an integer, and $n\geq3$. Let $\hat\alpha:=n-d-2$. Write $D_{\hat\alpha}=D_{\sigma,\hat\alpha}$ for the regularized distance (\ref{eq.regdist}) with $\mu=\sigma$. It turns out that in this setting, $D_{\hat\alpha}$ solves the equation $L_{\hat\alpha}u=L_{\sigma,\hat\alpha}u=0$ in $\bb R^n\backslash\Gamma$. Ultimately, this observation can be used to deduce that the harmonic measure of $L_{\hat\alpha}$ is equivalent to $\sigma$ (in the sense of pointwise equivalent bounds). Our Theorem \ref{thm.main} implies

\begin{corollary}[Open ball around ``magic $\alpha$'' example of \cite{dem}]\label{cor.magic} Let $n\geq3$ and let $\Gamma\subset\bb R^n$ be a closed \dADR set with $d\in[1,n-2)$. Suppose that $L$ is an operator of the form (\ref{eq.opintro}) with matrix $\m A$ for which the disagreement of $\m A$ with $\m A_{\sigma,\hat\alpha}$ satisfies the Carleson measure perturbation (\ref{eq.acarlesonintro}). Then $\omega_L\in A_{\infty}(\sigma)$.
\end{corollary}

In other words, the  counterexample of \cite{dem} to certain free boundary problems extends to give an open set of counterexamples.

Perturbation results are critical for the study of well-posedness, with applications to inverse problems and numerical analysis, because they allow some degree of relaxation of the assumptions on the coefficients needed for solvability. Indeed, necessary and sufficient conditions for solvability are very difficult to come by, whence a sensible strategy to establish a ``fat'' domain of operator invertibility is to first impose relatively strong conditions on the coefficients to make the problem tractable, and then in a second step exploit perturbation results. 

In the widely studied case of co-dimension $1$ with domain given by the half-space and transversal direction denoted by $t$, the robustness of well-posedness has generally been examined through two perspectives: $t-$independent $L^{\infty}$ perturbations, and Carleson measure perturbations. In a sense, these are complementary to one another, as the latter does not prescribe structural conformity but formally implies equality at the boundary of the perturbed coefficients to the original coefficients, while the former allows discrepancy at the boundary while imposing a rather strict structural condition (which is nevertheless natural; see \cite{cfk}). The $t-$independent sturdiness has been studied in \cite{fjk}, \cite{aah}, \cite{aamc}, \cite{aaahk} and more recently in \cite{bhlmp, bhlmp2} $t-$independent complex perturbations in the first-order and zeroth-order terms have been treated. In our setting, $t-$independent perturbation does not make sense since our operator must weigh the distance to the boundary, so that the structural advantage is lost unless $d=n-1$ which is the case that has classically been considered.

Therefore, robustness for our degenerate elliptic operators is better understood via Carleson measure perturbations. In the co-dimension $1$ setting, results in the spirit of Theorem \ref{thm.main} go back to \cite{fkp},\cite{kp},\cite{kp2} (see also \cite{hl}, \cite{hm1}, \cite{hmnote}), and have even been extended to complex coefficient perturbations under some conditions in \cite{aa}, \cite{hmm}. Analogues of the Carleson measure perturbation result Theorem \ref{thm.main} for the co-dimension 1 case have been seen to hold under progressively weaker geometric and topological assumptions on the boundary, going up to and including $1-$sided NTA domains with a capacity density condition; see \cite{mt10}, \cite{mpt}, \cite{chm}, \cite{chmt}, and \cite{ahmt}. That these perturbations are somehow optimal hypotheses to ensure transfer of the $A_{\infty}$ property of the harmonic measure (and hence, solvability of the Dirichlet problem on some $L^p$) was argued in \cite{fkp} for the co-dimension $1$ setting, where it was shown that, on the one hand, the analogous hypothesis of our Theorem \ref{thm.main} cannot provide a stronger result than $A_{\infty}\implies A_{\infty}$, and on the other hand, they show that no weaker assumption than the Carleson perturbation hypothesis can preserve $A_{\infty}$.  A key player in the study of the optimality of this hypothesis is the theory of quasiconformal mappings, which can be used to yield the existence of uniformly elliptic  operators on smooth domains for which the elliptic measure is not absolutely continuous with respect to the surface measure \cite{cfk} (for other such examples see \cite{mm, pog1} in the setting of co-dimension $1$, and \cite{dfmnew} for the case of low dimensional boundaries). In turn, Theorem \ref{thm.main2} is a higher co-dimensional analogue of \emph{small} Carleson perturbation theorems \cite{dah3}, \cite{mpt14}, \cite{chm} (see also \cite{esc}) which preserve solvability of $\operatorname{(D)}_p$ in the same $L^p$ space. 

Our method of proof follows the program of \cite{hm1} (see also \cite{chm}), where the main result of \cite{fkp} was given a new proof via an extrapolation theorem first presented by \cite{lm} and based on ideas of the Corona construction in \cite{carl}, \cite{cg}, where one aims to reduce matters to small perturbations on certain so-called sawtooth domains. Under this perspective, the overarching goal is to prove that the membership to a properly defined  $A_{\infty}$ class (see Definition \ref{def.ainftyw}) of the harmonic measure is preserved when the operator undergoes a perturbation ultimately controlled by its mass near the $d-$dimensional boundary. 

One of the main tools allowing one to attack uniformly rectifiable boundaries is the analysis on so-called sawtooth domains ``shielding'' bad parts of the boundary and providing a systematic comparison of our solutions with nice ones on a sawtooth domain. In contrast to the co-dimension $1$ case, our sawtooth domains will in general have mixed-dimensional boundaries, and therefore properties like Ahlfors-David regularity or uniform rectifiability cannot possibly transfer as-is to sawtooth domains. Nevertheless, we show that the sawtooth domains which we generate through a dyadic decomposition of the boundary \cite{christ, ds2} and which have been seen in \cite{hm2, mz} satisfy an axiomatic elliptic theory for sets with boundaries of mixed dimensions presented in \cite{dfm20}; in particular, we construct a measure on the boundary of the sawtooth domain which behaves sufficiently like a surface measure, as well as a suitable analogue of the harmonic measure on such sets. This allows us to show a global analogue of the Dahlberg-Jerison-Kenig sawtooth projection lemma \cite{djk}; in our setting, a similar result is shown in the upcoming work \cite{dm}.

Similarly, comparison principle techniques are much more subtle to use, due to the fact that the coefficients of our operators must see the boundary of the domain, and hence classical arguments in which restriction of the domain of the operator plays a crucial role are not available to us. Accordingly, our  arguments introduce  some new ideas even in the classical case of co-dimension $1$.

Furthermore,  we do not require either of the matrices to be symmetric in any of our main results. Though an analogue of Theorem \ref{thm.main} for the $1-$sided chord-arc domains was already shown in \cite{chmt} for the non-symmetric case, their methods are different and do not go through the extrapolation technique (besides, they do not show the non-symmetric case for their analogue of our Theorem \ref{thm.main2}). We will see in this article that the extrapolation technique does not need the symmetry of the matrices to work, essentially by being careful about the role of the adjoint operator. 

Finally, we remark on a couple more subtle technical differences in our approach to the proof of Theorem \ref{thm.main} than what has been seen in the literature. First, we circumvent the use of ``discrete'' tent spaces or the use of a dyadic averaging operator in our proof of Theorem \ref{thm.main} (though we have verified that both methods can work), in favor of a simple direct approach to exploiting the smallness of a ``discrete'' Carleson measure in a continuous setting (see (\ref{eq.calc3})).  Second, as mentioned above, we use a \emph{global} rather than local sawtooth projection lemma. This method allows us to only verify the axioms of the mixed-dimension elliptic theory for \emph{unbounded} sawtooth domains, and brings with it other simplifications in the geometric arguments; however, it introduces the complication that the globally constructed sawtooth domains do not locally coincide with the local sawtooth domains. Nevertheless, we shall see that this issue can be circumvented by realizing that the discrepancy between these sets is negligible from our point of view (see Lemma \ref{lm.poissoncube}).

In Section \ref{sec.bdrygeom}, we recap geometric results for our \dADR boundary. In Section \ref{sec.dyadic}, we present the dyadic decomposition of the boundary and related notions In Section \ref{sec.sawtooth}, we give a careful construction of the dyadically-generated sawtooth domains, reproving many results shown in the co-dimension $1$ setting. In Section \ref{sec.measure}, we construct a ``surface'' measure on the boundary of the sawtooth domain and prove that our sawtooth domains satisfy the axioms of the mixed-dimension elliptic theory of \cite{dfm20}. In Section \ref{sec.extrapolation}, we present the continuous Carleson measures, discrete Carleson measures, and the extrapolation theorem proved in \cite{dm}. In Section \ref{sec.theory}, we review the elliptic theory for sets of high co-dimension, considered in \cite{dfm1}. Section \ref{sec.djk} sees us proving the sawtooth projection lemma. Finally, in Sections \ref{sec.proof} and \ref{sec.smallthm}, we give the proofs of our main  results, Theorem \ref{thm.main} and Theorem \ref{thm.main2}.

\section{Geometry of domains with boundaries of high co-dimension}\label{sec.bdrygeom}

Throughout, our ambient space is $\bb R^n$, $n\geq3$. For $m\in\bb N$, we denote by $\n L^m$ the $m-$dimensional Lebesgue measure. For any $m\geq0$, we write $\n H^m$ for the $m-$dimensional Hausdorff measure (see \cite{federer}) . For integer $m$, we normalize $\n H^m$ so that it equals $\n L^m$. For $X\in\bb R^n$ and $r>0$, we write $B(X,r)\subset\bb R^n$ for the (open) ball of radius $r$ centered at $X$. If $A$ is a Borel set in $\bb R^n$ and $F\in L^1(A,\n L^n)$, we will often write $\int_AF\,\text{d}\n L^n(X)=\dint_AF\,dX$.  We will often write $a\lesssim b$ to mean that there exists a constant $C\geq1$ such that $a\leq Cb$, where $C$ may depend only on certain allowable parameters which will be identified with each statement. Likewise, we write $a\approx b$ if there exists a constant $C\geq1$ such that $\frac1Cb\leq a\leq Cb$. 

We now introduce the class of sets that the boundaries of our domains will reside in.

\begin{definition}[Ahlfors-David regular sets]\label{def.dadr} Let $\Gamma\subset\bb R^n$ be a closed set and $d\in(0,n)$. We say that $\Gamma$ is \emph{$d-$Ahlfors-David-regular} (or \dADR for short) if there exists a number $C_d\geq1$ such that for any $x\in\Gamma$ and $r>0$,
\begin{equation}\label{eq.ahlforsreg}
C_d^{-1}r^d\leq\n H^d(\Gamma\cap B(x,r))\leq C_dr^d,
\end{equation}
where $\n H^d$ is the $d-$dimensional Hausdorff measure. We shall often denote $\n H^d|_{\Gamma}$ by $\sigma$, and refer to it as the \emph{surface measure}. If $E\subset\bb R^n$ is a bounded, closed set, then we say that $E$ is \dADR if (\ref{eq.ahlforsreg}) holds for each $x\in E$ and $r\in(0,\diam E)$.
\end{definition}

Henceforth, we take $d\in(0,n-1)$ always. The set $\Gamma$ will always be a closed (unbounded) \dADR set, and $\Omega:=\bb R^n\backslash\Gamma$, so that $\Omega$ is an open set and $\partial\Omega=\Gamma$. Given $x\in\Gamma$ and $r>0$, we call 
\[
\Delta=\Delta(x,r)=\Gamma\cap B(x,r)
\]
a \emph{surface ball}. It is easy to see that, by virtue of (\ref{eq.ahlforsreg}), $\n H^d|_{\Gamma}$ is doubling on $\Gamma$ (see Definition \ref{def.doubling} below).  Thus, if $d|_{\Gamma}$ is the restriction of the Euclidean distance on $\bb R^n$ to $\Gamma$, then $(\Gamma,d|_{\Gamma})$ is a doubling metric space (that is, it admits a doubling measure), and hence the triple $(\Gamma,d|_{\Gamma},\sigma)$ is a space of homogeneous type (see \cite{christ}). It is obvious that if $\Gamma$ is a \dADR set, then so is any surface ball in $\Gamma$. Given a surface ball $\Delta=\Delta(x,r)$ and $c>0$, we denote by $c\Delta$ the set $\Gamma\cap B(x,cr)$.

We now see that the mass of a surface ball cannot be too concentrated at its center.

\begin{lemma}[Non-degeneracy of surface balls]\label{lm.surfdiam} Suppose that $\Gamma$ is a closed \dADR set. Then for any $x\in\Gamma$ and any $r>0$, 
\[
\diam\Delta(x,r)\geq\tfrac1{2^{1/d}}C_d^{-2/d}r.
\]
\end{lemma}

\noindent\emph{Proof.} We just need to prove that $\Delta(x,r)\backslash\Delta(x,C_d^{-2/d}r/2^{1/d})\neq\varnothing$, and this will be true provided that $\sigma\big(\Delta(x,r)\backslash\Delta(x,C_d^{-2/d}r/2^{1/d})\big)>0$. Observe that
\begin{multline*}
\sigma\big(\Delta(x,r)\backslash\Delta(x,C_d^{-2/d}r/2^{1/d})\big)=\sigma(\Delta(x,r))-\sigma(\Delta(x,C_d^{-2/d}r/2^{1/d}))\\ \geq C_d^{-1}r^d-C_d\big(C_d^{-2/d}r/2^{1/d}\big)^d=\big[1-\tfrac12\big]C_d^{-1}r^d>0,
\end{multline*}
as desired.\hfill{$\square$}

If $X\in\Omega$, then 
\[
\delta(X):=\text{dist}(X,\Gamma).
\]
We note in passing that since $\Gamma$ is closed, for each $X\in\Omega$ there exists $x\in\Gamma$ such that $|X-x|=\delta(X)$. It will be useful to denote 
\[
w(X):=\delta(X)^{-n+d+1},
\]
and let $m$ be the Borel measure on $\Omega$ given by $m(E):=\dint_Ew(X)\,dX$, $E\subseteq\Omega$. 

Next, we want to be able to use the openness of $\Omega$ in a quantitative way. The framework that we use is the following definition.

\begin{definition}[Corkscrew points]\label{def.corkscrew} Fix $x\in\Gamma$ and $r>0$. Then a point $X\in\Omega$ is called a \emph{Corkscrew point} (with \emph{Corkscrew constant} $c$) for the surface ball $\Delta(x,r)\subset\Gamma$ if there exists $c>0$ such that
\[
B(X,cr)\subset B(x,r)\cap\Omega.
\]
\end{definition}

For domains with a boundary of codimension less than or equal to $1$, it is not true in general that every surface ball has a Corkscrew point. The situation for \dADR  sets with $d\in(0,n-1)$ is different.

\begin{lemma}[Existence of Corkscrews; Lemma 11.6 of \cite{dfm1}]\label{lm.corkscrew} Suppose that $d<n-1$. Then there exists $c\in(0,1)$, depending only on $n$, $d$, and $C_d$, such that every surface ball in $\Gamma$ has a Corkscrew point with Corkscrew constant $c$.
\end{lemma}

Furthermore, domains with high co-dimensional \dADR boundaries enjoy quantitative connectedness as well, as given in

\begin{lemma}[Existence of Harnack Chains; Lemma 2.1 of \cite{dfm1}]\label{lm.harnackchain} Suppose that $d<n-1$. Then there exists a constant $c_{\m H}\in(0,1)$, that depends only on $d, n, C_d$, such that for $\Lambda\geq1$ and $X_1,X_2\in\Omega$ with $\delta(X_i)\geq s$ and $|X_1-X_2|\leq\Lambda s$, we can find two points $Y_i\in B(X_i, s/2)$ verifying that 
\[
\text{dist}([Y_1,Y_2],\Gamma)\geq (c_{\m H}\Lambda^{-\frac d{n-1-d}})s,
\]
where $[Y_1,Y_2]$ is the straight line segment in $\bb R^n$ with endpoints $Y_1$ and $Y_2$. That is, there is a thick tube in $\Omega$ that connects the balls $B(X_i,s/2)$.
\end{lemma}

Inspired by the prior lemma, we consider 

\begin{definition}[Harnack Chains for $\bb R^n\backslash\Gamma$]\label{def.hc} Let $c_{\m H}, s,\Lambda, X_1, X_2, Y_1$, and $Y_2$ be as described in Lemma \ref{lm.harnackchain}. Let $\{B_j\}_j$ be a finite collection of balls whose centers lie in the line segment $[Y_1, Y_2]$, whose radii are less than $\frac12c_{\m H}\Lambda^{-\frac d{n-1-d}}s$, and so that they cover $[Y_1,Y_2]$.   If the radii of the balls $B_j$ are all exactly $\frac12c_{\m H}\Lambda^{-\frac d{n-1-d}}s$, we call $\m H:=\cup_jB_j$ a \emph{(well-tempered) Harnack Chain connecting $X_1$ and $X_2$}.
\end{definition}

\begin{remark}\label{rm.closetohm} As was shown in \cite{dfm1} and mentioned above, the closed unbounded \dADR set $\Gamma$ that we consider has ample interior Corkscrew points and Harnack Chains. Therefore, our boundary $\Gamma$ is \emph{axiomatically} very similar to the boundaries of the so-called \emph{$1$-sided chord-arc domains}, which are open sets $\tilde\Omega\subset\bb R^n$ whose boundary $\partial\tilde\Omega$ is $n$-Ahlfors-David regular, and having interior Corkscrews and Harnack Chains. For this reason, as we shall see in the rest of the paper, many of the results here have direct analogues in the setting of $1$-sided chord-arc domains that was explored in the seminal paper \cite{hm2}, and often with very similar proofs, that we decide to omit in some cases, and in some other occasions, we decide to give different proofs of the expected results.

A key difference for us is that the sawtooth domains (to be defined further below) will have boundaries of mixed dimension, whence the usual global ADR notion is meaningless for them. Another issue is that we cannot rely on comparison principles for domains with different boundaries than $\Gamma$, because the coefficients of the operator explicitly depend on the distance to the boundary; instead we are restricted to work with ``global'' comparison principles. Still, we are able to overcome these issues. Hence we supply very careful and detailed proofs of our results leading up to the analogue of the Dahlberg-Jerison-Kenig sawtooth-projection lemma.
\end{remark}

We will need to study non-negative doubling Borel measures on $\Gamma$.

\begin{definition}[Doubling measures]\label{def.doubling} Fix a surface ball $\Delta(x,r)\subseteq\Gamma$, $r\in(0,+\infty]$ (by convention, $\Delta(x,+\infty)=\Gamma$).  We say that a non-negative Borel measure $\mu$ on $\Delta$ is \emph{doubling} on $\Delta$ if there exists a constant $M$ large enough such that for each surface ball $\Delta'$ with $2\Delta'\subset\Delta$, we have that
	\begin{equation}\label{eq.doublingcond}\notag
	\mu(2\Delta')\leq C_{\doubling}\mu(\Delta').
	\end{equation}
\end{definition}
 
The following definition gives a quantitative version of absolute continuity between measures.
 
\begin{definition}[$A_{\infty}$ measures]\label{def.ainfty} Given a doubling non-negative Borel measure $\nu$ on $\Gamma$, and a fixed surface ball $\Delta\subseteq\Gamma$, we say that the doubling measure $\mu$ is of class $A_{\infty}(\nu,\Delta)$ if for each $\ep>0$, there exists a number $\xi=\xi(\ep)>0$ such that for  every surface ball $\Delta'\subseteq\Delta$, and every Borel set $E\subset\Delta'$, we have that
\begin{equation}\label{eq.ainftycond}
\frac{\nu(E)}{\nu(\Delta')}<\xi\implies\frac{\mu(E)}{\mu(\Delta')}<\ep.
\end{equation}
\end{definition}

After reviewing the elliptic theory for the sets we are studying, we will need to study a more precise quantification of absolute continuity than just membership to $A_{\infty}$. Given a doubling Borel measure $\mu$ on $\Gamma$, a \emph{weight} $\f w$ on $\Gamma$ is a non-negative $L_{\textup{loc}}^1(\Gamma,\mu)$ function. A weight induces a non-negative Borel measure as follows: for any $\mu-$measurable set $E\subset\Gamma$ we write $\f w(E):=\int_E\f w\,d\mu$. 

\begin{definition}[The Reverse H\"older class $RH_p$]\label{def.rh} Given a non-negative doubling Borel measure $\mu$ on $\Gamma$, a fixed surface ball $\Delta_0\subset\Gamma$,  a weight $\f w\in L_{\textup{loc}}^1(\Delta_0,\mu)$, and $p\in(1,\infty)$, we say that $\f w\in RH_p(\mu,\Delta_0)$ if there exists a constant $C_p$ such that for every surface ball $\Delta\subset\Delta_0$,
\begin{equation}\label{eq.rhpcond}
\Big(\frac1{\mu(\Delta)}\int_\Delta\f w^p\,d\mu\Big)^{\frac1p}\leq C_p\frac1{\mu(\Delta)}\int_\Delta\f w\,d\mu.
\end{equation}
We denote by the \emph{$RH_p(\mu,\Delta_0)$ characteristic of }$\f w$ the smallest number $C_p$ such that (\ref{eq.rhpcond}) holds for all $\Delta\subset\Delta_0$. When $\mu$ is the surface measure $\sigma$, we simply write $RH_p(\sigma,\Delta_0)=RH_p(\Delta_0)$. Furthermore, if $\nu,\mu$ are two doubling non-negative Borel measures on $\Gamma$, $\Delta_0\subset\Gamma$, and $p\in(1,\infty)$, we say that $\mu\in RH_p(\nu,\Delta_0)$ if $\mu\ll\nu$ and the Radon-Nikodym derivative $\frac{d\mu}{d\nu}$ lies in $RH_p(\nu,\Delta)$. 
\end{definition}

Let us record some properties of the $A_{\infty}$   class. For our setting, the following results have appeared in \cite{jawerth} and \cite{stweighted}.

\begin{theorem}[Properties of $A_{\infty}$ measures; Theorem 1.4.13 of \cite{kbook}; \cite{stweighted}]\label{thm.ainftyprop}  Let $\mu,\nu$ be doubling non-negative Borel measures on $\Gamma$, and let $\Delta$ be a surface ball. The following statements hold.
\begin{enumerate}[(i)]
	\item If $\mu\in A_{\infty}(\nu,\Delta)$, then $\mu$ is absolutely continuous with respect to $\nu$ on $\Delta$.
	\item The class $A_{\infty}$ is an equivalence relation.
	\item We have that $\mu\in A_{\infty}(\nu,\Delta)$ if and only if there exist uniform constants $C>0,\theta>0$, such that for each surface ball $\Delta'\subseteq\Delta$ and each Borel set $E\subseteq\Delta'$, we have that
	\[
	\frac{\mu(E)}{\mu(\Delta')}\lesssim\Big(\frac{\nu(E)}{\nu(\Delta')}\Big)^{\theta}.
	\]
	\item We have that $\mu\in A_{\infty}(\nu,\Delta)$ if and only if there exist $\ep\in(0,1),\delta\in(0,1)$ so that (\ref{eq.ainftycond}) holds for all surface balls $\Delta'\subset\Delta$ and all Borel sets $E$ (see \cite{gr}).
	\item The characterization $A_{\infty}(\nu)=\bigcup_{p>1}RH_p(\nu)$ holds.
	\item We have that $\mu\in RH_p(\nu)$ if and only if the uncentered Hardy-Littlewood maximal function adapted to $\Gamma$,
	\[
	(M_{\mu}f)(x):=\sup_{\Delta\ni x}\frac1{\mu(\Delta)}\int_\Delta|f|\,d\mu,
	\]
	verifies the estimate
	\[
	\Vert M_{\mu}f\Vert_{L^{p'}(d\nu)}\lesssim\Vert f\Vert_{L^{p'}(d\nu)},
	\]
	where $p'$ is the H\"older conjugate of $p$, so that $\frac1p+\frac1{p'}=1$.
\end{enumerate}
\end{theorem}
%

To conclude this minimal set-up for \dADR sets, we give a meaning to the non-tangential maximal functions and square functions, which are essential concepts in theory of the Dirichlet problems with rough data.

\begin{definition}[Non-tangential maximal function and square function]\label{def.ntmax} For any $x\in\Gamma$ and $\alpha>0$, we define the \emph{non-tangential cone} $\gamma^{\alpha}(x)$ with vertex $x$ and aperture $\alpha$ as
\begin{equation}\label{eq.ntcone}\notag
\gamma^{\alpha}(x)=\Big\{X\in\Omega\,:\,|X-x|<(1+\alpha)\delta(X)\Big\},
\end{equation}
We often omit the superscript $\alpha$. Define the \emph{square function} as
\begin{equation}\label{eq.squarefn}\notag
Su(x)=\Big(\dint_{\gamma(x)}|\nabla u(X)|^2\delta(X)^{2-n}\,dX\Big)^{\frac12}.
\end{equation}
Finally, the \emph{non-tangential maximal function}  is given by
\begin{equation}\label{eq.ntmaximalfn}\notag
Nu(x)=\sup\limits_{X\in\gamma(x)}|u(X)|.
\end{equation}
Given a measurable function $f$ on $\Gamma$, we say that $u\ra f$ \emph{non-tangentially} if for $\sigma-$almost every $x\in\Gamma$, we have that $\lim\limits_{\gamma(x)\ni X\ra x}u(X)=f(x)$.
\end{definition}

\section{Dyadic decomposition of sets of high co-dimension}\label{sec.dyadic}

In the following lemma, we exhibit a family of partitions for $\Gamma$ which are analogous to dyadic cubes. The original construction of such a dyadic grid for \dADR sets with $d=n-1$ is found in  \cite{davidibero}; in the book \cite{davidwavelets} there is a simpler proof which adapts to our setting. See also \cite{christ} for a different proof in the even more general case of spaces of homogeneous types.

\begin{lemma}[Dyadic cubes for $d-$Ahlfors-David regular set]\label{lm.dyadiccubes} \cite{davidwavelets}.  There exist constants $a_0\in(0,1],A_0\in[1,\infty),\zeta\in(0,1)$, depending only on $d$, $n$, and the \dADR constant $C_d$, such that for each $k\in\bb Z$, there is a collection of Borel sets \emph{(``dyadic cubes'')}
	\[
	\bb D^k=\bb D^k(\Gamma):=\{Q_j^k\subset\Gamma\,:\,j\in\n J^k\},
	\]
	where $\n J^k$ denotes some indexing set depending on $k$, satisfying the following properties.
	\begin{enumerate}[(i)]
		\item\label{item.grid} For each $k\in\bb Z$, $\Gamma=\bigcup_{j\in\n J^k}Q_j^k$.
		\item\label{item.nonoverlap} If $m\geq k$ then either $Q_i^m\subset Q_j^k$ or $Q_i^m\cap Q_j^k=\varnothing$.
		\item\label{item.generations} For each pair $(j,k)$ and each $m<k$, there is a unique $i\in\n J^m$ such that $Q_j^k\subset Q_i^m$. When $m=k-1$, we call $Q_i^m$ the \emph{dyadic parent} of $Q_j^k$, and we say that $Q_j^k$ is a \emph{dyadic child} of $Q_i^m$.
		\item\label{item.diam} diam $Q_j^k<A_02^{-k}$.
		\item\label{item.inscribe} Each $Q_j^k$ contains some surface ball $\Delta(x_j^k,a_02^{-k})=B(x_j^k,a_02^{-k})\cap\Gamma$.
		\item\label{item.bdrythin} $\n H^d(\{x\in Q_j^k:\text{dist}(x,\Gamma\backslash Q_j^k)\leq\rho2^{-k}\})\leq A_0\rho^{\zeta}\n H^d(Q_j^k)$, for all $(j,k)$ and all $\rho\in(0,a_0)$.
	\end{enumerate}
\end{lemma}

Let us define some notions and state some useful properties of this construction.

$\bullet$ We shall denote by $\bb D=\bb D(\Gamma)$ the collection of all relevant $Q_j^k$; that is,
\[
\bb D=\bb D(\Gamma):=\bigcup_{k\in\bb Z}\bb D^k(\Gamma).
\]
Henceforth, we refer to the elements of $\bb D$ as \emph{dyadic cubes}, or \emph{cubes}. For $Q\in\bb D$, we write $\bb D_Q:=\{Q'\in\bb D\,:\,Q'\subseteq Q\}$, and $\bb D_Q^k=\bb D^k(\Gamma)\cap\bb D_Q$.

$\bullet$ Note carefully that if $Q_i^{k+1}$ is the dyadic parent of $Q_j^k$, then it is possible that, \emph{as sets}, $Q_i^{k+1}=Q_j^k$. In other words, if $Q\in\bb D$, then the set 
\[
\bb K(Q):=\{k\in\bb Z:Q=Q_j^k\text{ for some }j\}
\] 
may in general have cardinality greater than or equal to 1. We call $\bb K(Q)$ the \emph{generational bandwith} of $Q$. By Lemma \ref{lm.surfdiam} and properties \ref{item.diam} and \ref{item.inscribe} above, we have that if $k\in\bb K(Q)$, then
\begin{equation}\label{eq.qdiam}
2^{-1/d}C_d^{-2/d}a_02^{-k}\leq\diam Q\leq A_02^{-k},
\end{equation}
which implies that $\bb K(Q)$ is finite, and in fact,
\begin{equation}\label{eq.genband}
1\leq\card(\bb K(Q))\leq\log_2\big[2^{\frac{d+1}d}C_d^{2/d}A_0a_0^{-1}\big].
\end{equation}
Define the \emph{dyadic generation} of $Q\in\bb D$ as the oldest generation that $Q$ belongs to; that is,
\[
k(Q)=\min_{k\in\bb K(Q)}k,
\]
and note that, if $k\in\bb K(Q)$, then
\[
k(Q)\leq k\leq k(Q)+\log_2\big[2^{\frac{d+1}d}C_d^{2/d}A_0a_0^{-1}\big].
\]
We  call the number $\ell(Q)=2^{-k(Q)}$ the \emph{length} of $Q$. Given a fixed $Q\in\bb D$, we call a cube $Q'\in\bb D_Q\backslash\{Q\}$ a \emph{proper child} of $Q$ if $\ell(Q')<\ell(Q)$ and $\ell(Q')\geq\ell(Q'')$ for any other $Q''\in\bb D_Q$. Likewise, given $Q\in\bb D$, we call a cube $\hat Q\in\bb D$ with $\hat Q\supset Q$ a \emph{proper parent} of $Q$ if $\ell(\hat Q)>\ell(Q)$ and $\hat Q\subseteq\tilde Q$ for any $\tilde Q\in\bb D$ with $\tilde Q\supsetneq Q$. If $Q'\in\bb D_Q$ is a proper child of $Q$, then we have that
\begin{equation}\label{eq.largechild}
\ell(Q)>\ell(Q')\geq\frac{a_0}{2^{\frac{2d+1}d}C_d^{2/d}A_0}\ell(Q)=:c_{\bb K}\ell(Q).
\end{equation}
If $Q'$ is a proper child of $Q$, then by the partitioning property of the dyadic cubes we must have that there exists a collection $\{Q''\}$ of proper children of $Q$ such that $\cup Q''=Q$. In the sequel, if we say that $Q'$ is a \emph{child} of $Q$, we mean that $Q'$ is a dyadic child of $Q$, leaving open the possibility that $Q'=Q$ as sets.

$\bullet$ (Almost inscription and subscription of surface balls). Properties \ref{item.diam} and \ref{item.inscribe} also imply that for each cube $Q\in\bb D $, there is a point $x_Q\in\Gamma$ such that  
\begin{equation}\label{eq.centerofQ} 
\Delta(x_Q,a_0\ell(Q))\subset Q\subseteq \Delta(x_Q,A_0\ell(Q)).
\end{equation}
We call $x_Q$ the \emph{center} of $Q$. We note that (\ref{eq.centerofQ}) and (\ref{eq.ahlforsreg}) imply the following estimate on the surface measure of $Q$:
\begin{equation}\label{eq.measureofQ}
C_d^{-1}a_0^{d}\ell(Q)^d\leq\sigma(Q)\leq C_dA_0^d\ell(Q)^d.
\end{equation}

$\bullet$ (Number of children of $Q$). Fix $Q\in\bb D$ and let $\{Q_j\}_{j\in J}$ be the collection of all (dyadic) children of $Q$. It must be the case by property \ref{item.grid} that $Q=\cup_{j\in J}Q_j$. Observe the elementary estimate 
\begin{equation*}
\sigma(Q)=\sigma(\cup_{j\in J}Q_j)=\sum_{j\in J}\sigma(Q_j)\geq C_d^{-1}a_0^d\frac{\ell(Q)}2\card J, 
\end{equation*}
where we used (\ref{eq.measureofQ}) in the last inequality. Putting the previous estimate together with the upper bound in (\ref{eq.measureofQ}) gives that
\begin{equation}\label{eq.cardchildren}
\card(\{\text{children of }Q\})\leq C_d^2[a_0^{-1}A_02]^d.
\end{equation}

$\bullet$ (Corkscrew points for $Q$). We denote by $X_Q$ a point in $\Omega$ which is a Corkscrew point (with Corkscrew constant $\tilde c>0$) for the surface ball $\Delta(x_Q,a_0\ell(Q))$. Such a point is called a \emph{Corkscrew point for $Q$} (with Corkscrew constant $\tilde c>0$).

$\bullet$ The inequality in \ref{item.bdrythin} says that the boundary of a dyadic cube $Q_j^k$ is uniformly thin; indeed, one may easily deduce from it that $\n H^d(\partial Q)=0$ for any $Q\in\bb D$.

$\bullet$ (Dyadic maximal function). Given $Q_0\in\bb D$ and a function $f\in L^1_{\loc}(Q_0,\sigma)$, we write
\[
E^kf(x)=\sum\limits_{Q\in\bb D^k_{Q_0}}\Big(\frac1{\sigma(Q)}\int_Qf\,d\sigma\Big){\bf 1}_Q(x),
\]
and we define the \emph{dyadic Hardy-Littlewood maximal function} with respect to $Q_0$ as
\begin{equation}\label{eq.dyadichl}\notag
M_{Q_0}^{\dyadic}f(x)=\sup_k|E^kf(x)|.
\end{equation}

$\bullet$ By $\m F=\{Q_j\}_j$ we denote a family of pairwise disjoint dyadic cubes in $\bb D$, which we identify \emph{as subsets of $\Gamma$} and not as elements of $\cup_k\bb D^k$. Accordingly, if $Q_j\in\m F$, then its parent $\hat Q\in\bb D$ does not belong to $\m F$, and we have that $\ell(Q_j)<\ell(\hat Q)\leq c_{\bb K}^{-1}\ell(Q_j)$. We refer to such a collection $\m F$ as a \emph{disjoint family}.

$\bullet$ We define the projection operator $\m P_{\m F}:L^1_{\loc}(\Gamma,\sigma)\ra L^1_{\loc}(\Gamma,\sigma)$ by
\begin{equation}\label{eq.projection}
(\m P_{\m F}f)(x):=f(x){\bf 1}_{\Gamma\backslash(\cup_jQ_j)}(x)+\sum\limits_j\Big(\frac1{\sigma(Q_j)}\int_{Q_j}f\,d\sigma\Big){\bf 1}_{Q_j}(x),\qquad x\in\Gamma.
\end{equation}
One has that $\m P_{\m F}\circ\m P_{\m F}=\m P_{\m F}$, $\m P_{\m F}$ is self-adjoint, and $\Vert\m P_{\m F}f\Vert_{L^p(\Gamma,\sigma)}\leq\Vert f\Vert_{L^p(\Gamma,\sigma)}$ for every $p\in[1,\infty]$. Observe that if $\mu$ is a non-negative finite Borel measure on $\Gamma$ and $E\subseteq\Gamma$ is a Borel set, then we may naturally define the measure $\m P_{\m F}\mu$ as follows:
\[
\m P_{\m F}\mu(E):=\int_{\Gamma}\m P_{\m F}({\bf 1}_E)\,d\mu=\mu(E\backslash\cup_jQ_j)+\sum\limits_j\frac{\sigma(E\cap Q_j)}{\sigma(Q_j)}\mu(Q_j).
\]
In particular, $\m P_{\m F}\mu(\Gamma)=\mu(\Gamma)$. Notice that $\m P_{\m F}\mu$ is defined in such a way that it coincides with $\mu$ in $\Gamma\backslash(\cup_jQ_j)$ and in each $Q_j$ we replace $\mu$ by $\mu(Q_j)/\sigma(Q_j)\,d\sigma$.

\subsection{The theory of quantitative absolute continuity adapted to the dyadic grid} As an initial step, we will need to make sense of the doubling property adapted to our dyadic grids.

\begin{definition}[Dyadically doubling measures]\label{def.dyadicdoubling} We say that a Borel measure $\mu$ on $Q_0\in\bb D$ is \emph{dyadically doubling} in $Q_0$ if $0<\mu(Q)<\infty$ for every $Q\in\bb D_{Q_0}$ and there exists a constant $C_{\mu}\geq1$ such that $\mu(Q)\leq C_{\mu}\mu(Q')<\infty$ for every $Q\in\bb D_{Q_0}$ and for every dyadic child $Q'$ of $Q$.
\end{definition}

Let us relate the concepts of doubling and dyadically doubling measures.

\begin{lemma}[Doubling implies dyadically doubling]\label{lm.dimpliesdd} Fix $Q_0\in\bb D$ and suppose that $\mu$ is a doubling Borel measure on the surface ball $\Delta_0:=\Delta(x_{Q_0},2A_0\ell(Q_0))$. Then $\mu$ is dyadically doubling in $Q_0$.
\end{lemma}

\noindent\emph{Proof.} Fix $Q\in\bb D_{Q_0}$. Since $Q\supset\Delta_Q=\Delta(x_Q,a_0\ell(Q))$ and $\mu$ is doubling in $\Delta_0\supset\Delta_Q$, then $\mu(Q)\in(0,+\infty)$. Now fix $Q'\in\bb D_Q$ to be any dyadic child of $Q$. Let $y\in Q$, and observe that
\[
|y-x_{Q'}|\leq\diam Q\leq A_0\ell(Q)\leq2A_0\ell(Q'),
\]
so that $Q\subset\Delta(x_{Q'},2A_0\ell(Q'))$. It follows that
\[
\mu(Q)\leq\mu\big(\Delta(x_{Q'},2A_0\ell(Q'))\big)\lesssim\mu\big(\Delta(x_{Q'},a_0\ell(Q'))\big)\leq\mu(Q'),
\]
as desired.\hfill{$\square$}

The following lemma gives us that if a measure is dyadically doubling, then so is its projection. We omit its easy proof as it is essentially the same as that in \cite{hm2} (see Remark \ref{rm.closetohm}).

\begin{lemma}[Lemma B.1 of \cite{hm2}]\label{lm.pisdoubling} Fix $Q_0\in\bb D$, let $\m F\subset Q_0$ be a disjoint family, and let $\mu$ be a dyadically doubling measure in $Q_0$. Then   $\m P_{\m F}\mu$ is dyadically doubling in $Q_0$.
\end{lemma}

We now record the fact that we have a local Calder\'on-Zygmund decomposition for dyadically doubling weights. Its proof is standard.

\begin{lemma}[Lemma B.12 of \cite{hm2}]\label{lm.czdecomp} Given $Q_0\in\bb D$ and $\mu$ a dyadically doubling measure on $Q_0$ with doubling constant $C_{\mu}$, we consider the local dyadic Hardy-Littlewood maximal function with respect to $\mu$:
	\[
	\m M_{\mu}f(x)=\sup_{x\in Q\in\bb D_{Q_0}}\frac1{\mu(Q)}\int_Q|f|\,d\mu.
	\]
	For any $0\leq f\in L^1(Q_0,\mu)$ and $\tau\geq\frac1{\mu(Q_0)}\int_{Q_0}|f|\,d\mu$, there exists a collection of maximal and therefore disjoint dyadic cubes $\{Q_j\}_j\subset\bb D_{Q_0}$ such that
	\begin{gather*}
	\Omega_{\tau}=\{x\in Q_0\,:\,\m M_{\mu}f(x)>\tau\}=\bigcup_jQ_j,\label{eq.cz1}\\ f(x)\leq\tau,\quad\text{for }\mu-\text{a.e. }x\notin\Omega_{\tau}\label{eq.cz2}\\ \tau<\frac1{\mu(Q_j)}\int_{Q_j}f(y)\,d\mu(y)\leq C_{\mu}\tau.\label{eq.cz3}
	\end{gather*}
\end{lemma}

Now we define quantitative absolute continuity on our dyadic grid.

\begin{definition}[$A_{\infty}^{\dyadic}$ and $RH_p^{\dyadic}$]\label{def.dyadica} Given $Q_0\in\bb D$, we say that a Borel measure $\mu$ defined on $Q_0$ belongs to $A_{\infty}^{\dyadic}(Q_0)$ if there exist constants $0<\alpha,\beta<1$ such that for every $Q\in\bb D_{Q_0}$ and for every Borel set $F\subset Q$, we have that
\[
\frac{\sigma(F)}{\sigma(Q)}>\alpha\quad\implies\quad\frac{\mu(F)}{\mu(Q)}>\beta.
\]
Given $p\in(1,\infty)$, we say that $\mu\in RH_p^{\dyadic}(Q_0)$ if and only if $\mu\ll\sigma$ in $Q_0$ and there exists a constant $C_p\geq1$ such that for every $Q\in\bb D_{Q_0}$, the estimate
\[
\Big(\frac1{\sigma(Q)}\int_Qk^p\,d\sigma\Big)^{\frac1p}\leq C_p\frac1{\sigma(Q)}\int_Qk\,d\sigma
\]
holds, where $k=\frac{d\mu}{d\sigma}$.
\end{definition}

The next result gives that the $A_{\infty}^\dyadic$ property is passed on from a measure to its projection. Its proof is essentially the same as that of Lemma 4.1 in \cite{hmnote}  (see Remark \ref{rm.closetohm}); and so we omit the details. 

\begin{lemma}[Lemma 4.1 of \cite{hmnote}]\label{lm.projdyadic} Fix $Q_0\in\bb D$, let $\m F\subset Q_0$ be a disjoint family, and suppose that $\mu\in A_{\infty}^{\dyadic}(Q_0)$. Then $\m P_{\m F}\mu\in A_{\infty}^{\dyadic}(Q_0)$.
\end{lemma}

As expected, we have the symmetry of the $A_{\infty}^{\dyadic}$ class. It is proven in \cite{hm2} in a similar setting  (see Remark \ref{rm.closetohm}), but their proof generalizes to our situation immediately.

\begin{lemma}[Symmetry of $A_{\infty}^{\dyadic}$, Lemma B.7 of \cite{hm2}]\label{lm.ainftrans} Let $Q_0\in\bb D$  and let $\mu$, $\nu$ be two dyadically doubling measures on $Q_0$. Assume that there exist positive constants $C_0,\theta_0,$ such that for all $Q\in\bb D_{Q_0}$ and all Borel sets $F\subseteq Q$,
	\begin{equation}\label{eq.ainftrans1}\notag
	\frac{\nu(F)}{\nu(Q)}\leq C_0\Big(\frac{\mu(F)}{\mu(Q)}\Big)^{\theta_0}.
	\end{equation}
	Then, there exist positive constants $C_1,\theta_1$ such that for all $Q\in\bb D_{Q_0}$ and all Borel sets $F\subseteq Q$,
	\begin{equation}\label{eq.ainftrans2}\notag
	\frac{\mu(F)}{\mu(Q)}\leq C_1\Big(\frac{\nu(F)}{\nu(Q)}\Big)^{\theta_1}.
	\end{equation}
\end{lemma}

To finalize this section, we present a generalization of Lemma B.7 in \cite{hm1}, which allows us to conclude that a measure is $A_{\infty}^{\dyadic}$ if it satisfies a certain local Reverse H\"older inequality. The result has been shown in a setting very similar to ours in Lemma 3.1 of \cite{chm}  (see Remark \ref{rm.closetohm}); the proof is essentially the same and so we omit it.

\begin{lemma}[Local $RH$ implies $A_{\infty}^{\dyadic}$, Lemma 3.1 of \cite{chm}]\label{lm.technical} Fix $Q_0\in\bb D$ and $\ep\in(0,1)$. Let $v\in L^1(Q_0)$ be a function such that there exists $C_0\geq1$ verifying that $0<v(Q)\leq C_0v(\ep\Delta_Q)$ for every $Q\in\bb D_{Q_0}$. Assume also that there exist $C_1\geq1$ and $p\in(1,\infty)$ such that
	\begin{equation}\label{eq.vrhp}\notag
	\Big(\frac1{\sigma(\ep\Delta_Q)}\int_{\ep\Delta_Q}v^p\,d\sigma\Big)^{\frac1p}\leq C_1\frac1{\sigma(\ep\Delta_Q)}\int_{\ep\Delta_Q}v\,d\sigma,\qquad\text{for each }Q\in\bb D_{Q_0}.
	\end{equation}
Then $v\in RH_p^{\dyadic}(Q_0)$ with uniform constants depending only on $n$, $d$, $C_d$, $p$, $\ep$, $C_0$, and $C_1$. 
\end{lemma}

\section{Dyadically-generated sawtooth domains}\label{sec.sawtooth} 

\subsection{Construction of dyadically-generated sawtooth domains} In this subsection, we will construct the so-called sawtooth domains for \dADR sets with $d<n$. The abstract construction here was first considered for certain domains of co-dimension $1$ in \cite{hm2}, and developed for the setting $d<n-1$ in \cite{mz}. If $d<n-1$, we have no further assumptions, since the $\dADR$ property of $\Gamma$ gives the existence of Corkscrew points and Harnack Chains.  We will see that our dyadically-generated sawtooth domains will necessarily be mixed-dimensional, so that an elliptic theory for them is highly non-trivial. As such, we will give a  careful construction with the goal to prove in the following section that our dyadic sawtooth domains satisfy the mixed-dimension theory of \cite{dfm20}. Indeed, while in \cite{dfm20} it is shown that certain sawtooth domains over Lipschitz graphs satisfy their axioms, our dyadic sawtooth domains over arbitrary \dADR sets (with possibly fractional dimension) were not considered, and the verification is considerably more subtle.

Since $\Omega=\bb R^n\backslash\Gamma$ is an open set in $\bb R^n$, there exists a collection of closed dyadic \emph{Whitney boxes}, denoted by $\m W=\m W(\Omega)$, so that the interiors of the boxes never overlap pairwise, the boxes form a covering of $\Omega$, and moreover they satisfy the conditions
\begin{equation}\label{eq.whitney1}
4\text{diam }I\leq\text{dist}(4I,\Gamma)\leq\text{dist}(I,\Gamma)\leq40\text{diam }I,\qquad\text{for each }I\in\m W,
\end{equation}
and
\begin{equation}\label{eq.whitney2}\notag
\frac14\text{diam }I_1\leq\text{diam }I_2\leq4\text{diam }I_1
\end{equation}
whenever $\partial I_1\cap\partial I_2\neq\varnothing$ (see, for instance, \cite{ste2}). Let $X_I$ denote the \emph{center} of $I$ and $\ell(I)$ the \emph{side-length} of $I$, so that $\text{diam }I=\sqrt n\ell(I)$. We also write $k(I)=k$ if $\ell(I)=2^{-k}$. We say that two Whitney boxes \emph{touch}, or that they are \emph{adjacent}, if their boundaries intersect. If $X\in I$ and $I\in\m W$, then
\[
4\diam I\leq\delta(X)\leq41\diam I.
\]

Next, we want to associate to each $Q\in\bb D$ a ``Whitney region'' in $\Omega$, which we will construct by taking a union of certain dylated Whitney boxes. Hence, we ought to understand which Whitney boxes should be part of a Whitney region associated to the cube $Q$. The main two properties we desire to embed in such a region are, first, that it houses   the Corkscrew points for $Q$, and second, that these Corkscrew points are joined together by Harnack Chains that remain within the region. We will also want to fit in parameters that allow us to control the non-tangential aperture of these regions. In preparation to define these regions, we supply the technically relevant results.

\begin{lemma}[Whitney boxes contain Corkscrew points]\label{lm.corkscrewbox} Fix $Q\in\bb D$ and suppose\\ that $X\in\Omega$ is a Corkscrew point with Corkscrew constant $c>0$ for the surface ball $\Delta(x_Q,a_0\ell(Q)/2)$. Then there exists $I\in\m W$ such that $X\in I$ and satisfying
\begin{equation}\label{eq.whitneycork}
\frac{a_0c}{82\sqrt n}\ell(Q)\leq\ell(I)\leq\frac{a_0}{8\sqrt n}\ell(Q),\qquad \dist(I,Q)\leq a_0\ell(Q)/2.
\end{equation}
\end{lemma}

\noindent\emph{Proof.} Since $\m W$ is a covering of $\Omega$, it follows that there exists $I\in\m W$ such that $X\in I$. We now prove the bounds in (\ref{eq.whitneycork}). The upper estimate for $\ell(I)$ in (\ref{eq.whitneycork}) is deduced from the following chain of inequalities:
\begin{equation*}
4\sqrt n\ell(I)=4\diam I\leq \dist(I,\Gamma)\leq \dist(I,x_Q)\leq |X-x_Q|<a_0\ell(Q)/2,
\end{equation*}
where we have used (\ref{eq.whitney1}), $x_Q\in\Gamma$, $X\in I$, and that $X\in B(x_Q,a_0\ell(Q)/2)$. By observing that $\dist(I,Q)\leq\dist(I,x_Q)$, we arrive at the desired estimate for $\dist(I,Q)$. It remains now to give the lower bound for $\ell(I)$. Since $B(X,ca_0\ell(Q)/2)\subset\Omega$, we have that $\dist(X,\Gamma)\geq ca_0\ell(Q)/2$. Note that
\[
40\sqrt n\ell(I)=40\diam I\geq\dist(I,\Gamma)\geq\dist(X,\Gamma)-\diam I\geq ca_0\ell(Q)/2-\sqrt n\ell(I),
\]
whence the desired result follows.  \hfill{$\square$}

Thus, for each cube $Q\in\bb D$ the collection
\begin{equation}\label{eq.whitneyboxesc}
\m W_Q^{\operatorname{cs}}:=\Big\{I\in\m W\,:\,\frac{a_0c}{82\sqrt n}\ell(Q)\leq\ell(I)\leq\frac{a_0}{8\sqrt n}\ell(Q),\quad\text{dist}(I,Q)\leq a_0\ell(Q)/2\Big\}
\end{equation}
contains all the Corkscrew points for $\Delta(x_Q,a_0\ell(Q)/2)$ with Corkscrew constant $c$ (which evidently are Corkscrew points of $Q$). We note that without loss of generality, we may assume that a Corkscrew point for $Q$ is located at the center of some $I\in\m W_Q^{\operatorname{cs}}$ (with possibly smaller Corkscrew constant), as gives us the following result.

\begin{lemma}[Corkscrew points lie at the centers of Whitney boxes]\label{lm.centercork} Fix $Q\in\bb D$ and $X$ a Corkscrew point for the surface ball $\Delta(x_Q,a_0\ell(Q)/2)$ with Corkscrew constant $c>0$. If $I\in\m W_Q^{\operatorname{cs}}$ contains $X$, then $X_I$ is a Corkscrew point for $Q$ with Corkscrew constant $\tilde c=c/(1000\sqrt n)$. Moreover, $B(X_I,\tilde ca_0\ell(Q))\subset\opint(\frac12I)$.
\end{lemma}

\noindent\emph{Proof.} Suppose that $I\in\m W_Q^{\operatorname{cs}}$ contains $X$. Then $|X-X_I|\leq\diam I\leq a_0\ell(Q)/4$, and therefore,
\[
|X_I-x_Q|\leq|X_I-X|+|X-x_Q|<a_0\ell(Q)/4+a_0\ell(Q)/2<a_0\ell(Q).
\]
Hence $X_I\in B(x_Q,a_0\ell(Q))$. Now let $\tilde c$ be as in the statement of the lemma, and observe that for any $Y\in B(X_I,\tilde ca_0\ell(Q))$, we have that
\[
|X_I-Y|\leq\tilde ca_0\ell(Q)<\tfrac{a_0c}{82\sqrt n}\tfrac{\ell(Q)}8<\ell(I)/8,
\]
\[
|Y-x_Q|\leq\diam I+|X-x_Q|<a_0\ell(Q),
\]
as desired.\hfill{$\square$}

\begin{corollary}\label{cor.centercork} For any $Q\in\bb D$, there exists $I\in\m W_Q^{\operatorname{cs}}$ such that $X_I$ is a Corkscrew point for $Q$ with Corkscrew constant $\tilde c=\tilde c(c,n)$.
\end{corollary}
	
It may happen that $\m W_Q^{\operatorname{cs}}$ is too meager a region to use it to pass to ``continuous'' sawtooth-domains, or to pass between ``adjacent'' Whitney regions.  We introduce parameters $\eta\in(0,1)$ and $K\geq1$ and define
\begin{equation}\label{eq.whitneyboxes0}
\m W_Q^0:=\Big\{I\in\m W\,:\,\frac{a_0c}{82\sqrt n}\eta\ell(Q)\leq\ell(I)\leq\frac{a_0}{4\sqrt n}K\ell(Q),\quad\text{dist}(I,Q)\leq a_0K\ell(Q)\Big\}
\end{equation}
so that we may enlargen $\m W_Q^0$ according to aperture considerations. Immediately we have the following two technical results.

\begin{lemma}[Transversal adjacency of Whitney regions]\label{lm.whitneychild} If $\eta\in(0,c_{\bb K})$, then for any $Q\in\bb D$, the Whitney region $\m W_Q^0$ contains all Corkscrew points of the proper children of $Q$ (with Corkscrew constant $c$).
\end{lemma}

\noindent\emph{Proof.} Upon using (\ref{eq.largechild}), the proof is very similar to that of (\ref{eq.whitneycork}), and thus we omit it.\hfill{$\square$}

\begin{lemma}[Parallel adjacency of Whitney regions]\label{lm.closeregions} Fix $Q_1,Q_2\in\bb D$ and suppose that $\ell(Q_1)\leq\ell(Q_2)\leq c_{\bb K}^{-1}\ell(Q_1)$, and $\dist(Q_1,Q_2)\leq500\ell(Q_2)$. If $\eta\in(0,c_{\bb K})$ and $K\geq500A_0a_0^{-1}$, then $\m W_{Q_1}^0\cap\m W_{Q_2}^0\neq\varnothing$.
\end{lemma}

\noindent\emph{Proof.} Recall that $\m W_{Q_1}^{\operatorname{cs}}\neq\varnothing$, so fix $I\in\m W_{Q_1}^{\operatorname{cs}}$. It is easy to see that
\[
\frac{a_0c}{82\sqrt n}c_{\bb K}\ell(Q_2)\leq\ell(I)\leq\frac{a_0}{4\sqrt n}K\ell(Q_2),
\]
while the triangle inequality gives us that
\begin{multline*}
\dist(I,Q_2)\leq\dist(I,Q_1)+\diam Q_1+\dist(Q_1,Q_2)\leq a_0\ell(Q_1)/2+A_0\ell(Q_1)+500\ell(Q_2)\\ \leq 500A_0\ell(Q_2)\leq a_0K\ell(Q_2).
\end{multline*}
Thus $I$ verifies the conditions to be an element of $\m W_{Q_2}^0$.\hfill{$\square$}

Henceforth, we assume that
\begin{equation}\label{eq.parameters}\notag
\eta<c_{\bb K},\qquad K\geq500A_0a_0^{-1}.
\end{equation}

We need to augment $\m W_Q^0$ one final time: we must provide it with enough new boxes so that Harnack Chains connecting its old boxes are contained within a region that stays far from $\Gamma$. When $d<n-1$,  the technical result needed to accomplish this is

\begin{lemma}[Harnack Chains of $\m W_Q^0$ if $d<n-1$]\label{lm.hcwhitney} Fix $Q\in\bb D$, $X_Q$ a Corkscrew \\ point of $Q$ with Corkscrew constant $c>0$, and $I\in\m W^0_Q$ (note that $I$ may or may not contain $X_Q$). Then we may construct a well-tempered Harnack Chain  $\m H_I$ connecting $X_Q$ to $X_I$ consisting of a number at most $N_{\m H}$ of balls, where
\begin{equation}\label{eq.hccard}\notag
N_{\m H}=\tfrac4{c_{\m H}}\Big[\tfrac{41}{2c\eta}\big(1+\tfrac54K+\tfrac{A_0}{a_0}\big)\Big]^{\frac{n-1}{n-1-d}},
\end{equation}
and
\begin{equation}\label{eq.hcdist}\notag
\dist(\m H_I,\Gamma)\geq\frac12c_{\m H}\Big[\tfrac{41}{2c}\big(1+\tfrac54K+\tfrac{A_0}{a_0}\big)\Big]^{\frac{-d}{n-1-d}}\big(\tfrac2{41}a_0c\big)\eta^{\frac{n-1}{n-1-d}}\ell(Q).
\end{equation}
\end{lemma}

\noindent\emph{Proof.} Since $I\in\m W_Q^0$, we have that
\[
\delta(X_I)=\dist(X_I,\Gamma)\geq\dist(I,\Gamma)\geq4\sqrt n\ell(I)\geq\frac2{41}a_0c\eta\ell(Q).
\]
Furthermore, $\delta(X_Q)\geq cr_Q=a_0c\ell(Q)\geq\frac2{41}a_0c\ell(Q)$. Thus we take $s:=\frac2{41}a_0c\eta\ell(Q)$ in the setup of Lemma \ref{lm.harnackchain}. Next, we estimate
\begin{multline*}
|X_I-X_Q|\leq\diam I+\dist(I,Q)+\diam Q+|X_Q-x_Q|\\ \leq\sqrt n\ell(I)+a_0K\ell(Q)+A_0\ell(Q)+a_0\ell(Q)\\ \leq\big((\tfrac54K+1)a_0+A_0\big)\ell(Q)=\Big(\tfrac{41}{2c\eta}\big[1+\tfrac54K+\tfrac{A_0}{a_0}\big]\Big)\Big(\tfrac2{41}a_0c\eta\ell(Q)\Big).
\end{multline*}
Hence we take $\Lambda:=\frac{41}{2c\eta}\big[1+\frac54K+\frac{A_0}{a_0}\big]$. We now invoke the conclusion of Lemma \ref{lm.harnackchain} to find two points $Y_I\in B(X_I,a_0c\eta\ell(Q)/41), Y_Q\in B(X_Q,a_0c\eta\ell(Q)/41)$ such that 
\begin{equation}\label{eq.xi}
\dist([Y_I,Y_Q],\Gamma)\geq c_{\m H}\Big[\tfrac{41}{2c}\big(1+\tfrac54K+\tfrac{A_0}{a_0}\big)\Big]^{\frac{-d}{n-1-d}}\big(\tfrac2{41}a_0c\big)\eta^{\frac{n-1}{n-1-d}}\ell(Q)=:a_1\eta^{\frac{n-1}{n-1-d}}\ell(Q).
\end{equation}
Consider a finite covering of $[Y_I,Y_Q]$ by balls $B_j$ with centers in $[Y_I,Y_Q]$, radii all equal to $\frac{a_1}2\eta^{\frac{n-1}{n-1-d}}\ell(Q)$, and centers spaced by the radii. It is clear that the union of the balls $B_j$ is a well-tempered Harnack Chain $\m H_I$ satisfying the second desired estimate.  Let $N'_{\m H}$ be the cardinality of the covering. Since
\begin{equation*}
\frac{|Y_I-Y_Q|}{\frac{a_1}2\eta^{\frac{n-1}{n-1-d}}\ell(Q)}\leq\frac{|Y_I-X_I|+|X_I-X_Q|+|X_Q-Y_Q|}{\frac{a_1}2\eta^{\frac{n-1}{n-1-d}}\ell(Q)}\leq\frac{2\Lambda s}{\frac{c_{\m H}}2\Lambda^{\frac{-d}{n-1-d}}s},
\end{equation*}
and $N'_{\m H}a_1\eta^{\frac{n-1}{n-1-d}}\ell(Q)/2\leq|Y_I-Y_Q|+a_1\eta^{\frac{n-1}{n-1-d}}\ell(Q)/2$, the first desired estimate follows.\hfill{$\square$}

Actually, Harnack Chains cannot stray too far from the boxes in $\m W_Q^0$, as gives the following result.

\begin{lemma}[Non-degeneracy of boxes in Harnack Chains of $\m W_Q^0$]\label{lm.hcwnondeg} Suppose that\\ $J\in\m W$ intersects the Harnack Chain of Lemma \ref{lm.hcwhitney}. Then
\begin{equation}\label{eq.whitneyhc}
\frac{a_1}{82\sqrt n}\eta^{\frac{n-1}{n-1-d}}\ell(Q)\leq\ell(J)\leq\frac{2a_0K+A_0}{4\sqrt n}\ell(Q),\qquad\dist(J,Q)\leq(2a_0K+A_0)\ell(Q),
\end{equation}
where $a_1=a_1(n,d,C_d,c,c_{\m H},K)$ is the quantity defined in (\ref{eq.xi}).
\end{lemma}

\noindent\emph{Proof.} Given that $J\cap\m H_I\neq\varnothing$, then there exists $X\in J$ and $X\in\m H_I$, so that in particular $\dist(X,\Gamma)\geq\frac12a_1\eta^{\frac{n-1}{n-1-d}}\ell(Q)$. On the other hand,
\[
\dist(X,\Gamma)\leq\diam J+\dist(J,\Gamma)\leq\sqrt n\ell(J)+40\sqrt n\ell(J)=41\sqrt n\ell(J),
\]
so that the lower bound for $\ell(J)$ in (\ref{eq.whitneyhc}) follows immediately. Now, note that $4\sqrt n\ell(J)\leq\dist(J,\Gamma)\leq\dist(X,\Gamma)$. Since $X\in\m H_I$, there exists a ball $B$ such that $X\in B$, where $B$ has center $Y_B\in[Y_I,Y_Q]$ and has radius $\frac{a_1}2\eta^{\frac{n-1}{n-1-d}}\ell(Q)$. Consider the estimate
\begin{multline*}
\dist(Y_B,\Gamma)\leq\max\big\{|Y_I-x_Q|,|Y_Q-x_Q|\big\}\leq\tfrac1{41}a_0c\eta\ell(Q)+\max\big\{|X_I-x_Q|,|X_Q-x_Q|\big\}\\ \leq\tfrac1{41}a_0c\eta\ell(Q)+\dist(I,Q)+\diam Q+\diam I\leq\big[\tfrac1{41}c\eta+\tfrac54K+\tfrac{A_0}{a_0}\big]a_0\ell(Q),
\end{multline*}
so that
\begin{multline*}
\dist(J,\Gamma)\leq\dist(X,\Gamma)\leq\dist(Y_B,\Gamma)+\tfrac{a_1}2\eta^{\frac{n-1}{n-1-d}}\ell(Q)\\ \leq\big[\tfrac1{41}c\eta+\tfrac54K+\tfrac{A_0}{a_0}\big]a_0\ell(Q)+\tfrac{a_1}2\eta^{\frac{n-1}{n-1-d}}\ell(Q)\leq(2a_0K+A_0)\ell(Q).
\end{multline*}
This last estimate gives the rest of the bounds in (\ref{eq.whitneyhc}).\hfill{$\square$}

We proceed with the construction of the sawtooth domains. Fix $Q\in\bb D$ and let $X_Q$ be a Corkscrew point for $Q$, which belongs to some Whitney box in $\m W_Q^0$. For each $I\in\m W_Q^0$, we let $\m H_I$ be any well-tempered Harnack Chain connecting $X_I$ to $X_Q$ manufactured in Lemma \ref{lm.hcwhitney}. Then we let $\m W_Q$ be the set of all $J\in\m W$ which meet at least one of the Harnack Chains $\m H_I$ with $I\in\m W_Q^0$; that is, 
\begin{equation}\label{eq.whitneyboxes}
\m W_Q:=\big\{J\in\m W\,:\,\text{there exists }I\in\m W_Q^0\text{ for which }\m H_I\cap J\neq\varnothing\big\}.
\end{equation}
We clearly have that $\m W_Q^0\subset\m W_Q$, and we note that if $J\in\m W_Q$, then $J$ satisfies the assumptions of Lemma \ref{lm.hcwnondeg} and therefore (\ref{eq.whitneyhc}) holds, giving that
\begin{equation}\label{eq.finalwhitney}
\m W_Q\subseteq\big\{I\in\m W~:~a_2\ell(Q)\leq\ell(I)\leq A_2\ell(Q),\qquad\dist(I,Q)\leq4\sqrt nA_2\ell(Q)\big\},
\end{equation}
where $a_2, A_2$ are the corresponding multiplicative constants in the first inequality chain in (\ref{eq.whitneyhc}). In particular, once $\eta, K$ are fixed, for any $Q\in\bb D$, the cardinality of $\m W_Q$ is uniformly bounded, which is a corollary of the following result.

\begin{lemma}[Number of Whitney boxes in a Whitney region]\label{lm.numbergenw} Fix $Q\in\bb D$ and for positive numbers $\alpha,\beta,\gamma$ with $\alpha<\beta$, define the set
\[
\m W(\alpha,\beta,\gamma)=\big\{I\in\m W~:~\alpha\ell(Q)\leq\ell(I)\leq\beta\ell(Q),\qquad\dist(I,Q)\leq\gamma\ell(Q)\big\}.
\]
Then  
\begin{equation}\label{eq.cardgenw}
\operatorname{card}(\m W(\alpha,\beta,\gamma))\leq\Big[\frac{\sqrt n\beta+\gamma+A_0}{\alpha}\Big]^n|B(0,1)|
\end{equation}
\end{lemma}

\noindent\emph{Proof.} Fix $Q\in\bb D$, let $X\in I$ be arbitrary with $I\in\m W(\alpha,\beta,\gamma)$, and consider the estimate
\[
|X-x_Q|\leq\diam I+\dist(I,Q)+\diam Q\leq\sqrt n\beta\ell(Q)+\gamma\ell(Q)+A_0\ell(Q).
\]
It follows that $\cup_{I\in\m W(\alpha,\beta,\gamma)}I\subseteq B(x_Q,(\sqrt n\beta+\gamma+A_0)\ell(Q))$. Hence, note that
\begin{multline*}
\operatorname{card}(\m W(\alpha,\beta,\gamma))[\alpha\ell(Q)]^n=\sum_{I\in\m W(\alpha,\beta,\gamma)}[\alpha\ell(Q)]^n\leq\sum_{I\in\m W(\alpha,\beta,\gamma)}\ell(I)^n\\=\big|\textstyle\bigcup_{I\in\m W(\alpha,\beta,\gamma)}I\big|\leq|B(x_Q,(\sqrt n\beta+\gamma+A_0)\ell(Q))|=[(\sqrt n\beta+\gamma+A_0)\ell(Q)]^n|B(0,1)|.
\end{multline*}
The desired result follows immediately.\hfill{$\square$}

\begin{corollary}[Cardinality of $\m W_Q$]\label{cor.numberwq} For each  $Q\in\bb D$, we have that
\begin{equation}\label{eq.cardw}\notag
\operatorname{card}(\m W_Q)\leq\Big[\frac{5\sqrt nA_2+A_0}{a_2}\Big]^n|B(0,1)|=:N_0.
\end{equation}
\end{corollary}

Next, we choose a small dilation parameter $\theta\in(0,1)$ so that for any $I\in\m W$, the concentric dilation $I^*=(1+\theta)I$ still satisfies the Whitney property
\begin{equation}\label{eq.whitney4}\notag
\text{diam }I\approx\text{diam }I^*\approx\text{dist}(I^*,\Gamma)\approx\text{dist}(I,\Gamma),
\end{equation}
with uniform constants not depending on the choices of $\eta, K$. More precisely, it can be easily shown that, as long as $\theta\in(0,8)$, we have for each $I\in\m W$ the estimates
\begin{gather*}
\diam I\leq\diam I^*\leq(1+\theta)\diam I,\\ (1-\tfrac{\theta}8)\dist(I,\Gamma)\leq\dist(I^*,\Gamma)\leq\dist(I,\Gamma),\\ \tfrac1{40}\dist(I^*,\Gamma)\leq\diam I^*\leq2\tfrac{1+\theta}{8-\theta}\dist(I^*,\Gamma).
\end{gather*}

For later use, we record also that if $X\in\partial I^*$ and $\theta\in(0,1)$, then
\begin{equation}\label{eq.bdrystar}
2\diam I\leq\delta(X)\leq82\diam I.
\end{equation}

Moreover, by taking $\theta$ small enough we can guarantee that $\text{dist}(I^*,J^*)\approx\text{dist}(I,J)$ for every $I,J\in\m W$:

\begin{lemma}[Distances between dilated Whitney boxes]\label{lm.dylated} Suppose that $I,J\in\m W$. Then for each $\theta\in(0,1/(4\sqrt n))$, we have the estimates
\begin{equation}\label{eq.dylateddist}
\big[1-4\sqrt n\theta]\dist(I,J)\leq\dist(I^*,J^*)\leq\dist(I,J).
\end{equation}
Furthermore, if $\theta$ is as above and $I,J\in\m W$ are distinct, then $I^*\cap\frac12J=\varnothing$.
\end{lemma}

\noindent\emph{Proof.} Without loss of generality, suppose that $\ell(I)\geq\ell(J)$. If $I$ and $J$ are adjacent, so that $\partial I\cap\partial J\neq\varnothing$, then $\dist(I,J)=0$ and $\dist(I^*,J^*)=0$, so that (\ref{eq.dylateddist}) holds trivially. Now suppose that $I$ and $J$ are not adjacent. The upper bound in (\ref{eq.dylateddist}) holds trivially. By the Pidgeonholing Principle, we must have that $\dist(I,J)\geq\frac14\ell(I)$, for otherwise a point of $J$ lies in a Whitney box adjacent to $I$, which implies that $J$ is adjacent to $I$ and we have a contradiction. Next, let $X^*\in I^*$ and $Y^*\in J^*$ be the points such that $\dist(I^*,J^*)=|X^*-Y^*|$. For each $X\in I$ and $Y\in J$, observe the basic estimate
\begin{equation*}
|X^*-Y^*|\geq|X-Y|-|Y-Y^*|-|X^*-X|\geq\dist(I,J)-|Y-Y^*|-|X^*-X|.
\end{equation*}
Choose $X$ so that $|X-X^*|=\dist(X^*,I)$ and $Y$ so that $|Y-Y^*|=\dist(Y^*,J)$. Note that $\dist(X^*,I)\leq\frac{\theta}2\sqrt n\ell(I)$ and $\dist(Y^*,J)\leq\frac{\theta}2\sqrt n\ell(J)\leq\frac{\theta}2\sqrt n\ell(I)$. It follows that
\[
\dist(I^*,J^*)\geq\big[1-4\sqrt n\theta\big]\dist(I,J),
\]
which ends the proof of (\ref{eq.dylateddist}). Now suppose that $I,J\in\m W$ are distinct and without loss of generality say that $\ell(I)\geq\ell(J)$. If they are not adjacent, then $I^*\cap\frac12J=\varnothing$ follows immediately from (\ref{eq.dylateddist}). Hence we need only consider the case that $I$ and $J$ are adjacent, and in this case we have that $\ell(J)\geq\frac14\ell(I)$. Let $X^*\in I^*$ and $Y\in\frac12J$ be points such that $\dist(I^*,\frac12J)=|X^*-Y|$, and choose $X\in I$ so that $|X^*-X|=\dist(X^*,I)\leq\frac{\theta}2\sqrt n\ell(I)$. Reckon the elementary estimate
\begin{equation*}
|X^*-Y|\geq|Y-X|-|X-X^*|\geq\tfrac12\ell(J)-\tfrac{\theta}2\sqrt n\ell(I)\geq\big[\tfrac18-\tfrac{\theta}2\sqrt n\big]\ell(I)>0,
\end{equation*}
which yields the desired result.\hfill{$\square$}

\begin{corollary}\label{cor.adjacentdylated} Suppose that $I,J\in\m W$ and $\theta\in(0,1/(4\sqrt n))$. Then $I^*\cap J^*\neq\varnothing$ if and only if $\partial I\cap\partial J\neq\varnothing$.
\end{corollary}

Given   $Q\in\bb D$, we   define an associated \emph{Whitney region} $U_Q$ as
\begin{equation}\label{eq.whitneyregion}\notag
U_Q:=\medcup_{I\in\m W_Q}I^*.
\end{equation}

For any disjoint family $\m F=\{Q_j\}_j$, we define the \emph{discretized sawtooth relative to $\m F$} by
\[
\bb D_{\m F}:=\bb D\backslash\medcup_{\m F}\bb D_{Q_j},
\]
so that $\bb D_{\m F}$ is the collection of all $Q\in\bb D$ which are not contained in any $Q_j\in\m F$. We may also need to consider a local version of the sawtooth. If $\m F\subset\bb D_{Q_0}\backslash\{Q_0\}$ is as above, then we define the \emph{(local) discretized sawtooth relative to $\m F$} by
\begin{equation}\label{eq.localsawtooth}\notag
\bb D_{\m F,{Q_0}}:=\bb D_{Q_0}\backslash\medcup_{\m F}\bb D_{Q_j}.
\end{equation}
Finally, we define the global and local \emph{sawtooth domains} relative to $\m F$ via
\begin{equation}\label{eq.sawtoothdomain}\notag
\Omega_{\m F}:=\textup{int}\Big(\bigcup_{Q\in\bb D_{\m F}}U_{Q}\Big),\qquad\qquad\Omega_{\m F,{Q_0}}:=\text{int}\Big(\bigcup_{Q\in\bb D_{\m F,{Q_0}}}U_{Q}\Big).
\end{equation}
For convenience, we set
\begin{equation}\label{eq.whitneyf}\notag
\m W_{\m F}:=\medcup_{Q\in\bb D_{\m F}}\m W_{Q},\qquad\qquad\m W_{\m F,{Q_0}}:=\bigcup_{Q\in\bb D_{\m F,{Q_0}}}\m W_{Q},
\end{equation}
so that in particular, we may write
\begin{equation}\label{eq.domainf}
\Omega_{\m F}=\text{int}\Big(\bigcup_{I\in\m W_{\m F}}I^*\Big),\qquad\qquad\Omega_{\m F,{Q_0}}=\text{int}\Big(\bigcup_{I\in\m W_{\m F,{Q_0}}}I^*\Big).
\end{equation}

We remark that
\begin{equation}\label{eq.sawtoothbdd}
\Omega_{\m F,Q_0}\subset B(x_{Q_0}, 7\sqrt n A_2\ell(Q_0))\cap\Omega
\end{equation}
for any $Q_0\in\bb D$ and any family $\m F$ as above. Indeed, let $X\in\Omega_{\m F,Q_0}$, so that there exists $Q\in\bb D_{\m F,Q_0}$ and $I\in\m W_Q$ with $X\in I^*$. By (\ref{eq.finalwhitney}) and the triangle inequality, we see that
\begin{multline*}
|X-x_{Q_0}|\leq\diam I^*+\dist(I^*,Q)+\diam Q_0\leq2\sqrt n\ell(I)+4\sqrt nA_2\ell(Q)+A_0\ell(Q)\\ \leq7\sqrt nA_2\ell(Q_0),
\end{multline*}
as desired.

\subsection{Some further notation and auxiliary results} Let $\m F=\{Q_j\}_j$ be a disjoint family. The definitions below will be stated for the global discretized sawtooth relative to $\m F$, but it is clear that we have direct analogues for the local discrete sawtooth.

$\bullet$ We denote by $\Delta_\star$ a surface ball on $\partial\Omega_{\m F}$. More precisely, suppose that $x_\star\in\partial\Omega_{\m F}$ and $r>0$. Then
\[
\Delta_\star(x_\star,r):=B(x_\star,r)\cap\partial\Omega_{\m F}.
\]

$\bullet$ Let $\delta_\star:\Omega_{\m F}\ra[0,\infty)$ be the distance to $\partial\Omega_{\m F}$; that is,
\[
\delta_\star(X):=\dist(X,\partial\Omega_{\m F}).
\]

$\bullet$ We denote $\Sigma:=\partial\Omega_{\m F}\backslash\Gamma$, and observe that
\begin{equation}\label{eq.sigma}\notag
\Sigma=\partial\Omega_{\m F}\backslash\Gamma=\partial\big(\textstyle\bigcup_{I\in\m W_{\m F}}I^*\big)\backslash\Gamma\subset\displaystyle\bigcup\limits_{\tiny\begin{matrix}I\in\m W_{\m F}\\ I\cap\partial\Omega_{\m F}\neq\varnothing\end{matrix}}\partial I^*,
\end{equation}
so that $\Sigma$ consists of subsets of $(n-1)-$dimensional faces of Whitney boxes $I\in\m W_{\m F}$.

$\bullet$ For each $Q_0\in\bb D$, we let the \emph{Carleson collection associated to $Q_0$} be
\begin{equation}\label{eq.carlesonfam}
\m R_{Q_0}:=\medcup_{Q\in\bb D_{Q_0}}\m W_{Q}.
\end{equation}
We define the \emph{Carleson region} $R_{Q_0}$ as
\begin{equation}\label{eq.carlesonbox}\notag
R_{Q_0}:=\opint\Big(\bigcup_{Q\in\bb D_{Q_0}}U_{Q}\Big)=\opint\Big(\bigcup_{I\in\m R_{Q_0}}I^*\Big).
\end{equation}

$\bullet$ We consider a discrete dyadic version of the approach region rather than the standard non-tangential cone described in Definition \ref{def.ntmax}. For every $x\in\Gamma$, we define the (global and local) \emph{dyadic non-tangential cones} as
\begin{equation}\label{eq.localdyadicntcone}
\gamma_d(x)=\bigcup_{Q\in\bb D:Q\ni x}U_Q,\qquad \gamma_d^{Q_0}(x)=\bigcup_{Q\in\bb D_{Q_0}:Q\ni x}U_{Q}
\end{equation}
Given an aperture $\alpha>0$, there exists $K$ (in the definition (\ref{eq.whitneyboxes0})) sufficiently large such that the standard non-tangential cone $\gamma^{\alpha}(x)\subset\gamma_d(x)$ for all $x\in\Gamma$; and vice versa, for fixed values of $\eta, K$ and the dilation constant $\theta$, there exists $\alpha_1>0$ such that the dyadic cone $\gamma_d(x)\subset\gamma^{\alpha_1}(x)$ for all $x\in\Gamma$.

\subsection{Properties of the dyadically-generated sawtooth domains}

In this subsection, we collect a number of technical results, many of which are direct analogues of results shown in \cite{hm2}. There, the authors work with chord-arc domains with boundaries of co-dimension $1$. We are interested in borrowing their setup; the proofs of many results here are very similar to theirs, with some small modifications.

The first lemma we wish to present says that the boundary of a union of Whitney boxes consists of hyper-rectangles that do not degenerate.

\begin{lemma}[Non-degeneracy of the faces of $\Sigma$]\label{lm.faces} For all $\theta\in(0,1/(16\sqrt n))$ and for each $I\in\m W_{\m F}$ intersecting $\partial\Omega_{\m F}$, the set 
\[
\partial_{\Sigma}I^*:=\partial I^*\backslash\cup_{J\in\m W_{\m F}, J\neq I}\opint J^*=(\partial I^*\cap\partial\Omega_{\m F})\backslash\cup_{J\in\m W_{\m F}, J\neq I}\opint J^*
\]
 is a non-empty union of $(n-1)-$dimensional rectangles $R$ embedded in the $(n-1)$-dimensional faces of $\partial I^*$, such that no sidelength of any $R$ is smaller than $\theta\ell(I)/4$, thus verifying
\begin{equation}\label{eq.rectangle}\notag
\n H^{n-1}(R)\geq c_n\theta^{n-1}\ell(I)^{n-1},
\end{equation}
where $c_n\in(0,1)$ is a uniform constant depending only on $n$.
\end{lemma}

\noindent\emph{Proof.} Suppose that $I\in\m W_{\m F}$ intersects $\partial\Omega_{\m F}$, whence it follows that some face of $I^*$ intersects $\partial\Omega_{\m F}$. In the union defining $\partial_{\Sigma}I^*$, we need only consider those $J\in\m W_{\m F}$ which are adjacent to $I$, by Corollary \ref{cor.adjacentdylated}. That $\partial_{\Sigma}I^*$ is a union of $(n-1)-$dimensional rectangles is clear by the construction of the sawtooth domain, as the boundary can be written as a union of faces of cubes intersecting the complements of cubes.

Now let $F^*$ be any face of $I^*$ so that $F^*\cap\partial\Omega_{\m F}\neq\varnothing$, and let $F$ be the face of $I$ corresponding to $F^*$, defined as the unique face of $I$ such that $\dist(\opint F,F^*)=\frac12\theta\ell(I)$. Fix a maximal rectangle $R\subset F^*\cap\partial\Omega_{\m F}$ (maximal in the sense that increasing any of its side-lengths makes it stop being a subset of $\partial\Omega_{\m F}$) and consider two cases.

{\bf Case a)} There exists $x\in R$ and $x'\in\overline F$ such that $\dist(x,x')=\frac12\theta\ell(I)$. In this case, for each $\theta$ small enough, we must have that $x'\in\overline{J'}$, where $J'\notin\m W_{\m F}$ is a Whitney box adjacent to $I$. Since $J'$ is adjacent to $I$, we have that $\ell(J)\geq\ell(I)/4$, and recall that the Whitney boxes are dyadically aligned. It follows that $F\backslash_{J\in\m W_{\m F}, J\neq I}\overline J$ contains an $(n-1)-$dimensional cube $F'\subset\partial J'$ of length $\ell(F')\geq\ell(I)/4$. For each $y'\in F'$, there exists a unique $y\in F^*$ such that $\dist(y,y')=\frac12\theta\ell(I)$; accordingly we set ${F'}^*$ to be the collection of $y\in F^*$ constructed in this way, and observe that $x\in{F'}^*$. Then we must have that $R\supset{F'}^*\backslash\cup_{J\in\m W_{\m F}, J\neq I}\overline J$ contains an $(n-1)-$dimensional cube of side-length $(\frac14-4\theta)\ell(I)$, because for any straight line segment $L$ in $F^*$ parallel to a coordinate axis and passing through the center of ${F'}^*$, $L$ intersects at most two Whitney boxes $J_1,J_2\in\m W_{\m F}$ different from $I$ and such that $J_i^*\cap\overline{{F'}^*}\neq\varnothing$; both of which have $4\ell(I)\geq\ell(J_i)\geq\ell(I)/4$. Hence, in case a) we have established the desired result.

{\bf Case b)} There is no such $x'\in\overline F$ as in Case a). It follows that $R\subset F^*\backslash\frac1{1+\theta}F^*$, and therefore $R\cap\opint J'\neq\varnothing$ for some Whitney box $J'\notin\m W_{\m F}$ touching $I$. If $J\in\m W_{\m F}$ is any Whitney box adjacent to $I$ with $\ell(J)=2^{k}\ell(I)$ and such that $J^*\cap J'\neq\varnothing$, then $k\in\{-2,-1\}$. Indeed, if $k\in\bb Z$ is larger then $J^*$ protrudes a distance (perpendicular to the face $F$ whose boundary is intersected by $\overline J$) greater than or equal to $\frac12\theta\ell(I)$, so that $R\subset\overline{J^*}$. If $R\subset\opint J^*$ we have a contradiction to the fact that $R\subset\partial\Omega_{\m F}$; whereas if $R\subset\partial J^*$, then there is a face $F_{J^*}$ which is adjacent to $F^*$ and such that $\opint(F_{J^*}\cup F^*)$ is a connected set; we may then reduce to Case a) by considering that $R\subset F_{J^*}$. Finally, since $k\leq-1$, then $J^*$ protrudes a distance at most $\frac14\theta\ell(I)$. It follows that all the sides of $R$ have length larger than or equal to $\frac12\theta\ell(I)-\frac14\theta\ell(I)=\frac14\theta\ell(I)$, giving the result.\hfill{$\square$}

\begin{proposition}[Characterization of non-hidden regions, Proposition 6.1 in \cite{hm2}]\label{prop.containment} Let $\m F$ be a disjoint family. Then
	\begin{equation}\label{eq.containment}
	\Gamma\backslash\big(\medcup_{\m F}\,Q_j\big)\subseteq\Gamma\cap\partial\Omega_{\m F}\subseteq\Gamma\backslash\big(\medcup_{\m F}\opint Q_j\big).
	\end{equation}
\end{proposition}

\noindent\emph{Proof.} Let us show by contradiction the second containment first, thus assume that there exists $x\in\Gamma\cap\partial\Omega_{\m F}\cap\opint Q_j$ for some $Q_j\in\m F$. Hence, there exists $\ep>0$ for which $B(x,\ep)\cap\Gamma\subset Q_j$. In particular, $B(x,\ep)\cap Q=\varnothing$ for any $Q\in\bb D_{\m F}$ that does not contain $Q_j$. Since $x\in\partial\Omega_{\m F}$, there exist $X_k\in\Omega_{\m F}$ such that $|X_k-x|\ra0$ as $k\ra\infty$. Accordingly, for each $k\in\bb N$ there exists a Whitney box $I_k$ and a dyadic cube $Q_k\in\bb D_{\m F}$ such that $I_k\in\m W_{Q_k}$ and $X_k\in I_k^*$. Since $x\in\Gamma$ and $X_k\ra x$, then $\delta(X_k)\ra0$ and therefore $\ell(I_k^*)\ra0$, which further implies that $\ell(Q_k)\ra0$ by (\ref{eq.finalwhitney}). It follows that for all $k$ large enough, $\ell(Q_k)\ll\ell(Q_j)$, so that $B(x,\ep)\cap Q_k=\varnothing$. On the other hand,
\[
\dist(Q_k,x)\leq\dist(Q_k,I_k^*)+\diam I_k^*+|X_k-x|\lesssim\ell(Q_k)+\ell(I_k)+|X_k-x|\longrightarrow0,
\]
which implies that there exists $k_0\in\bb N$ large enough and a point $q\in Q_{k_0}$ so that $|q-x|<\ep$. Hence $q\in B(x,\ep)\cap Q_{k_0}$, but this is a contradiction. The second containment is thus established.

We now prove the first containment. Suppose that $x\in\Gamma\backslash(\cup_{\m F}Q_j)$, and note that obviously $x\notin\Omega_{\m F}$. Hence, for any generation $k\in\bb Z$, $x\in Q_k$ for some $Q_k\in\bb D_{\m F}\cap\bb D^k$. According to each $Q_k\in\bb D_{\m F}$, there exists $I_k\in\m W_{Q_k}$ and $I_k\subset\Omega_{\m F}$. Then the centers $X_{I_k}\in I_k$ satisfy
\[
|X_{I_k}-x|\leq\diam Q_k+\dist(Q_k,I_k)+\diam I_k\lesssim\ell(Q_k).
\]
Taking $k\ra\infty$ gives that $\ell(Q_k)\approx2^{-k}\ra0$, and hence that $|X_{I_k}-x|\ra0$. It follows that $x\in\partial\Omega_{\m F}$.\hfill{$\square$}

Although we may not quite say that $\Gamma\cap\partial\Omega_{\m F}\subseteq\Gamma\backslash\opint(\cup_{\m F}Q)$, we do have a technical improvement to the upper containment in (\ref{eq.containment}).


\begin{lemma}[A refinement to the characterization of Proposition \ref{prop.containment}]\label{lm.kcontain} Suppose that $\m F$ is a disjoint family and let $x\in\Gamma\cap\partial\Omega_{\m F}$. Then for each $k\in\bb Z$, there exists $Q^k\in\bb D_{\m F}\cap\bb D^k$ verifying that $x\in\overline{Q^k}$.
\end{lemma}

\noindent\emph{Proof}. Fix $k\in\bb Z$ and recall that $\{Q_m^k\}_{m\in\n J^k}$ is a disjoint covering of $\Gamma$. We first show that $x\in\overline{\medcup_{\bb D^k\cap\bb D_{\m F}}Q_m^k}$. Indeed, suppose otherwise, so that $x\in\opint(\medcup_{\bb D^k\backslash\bb D_{\m F}}Q_m^k)$. Then there exists $\ep>0$ so that $B(x,\ep)\cap\Gamma\subset\medcup_{\bb D^k\backslash\bb D_{\m F}}Q^k_m$. Since $x\in\Gamma\cap\partial\Omega_{\m F}$, as in the proof of Proposition \ref{prop.containment} we can procure dyadic cubes $Q'_i\in\bb D_{\m F}$ such that $\dist(Q'_i,x)\ra0$ and $\ell(Q'_i)\ra0$ as $i\ra\infty$. Accordingly, for all $i$ large enough, we have that $Q'_i\subset B(x,\ep)\cap\Gamma$ and $\ell(Q'_i)\ll2^{-k}\approx Q_m^k$ for any $Q_m^k\in\bb D^k\backslash\bb D_{\m F}$, which prohibits $Q_i'$ from being an ancestor of any $Q_m^k\in\bb D^k\cap\bb D_{\m F}$, and this implies that for all $i$ large enough, $Q_i'\subseteq Q_j$ for some $Q_j\in\m F$, yielding a contradiction to $Q_i'\in\bb D_{\m F}$.

To finish the proof, we observe that for any $y\in\Gamma$, the cardinality of the set 
\[
S(y):=\big\{Q_m^k\in\bb D_{\m F}\cap\bb D^k~:~\dist(Q_m^k,y)<2^{-k}\big\}
\]
is uniformly finite, implying in particular that $\overline{\medcup_{\bb D^k\cap\bb D_{\m F}}Q_m^k}=\medcup_{\bb D^k\cap\bb D_{\m F}}\overline{Q_m^k}$ and thus yielding the desired result. To see that $\card S(y)<+\infty$, reckon the estimate
\begin{multline*}
C_d^{-1}(a_02^{-k})^d\card(S(y))\leq\sum_{Q_m^k\in S(y)}\sigma(\Delta(x_{Q_m^k},a_02^{-k})\leq\sum_{Q_m^k\in S(y)}\sigma(Q_m^k)=\sigma\big(\medcup_{S(y)}Q_m^k\big)\\ \leq\sigma\big(B(x,2A_02^{-k})\cap\Gamma\big)\leq C_d(2A_02^{-k})^d,
\end{multline*}
so that
\[
\card(S(y))\leq 2^dC_d^2A_0^da_0^{-d}.
\] \hfill{$\square$}

The following lemma establishes that Carleson regions are quantitatively fat.

\begin{lemma}[Quantitative fatness of Carleson regions, Lemma 3.55 of \cite{hm2}]\label{lm.carlesonball} The following statements are true.
	\begin{enumerate}[(i)]
		\item\label{item.carleson} For each $Q\in\bb D$, there exists a ball $B_s:=B(x_Q,s)\subseteq B(x_Q,a_0\ell(Q))$ with $s\approx\ell(Q)$ (we may, in fact, take $s=a_0\ell(Q)/4$), such that
		\begin{equation*}\label{eq.carlesonball}
		B_s\cap\Omega\subset R_Q.
		\end{equation*}
		\item\label{item.carleson2} Moreover, for a somewhat smaller choice of $s\approx\ell(Q)$ (in fact, we may choose $s=a_0a_2A_2^{-1}\ell(Q)/10$), we have for every pairwise disjoint family $\m F\subset\bb D$, and for each $Q_0\in\bb D$ containing $Q$, that
		\begin{equation}\label{eq.carlesonball2}
		B_s\cap\Omega_{\m F,Q_0}=B_s\cap\Omega_{\m F,Q}.
		\end{equation}
	\end{enumerate}
\end{lemma}

\noindent\emph{Proof.} We consider \ref{item.carleson} first. Fix $Y\in B_s\cap\Omega$, and let $I\in\m W$ be a  Whitney box that contains $Y$. Choose $y\in\Gamma$ such that $|Y-y|=\dist(Y,\Gamma)=\delta(Y)\approx\ell(I)$, and observe that
\[
|y-x_Q|\leq|Y-y|+|Y-x_Q|\leq\delta(Y)+s\leq2s<a_0\ell(Q),
\]
provided that $s<a_0\ell(Q)/2$. Hence $y\in\Delta(x_Q,a_0\ell(Q))\subseteq Q$. Now let $Q^d\subset\bb D_Q$ be a descendant of $Q$ (which is unique \emph{as a set}) that contains $y$ and that verifies the inequalities
\[
\tfrac{a_0cc_{\bb K}}{82\sqrt n}\ell(Q^d)\leq\ell(I)\leq\tfrac{a_0}{8\sqrt n}\ell(Q^d).
\]
Since we must have that
\[
\dist(I,Q^d)\leq|Y-y|=\delta(Y)\leq41\diam I\leq\tfrac{41}8a_0\ell(Q^d)\leq a_0K\ell(Q^d),
\]
we conclude that $I\in\m W_{Q^d}^0$, and therefore $Y\in U_{Q^d}\subset R_Q$, as desired.

Now we consider \ref{item.carleson2}. Since $Q\subseteq Q_0$, the containment  $B_s\cap\Omega_{\m F,Q_0}\supseteq B_s\cap\Omega_{\m F,Q}$ is trivial, and thus we need only verify the opposite containment. Fix $Y\in B_s\cap\Omega_{\m F,Q_0}$, and as such there exist $Q'\in\bb D_{\m F,Q_0}$ and $I\in\m W_{Q'}$ verifying $Y\in I^*$. Note that
\[
|Y-x_Q|\geq\dist(I^*,\Gamma)\geq\tfrac12\dist(I,\Gamma)\geq2\diam I=2\sqrt n\ell(I),
\]
so that $\ell(I)<s/(2\sqrt n)$. We claim that $Q'\subseteq Q$. To see this, let $z\in Q'$ and reckon the estimate
\begin{multline*}
|z-x_Q|\leq\diam Q'+\dist(I^*,Q')+\diam I^*+|Y-x_Q|\\ < A_0a_2^{-1}\ell(I)+4\sqrt n A_2a_2^{-1}\ell(I)+2\sqrt n\ell(I)+s~ < 5A_2a_2^{-1}s<a_0\ell(Q),
\end{multline*} 
provided that $s<\frac15a_0a_2A_2^{-1}\ell(Q)$. Hence $Q'\subseteq\Delta_Q\subseteq Q$, as claimed. But then, $Q'\in\bb D_{\m F,Q}$, so that $I\in\m W_{\m F,Q}$ and therefore $I^*\subset\Omega_{\m F,Q}$, which implies that $Y\in\Omega_{\m F,Q}$, as desired.\hfill{$\square$}

The next proposition gives us that the sets in (\ref{eq.containment}) are the same from the perspective of a doubling measure. Its proof is essentially the same as in \cite{hm2}, thus we omit the details (see Remark \ref{rm.closetohm}).

\begin{proposition}[Negligibility of pathologies in (\ref{eq.containment}) for doubling measures, Proposition 6.3 in \cite{hm2}]\label{prop.mumeasure0}  Suppose that $\mu$ is a doubling measure on $\Gamma$. Then $\partial Q:=\overline Q\backslash\textup{int }Q$ has $\mu-$measure $0$, for every $Q\in\bb D$. In particular, the sets in (\ref{eq.containment}) have the same $\mu$ measure.
\end{proposition}

We will now elicit the existence of a point that acts as a Corkscrew point simultaneously in $\Omega$ and in the sawtooth domain $\Omega_{\m F,Q_0}$.

\begin{proposition}[Existence of simultaneous Corkscrews, Proposition 6.4 of \cite{hm2}]\label{prop.simultcs}  Fix $Q_0\in\bb D$, and let $\m F\subset\bb D_{Q_0}$ be a disjoint family. Then for each $Q\in\bb D_{\m F,Q_0}$, there is a radius $r_Q\approx\ell(Q)$ (in fact, we may take $r_Q=6\sqrt nA_2\ell(Q)$), and a point $X_Q\in\Omega_{\m F,Q_0}$ which serves as a Corskcrew point (with Corkscrew constant $\tilde c:=\frac1{12\sqrt n}a_2A_2^{-1}$) simultaneously for $\Omega_{\m F,Q_0}$, with respect to the surface ball $\Delta_\star(y_Q,r_Q)$, for some $y_Q\in\partial\Omega_{\m F,Q_0}$, and for $\Omega$, with respect to each surface ball $\Delta(x,r_Q)$, for every $x\in Q$.
\end{proposition}

\noindent\emph{Proof.} Fix $Q\in\bb D_{\m F,Q_0}$, and so note that there exists $I\in\m W_Q$ with $\opint I^*\subset\Omega_{\m F,Q_0}$. We fix this Whitney cube $I$. Observe that 
\[
\dist(X_I,\partial\Omega_{\m F,Q_0})\geq\dist(X_I,\partial I^*)=\tfrac12(1+\theta)\ell(I)\geq\tfrac12a_2\ell(Q),
\]
while on the other hand, since $X_I\in\opint\Omega_{\m F,Q_0}$ and $\Gamma\subset\bb R^n\backslash\Omega_{\m F,Q_0}$, we have that
\[
\dist(X_I,\partial\Omega_{\m F,Q_0})=\dist(X_I,\bb R^n\backslash\Omega_{\m F,Q_0})\leq\dist(X_I,\Gamma)\leq\dist(X_I,Q)\leq5\sqrt nA_2\ell(Q).
\]
It follows that, if we let $y\in\partial\Omega_{\m F,Q_0}$ be a point that satisfies $\dist(X_I,\partial\Omega_{\m F,Q_0})=|X_I-y|$, then $B(y,r_Q)\supseteq I\supseteq B(X_I,\tilde cr_Q)$, where $\tilde c$ and $r_Q$ are as in the statement of the proposition. Moreover, for any $x\in Q$, note that
\begin{multline*}
|x-X_I|\leq\diam Q+\dist(I,Q)+\diam I\leq A_0\ell(Q)+4\sqrt nA_2\ell(Q)+\sqrt n A_2\ell(Q)\\ \leq 6\sqrt n A_2\ell(Q),
\end{multline*}
whence it is easy to see that $B(x,r_Q)\supset B(X_I,\tilde cr_Q)$. Letting $X_Q=X_I$ and $y_Q=y$ finishes the proof.\hfill{$\square$}

Owing to (\ref{eq.sawtoothbdd}), when $Q=Q_0$ in the above proposition, we have

\begin{corollary}[A uniform Corkscrew point]\label{cor.bigcorkscrew} The point $X_{Q_0}$ given by Proposition \ref{prop.simultcs} is a Corkscrew point (with Corkscrew constant $\tilde c=\frac1{12\sqrt n}a_2A_2^{-1}$) with respect to $\Delta_\star(y,r_{Q_0})$ for all $y\in\partial\Omega_{\m F,Q_0}$, and for $\Delta(x,r_{Q_0})$, for all $x\in Q_0$, with $r_{Q_0}=7\sqrt n A_2\ell(Q_0)$.
\end{corollary}
 
The next lemma establishes a fatness of the region ``hidden'' by the sawtooth boundary; hence the next result has a similar spirit to (i) of Lemma \ref{lm.carlesonball}.

\begin{lemma}[Quantitative fatness of hidden regions, Lemma 5.9 in \cite{hm2}]\label{lm.Bbelowsawtooth}  Let $\m F\subset\bb D$ be a disjoint family. Then for every $Q\subseteq Q_j\in\m F$, there is a ball $B'\subset\bb R^n\backslash\Omega_{\m F}$, centered at $\Gamma$, with radius $r'\approx\ell(Q)$, and $\Delta':=B'\cap\Gamma\subset Q$. In fact, we may take $B'=B(x_Q,r')$ and $r'=a_0\ell(Q)/(5A_2a_2^{-1})$.
\end{lemma}

\noindent\emph{Proof.} Recall that $\Delta_Q=B(x_Q,a_0\ell(Q))\cap\Gamma\subset Q$. Let $B_M:=B(x_Q,a_0\ell(Q)/M)$, where $M=5A_2a_2^{-1}$. We claim that $B_M$ is the ball $B'$ with the desired properties. We need only check that $B_M\subset\bb R^n\backslash\Omega_{\m F}$, and we proceed via proof by contradiction. Thus suppose that there exists $I\in\m W_{\m F}$ with $I^*\cap B_M\neq\varnothing$, so that we may find $Y\in I^*\cap B_M$ and $Q_I\in\bb D_{\m F}$ with $I\in\m W_{Q_I}$. Then $\delta(y)< a_0\ell(Q)/M$, and therefore
\begin{gather*}
\diam I^*\leq2\diam I\leq \dist(I,\Gamma)/2\leq\dist(I^*,\Gamma)\leq\delta(y)<a_0\ell(Q)/M,\\ \dist(I^*,Q_I)\leq\dist(I,Q_I)\leq4\sqrt nA_2\ell(Q_I)\leq4\sqrt nA_2a_2^{-1}\ell(I)\leq2A_2a_2^{-1}a_0\ell(Q)/M.
\end{gather*}
It follows that
\begin{equation*}
\dist(Q_I,x_Q)\leq\dist(I^*,Q_I)+\diam I^*+|y-x_Q|\leq\tfrac{4A_2a_2^{-1}}Ma_0\ell(Q),
\end{equation*}
and so for any $q_I\in Q_I$, we have that
\begin{multline*}
|q_I-x_Q|\leq\diam Q_I+\dist(Q_I,x_Q)\leq\tfrac{A_0a_2^{-1}}{\sqrt n}\diam I+\tfrac{4A_2a_2^{-1}}Ma_0\ell(Q)\\ \leq\tfrac{9A_2a_2^{-1}/2}{M}a_0\ell(Q) <a_0\ell(Q),
\end{multline*}
which implies that $Q_I\subset\Delta_Q\subset Q\subset Q_j\in\m F$, but this is a contradiction to the assumption that $Q_I\in\bb D_{\m F}$. The desired result ensues.\hfill{$\square$}

We would now like to fix an $(n-1)-$ dimensional rectangle of the boundary of the Carleson region $R_{Q_j}$ which is morally a ``lift'' of $\Delta_{Q_j}$. The precise statement is as follows.

\begin{proposition}[Lift of $\Delta_{Q_j}$, Proposition 6.7 in \cite{hm2}]\label{prop.cubeinface} Let $\m F$ be a disjoint family. Then for each $Q_j\in\m F$, there is an $(n-1)-$dimensional cube $P_j\subset\partial\Omega_{\m F}$, which is contained in a face of $I^*$ for some $I\in\m W$, and that satisfies
	\begin{equation}\label{eq.cubeinface}
	\ell(P_j)\approx\dist(P_j,Q_j)\approx\dist(P_j,\Gamma)\approx\ell(I)\approx\ell(Q_j),
	\end{equation}
	with uniform constants.
\end{proposition}

\noindent\emph{Proof.} Fix $Q_j\in\m F$ and let $\hat Q$ be its proper parent, so that $\ell(Q_j)<\ell(\hat Q)\leq c_{\bb K}^{-1}\ell(Q_j)$ and $\hat Q\in\bb D_{\m F}$. Let $I\in\m W_{\hat Q}^{\operatorname{cs}}$, and in particular $\opint I^*\subset\Omega_{\m F}$. On the other hand, by Lemma \ref{lm.Bbelowsawtooth}, the ball $B=B(x_{Q_j}, a_0a_2\ell(Q_j)/(5A_2))$ lies in $\bb R^n\backslash\Omega_{\m F}$, so that if $X$ is a Corkscrew point (with Corkscrew constant $c$) for the surface ball $B\cap\Gamma$, then there exists $J\in\m W_{Q_j}$ such that $X\in J$ and $J\backslash\Omega_{\m F}\neq\varnothing$, whence by the construction of the sawtooth domain we have that $\Omega_{\m F}\cap\frac12J=\varnothing$.

Let $X_I$ and $X_J$  be the centers of $I$ and $J$, respectively. We will connect these points via a Harnack Chain. Reckon the estimates
\begin{gather*}
\tfrac2{41}a_0c\ell(Q_j)\leq\delta(X_I)\leq\tfrac{41}4a_0c_{\bb K}^{-1}\ell(Q_j),\\ \tfrac1{100}a_0a_2cA_2^{-1}\ell(Q_j)\leq\delta(X_J)\leq\tfrac14a_0a_2A_2^{-1}\ell(Q_j),
\end{gather*}
\begin{multline*}
|X_I-X_J|\leq\diam I+\dist(I,\hat Q)+\diam\hat Q+|x_{Q_j}-X|+|X-X_J|\\ \leq\tfrac14a_0\ell(\hat Q)+a_0\ell(\hat Q)+A_0\ell(\hat Q)+\tfrac14a_0a_2A_2^{-1}\ell(Q_j)\leq 4A_0c_{\bb K}^{-1}\ell(Q_j)\\ \leq\Big[\frac{400A_2^2}{a_2^2cc_{\bb K}}\Big]\big(\tfrac1{100}a_0a_2cA_2^{-1}\ell(Q_j)\big).
\end{multline*}
We thus apply Lemma \ref{lm.harnackchain} with $s:=\tfrac1{100}a_0a_2cA_2^{-1}\ell(Q_j)$ and $\Lambda:=400A_2^2/(a_2^2cc_{\bb K})$ and fix a well-tempered Harnack Chain $\m H$ (cf. Definition \ref{def.hc}) connecting $X_I$ and $X_j$ with $\Lambda, s$ as above. Since $X_I\in\opint\Omega_{\m F}$, $X_J\in\opint(\bb R^n\backslash\Omega_{\m F})$, and $\m H$ is centered on a straight line segment, it follows that $\m H$ intersects $\partial\Omega_{\m F}$ at some $Z\in\partial\Omega_{\m F}$. By the construction of the well-tempered Harnack Chain, $Z$ verifies
\begin{equation}\label{eq.zver}
\tfrac1{200}c_{\m H}\Lambda^{\frac{-d}{n-1-d}}a_0a_2cA_2^{-1}\ell(Q_j)\leq\delta(Z)\leq13a_0c_{\bb K}^{-1}\ell(Q_j).
\end{equation}
Since $Z\in\partial\Omega_{\m F}$, there exists $I_Z\in\m W_{\m F}$ such that $Z\in\partial I_Z^*$, which implies by (\ref{eq.bdrystar}) and (\ref{eq.zver}) that
\begin{equation}\label{eq.izver}
\tfrac1{82\sqrt n}\tfrac1{200}c_{\m H}\Lambda^{\frac{-d}{n-1-d}}a_0a_2cA_2^{-1}\ell(Q_j)\leq\ell(I_Z)\leq\tfrac{13}{2\sqrt n}a_0c_{\bb K}^{-1}\ell(Q_j),
\end{equation}
\begin{equation*} 
\tfrac1{(50)(82)}c_{\m H}\Lambda^{\frac{-d}{n-1-d}}a_0a_2cA_2^{-1}\ell(Q_j)\leq\dist(I_Z,\Gamma)\leq260a_0c_{\bb K}^{-1}\ell(Q_j)
\end{equation*}
By Lemma \ref{lm.faces}, there exists an $(n-1)-$dimensional cube $P_j\subset\partial I_Z^*$ of side-length $\theta\ell(I_Z)/4$ that contains $Z$, and hence
\begin{equation}\label{eq.pjver}
\tfrac1{82\sqrt n}\tfrac1{800}\theta c_{\m H}\Lambda^{\frac{-d}{n-1-d}}a_0a_2cA_2^{-1}\ell(Q_j)\leq\ell(P_j)\leq\tfrac{13}{8\sqrt n}\theta a_0c_{\bb K}^{-1}\ell(Q_j)
\end{equation}
Since $\dist(I_Z,\Gamma)\leq\dist(P_j,\Gamma)\leq\dist(I_Z,\Gamma)+\diam I_Z^*$, it follows that 
\begin{equation}\label{eq.pjver2}
\tfrac1{500}c_{\m H}\Lambda^{\frac{-d}{n-1-d}}a_0a_2cA_2^{-1}\ell(Q_j)\leq\dist(P_j,\Gamma)\leq300a_0c_{\bb K}^{-1}\ell(Q_j).
\end{equation}
Finally, we consider $\dist(P_j,Q_j)$. The lower bound is immediate from (\ref{eq.pjver2}) and the fact that $\dist(P_j,Q_j)\geq\dist(P_j,\Gamma)$. As for the upper bound, we first note that $\dist(X_I, Q_j)\geq\dist(X_J,Q_j)$ and $\diam I\geq\diam J$, so that by the construction of the Harnack Chain $\m H$, we have that 
\begin{multline*}
\dist(I_Z,Q_j)\leq\dist(Z,Q_j)\leq\tfrac12c_{\m H}\Lambda^{-\frac d{n-1-d}}s+\dist(I,Q_j)+\diam I \\ \leq \tfrac12c_{\m H}\Lambda^{-\frac d{n-1-d}}s+\dist(I,\hat Q)+\diam\hat Q+\diam I\leq 4A_0c_{\bb K}^{-1}\ell(Q_j).
\end{multline*}
Therefore, we conclude that
\begin{equation}\label{eq.pjver3}
\tfrac1{5000}c_{\m H}\Lambda^{\frac{-d}{n-1-d}}a_0a_2cA_2^{-1}\ell(Q_j)\leq\dist(P_j,Q_j)\leq11A_0c_{\bb K}^{-1}\ell(Q_j).
\end{equation}
The estimates (\ref{eq.izver}), (\ref{eq.pjver}), (\ref{eq.pjver2}), and (\ref{eq.pjver3}) readily imply (\ref{eq.cubeinface}).\hfill{$\square$}


\begin{notation}\label{not.centerofpj} Let $\m F$ be a disjoint family and $Q_j\in\m F$. Let $P_j\subset\partial\Omega_{\m F}$ be the $(n-1)-$dimensional cube constructed in Proposition \ref{prop.cubeinface} and satsifying (\ref{eq.cubeinface}). We denote by $x^\star_j$ the \emph{center} of $P_j$, and we write $r_j:=7\sqrt nA_2\ell(Q_j)$. With these choices, we have that $P_j\subseteq\Delta_\star(x^\star_j,r_j)$ and that
	\begin{equation}\label{eq.bigball}\notag
	\overline{R_{Q_j}}\subset B(x_j^\star,r_j),
	\end{equation}
by an argument similar to the one giving (\ref{eq.sawtoothbdd}). Moreover, given $Q\in\bb D_{\m F}$ and $y_Q\in\partial\Omega_{\m F}$ as in (an analogous global version of) Proposition \ref{prop.simultcs}, then we may choose $\hat r_Q\approx r_Q$ (in fact, we may take $\hat r_Q=42\sqrt n c_{\bb K}^{-1}A_2\ell(Q)$) so that the containment
\begin{equation}\label{eq.bigqball}
Q\medcup\Big(\medcup_{Q_j\in\m F: Q_j\subset Q}B(x_j^\star,r_j)\Big)\subset B(y_Q,\hat r_Q)
\end{equation}
holds. Indeed, it is easy to see that $Q\subset B(y_Q,\hat r_Q)$, while if $Q_j\in\m F$ with $Q_j\subset Q$, then for any $z\in B(x_j^\star,r_j)$ we use the bound
\[
|z-y_Q|\leq\diam B(x_j^\star,r_j)+\dist(P_j,Q_j)+\diam Q+|x_Q-y_Q|\lesssim r_Q.
\]
Henceforth, we fix $y_Q,\hat r_Q$ as in this paragraph.
\end{notation}

We conclude this section with the fact that there is a ``lift'' of any $\Delta_Q, Q\in\bb D_{\m F}$ which does not intersect Carleson regions of $Q_j\in\m F$, $Q_j$ not contained in $Q$.

\begin{proposition}[A lift of $\Delta_Q$, Proposition 6.12 of \cite{hm2}]\label{prop.starball} Let $\m F$ be a disjoint family. For $Q_j\in\m F$, let $B(x^\star_j,r_j)$ be the ball described in Notation \ref{not.centerofpj}. Then for each $Q\in\bb D_{\m F}$, there is a surface ball
	\begin{equation}\label{eq.starball}\notag
	\Delta_\star^Q:=\Delta_\star(x_Q^\star,t_Q)\subset\big(Q\cap\partial\Omega_{\m F}\big)\medcup\Big(\medcup_{Q_j\in\m F:Q_j\subset Q}\big(B(x_j^\star,r_j)\cap\partial\Omega_{\m F}\big)\Big),
	\end{equation}
	with $t_Q\approx\ell(Q), x^\star_Q\in\partial\Omega_{\m F}$, and $\dist(Q,\Delta_\star^Q)\lesssim\ell(Q)$.
\end{proposition} 

\noindent\emph{Proof.} Fix $M\geq$ a sufficiently large number to be chosen momentarily. We split the proof in two cases.

{\bf Case 1.} There exists $Q_j\subset Q$, for which $\ell(Q_j)\geq\ell(Q)/M$. In this case, we set $\Delta_\star^Q=\Delta_\star(x_k^\star,\ell(P_j)/2)$, where $P_j$ is the $(n-1)-$dimensional cube established in Proposition \ref{prop.cubeinface}.

{\bf Case 2.} There is no $Q_j$ as in Case 1. In this case, if $Q_j\cap Q\neq\varnothing$, then $Q_j\subset Q$ and $\ell(Q_j)<\ell(Q)/M$.

\emph{Sub-case 1.} No $Q_j\in\m F$ meets the surface ball $\Delta(x_Q,a_0\ell(Q)/(4\sqrt M)$. Then, we simply set $\Delta_\star^Q:=\Delta(x_Q,\ell(Q)/(4\sqrt M))$, and we reckon that $\Delta_\star^Q\subset\Delta_Q\subset Q\cap\partial\Omega_{\m F}$ by Proposition \ref{prop.containment}.

\emph{Sub-case 2.} There exists $Q_k\in\m F$ which meets the surface ball $\Delta(x_Q,a_0\ell(Q)/(4\sqrt M)$.  We claim that for all $M$ large enough, we have that $P_k\subset B(x_Q,a_0\ell(Q)/(2\sqrt M))$. Indeed, suppose that $y\in P_k$. Then, using  (\ref{eq.cubeinface}), we deduce that
\begin{multline*}
|y-x_Q|\leq\diam P_k+\dist(P_k,Q_k)+\diam Q_k+\dist(Q_k,x_Q)\\ \leq\tfrac{13}8\theta a_0c_{\bb K}^{-1}\ell(Q_k)+11A_0c_{\bb K}^{-1}\ell(Q_k)+A_0\ell(Q_k)+\tfrac1{4\sqrt M}a_0\ell(Q)\\ <13A_0c_{\bb K}^{-1}\tfrac1M\ell(Q)+\tfrac1{4\sqrt M}a_0\ell(Q)=\tfrac{13A_0a_0^{-1}c_{\bb K}^{-1}}{\sqrt M}\tfrac1{4\sqrt M}a_0\ell(Q)+\tfrac1{4\sqrt M}a_0\ell(Q)\\ <\frac1{2\sqrt M}a_0\ell(Q),
\end{multline*}
provided that $M>(13A_0a_0^{-1}c_{\bb K}^{-1})^2$. Thus the claim is shown, and accordingly, we have the inclusion
\[
\Delta_\star^Q:=B\big(x_k^\star,\tfrac{a_0\ell(Q)}{2\sqrt M}\big)\cap\partial\Omega_{\m F}\subseteq B\big(x_Q,\tfrac{a_0\ell(Q)}{\sqrt M}\big)\cap\Omega_{\m F}.
\]
In particular, note that $\Gamma\cap\Delta_{\star}^Q\subseteq Q$. It remains only to show that $\Delta_\star^Q$ has the desired properties. To do this, we claim that the inclusion
\begin{equation}\label{eq.bddst}
\partial\Omega_{\m F}\subseteq(\partial\Omega_{\m F}\cap\Gamma)\medcup\Big(\medcup_{Q_j\in\m F}\big(\partial\Omega_{\m F}\cap\overline{R_{Q_j}}\big)\Big)
\end{equation}
holds. Indeed, observe the following elementary set-theoretic calculations:
\[
\bb R^n\backslash\Omega_{\m F}\,=\,\bigcup_{\m W}I^*\backslash\opint\Big(\bigcup_{\m W_{\m F}}I^*\Big)\,\subseteq\,\bigcup_{\m W\backslash\m W_{\m F}}I^*\,=\,\displaystyle\bigcup_{\m R_{Q_j}, Q_j\in\m F}I^*\,\subseteq\,\bigcup_{Q_j\in\m F}\overline{R_{Q_j}},
\]
\begin{equation*}
\partial\Omega_{\m F}\backslash\Gamma\,=\,\partial(\bb R^n\backslash\Omega_{\m F})\backslash\Gamma\,\subset\,\big(\overline{\medcup_{Q_j\in\m F}\overline{R_{Q_j}}}\big)\backslash\Gamma\,=\,\big(\medcup_{Q_j\in\m F}\overline{R_{Q_j}}\big)\backslash\Gamma,
\end{equation*}
where in the last equality we used that the boundary points of $\medcup_{Q_j\in\m F}\overline{R_{Q_j}}$  which are not contained in the union, necessarily lie in $\Gamma$. From these calculations, (\ref{eq.bddst}) follows.

Since $B(x_j^\star, r_j)\supset\overline{R_{Q_j}}$, (\ref{eq.bddst}) holds, and $\Gamma\cap\Delta_{\star}^Q\subset Q$, we will have the desired result as soon as we show that for $Q_j\in\m F$, if $\overline{R_{Q_j}}$ meets $B(x_Q,a_0\ell(Q)/\sqrt M)\supset\Delta_\star^Q$, then $Q_j\subseteq Q$. Thus we show the latter. Suppose that $\overline{R_{Q_j}}\cap B(x_Q,a_0\ell(Q)/\sqrt M)\neq\varnothing$, whence there exists $Q'\in\bb D_{Q_j}$, $I\in\m W_{Q'}$ and $X\in\overline{I^*}$ such that $X\in\overline{I^*}\cap B(x_Q,a_0\ell(Q)/\sqrt M)$. Thus, we note that $\delta(X)\leq|X-x_Q|$, and
\begin{multline*}
\dist(Q',x_Q)\leq\diam Q'+\dist(I,Q')+\diam I^*+|X-x_Q|\\ \leq A_0\ell(Q')+4\sqrt n A_2\ell(Q')+2\sqrt n\ell(I)+\tfrac{a_0\ell(Q)}{\sqrt M}\leq \tfrac{A_0a_2^{-1}}{2\sqrt n}\delta(X)+2A_2a_2^{-1}\delta(X)+\delta(X)+\tfrac{a_0\ell(Q)}{\sqrt M}\\ \leq\tfrac{5A_2a_2^{-1}}{\sqrt M}a_0\ell(Q)<a_0\ell(Q),
\end{multline*}
provided that $M>(5A_2a_2^{-1})^2$. It follows that $Q'\subset\Delta_Q\subset Q$, which implies that $Q_j\cap Q\neq\varnothing$, and so $Q_j\subset Q$ since $Q\in\bb D_{\m F}$. As explained above, this calculation ends the proof.\hfill{$\square$} 

\section{A surface measure on the boundary of a dyadically-generated sawtooth}\label{sec.measure} The goal of this section is to construct a non-negative locally finite Borel measure $\sigma_{\star}$ on $\partial\Omega_{\m F}$ which is doubling and well-suited to work with the elliptic theory of \cite{dfm1}; so that it supplants the role of a ``surface measure'' on the boundary of the sawtooth domain. When $d<n-1$ and $\m F\neq\varnothing$, $\partial\Omega_{\m F}$ is necessarily of mixed dimension. In \cite{dfm20}, an author of this manuscript and other coauthors established an axiomatic elliptic theory for domains with boundaries of mixed dimension. Recall that $m$ is the non-negative Borel measure on $\Omega$ given by $m(E)=\dint_Ew(X)\,dX$, where $w(X)=\delta(X)^{-n+d+1}$. We must construct $\sigma_{\star}$ so that the triple $(\Omega_{\m F},m,\sigma_{\star})$ satisfies the axioms (H1)-(H6) outlined in \cite{dfm20}. 

Our candidate for the measure $\sigma_{\star}$ on $\partial\Omega_{\m F}$ is defined as follows: for each Borel set $E\subset\partial\Omega_{\m F}$, let
\begin{equation}\label{eq.nu}
\sigma_{\star}(E)=\n H^d|_{\Gamma}(E\cap\Gamma)+\int_{E\backslash\Gamma}\dist(X,\Gamma)^{d+1-n}\,\text{d}\n H^{n-1}|_{\partial\Omega_{\m F}\backslash\Gamma}(X).
\end{equation}
We see that $\sigma_{\star}=\sigma+\sigma_{\star}|_{\Sigma}$, where $\sigma=\n H^d|_{\Gamma}$ and $\sigma_{\star}|_{\Sigma}$ are mutually absolutely continuous, and
\[
\sigma_{\star}|_{\Sigma}(E):=\int_{E\backslash\Gamma}\dist(X,\Gamma)^{d+1-n}\,\text{d}\n H^{n-1}|_{\partial\Omega_{\m F}\backslash\Gamma}(X).
\]

\begin{theorem}[Dyadically-generated sawtooth domains admit an elliptic theory]\label{thm.sawtooth} The triple $(\Omega_{\m F}, m, \sigma_{\star})$ satisfies the following axioms.
\begin{enumerate}[(H1)]
\item\label{ax.h1} There exists a constant $c_1>0$ such that for each $x\in\partial\Omega_{\m F}$ and each $r>0$, there exists a point $X\in B(x,r)$ satisfying that $B(X,c_1r)\subset\Omega_{\m F}$.
\item\label{ax.harnack} There exists a positive integer $C_2=N+1$ such that for each $X_1,X_2\in\Omega_{\m F}$ with $\delta_{\star}(X_i)>r$, $i=1,2$, and $|X-Y|\leq7c_1^{-1}r$, there exist $N+1$ points $Z_0:=X_1,$ $Z_2,\ldots, Z_N:=X_2$ in $\Omega_{\m F}$ and verifying $|Z_i-Z_{i+1}|<\delta_{\star}(X)/2$.
\item\label{ax.doublingmeas} The support of $\sigma_{\star}$ is $\partial\Omega_{\m F}$, and $\sigma_{\star}$ is doubling. That is, there exists a constant $C_3>1$ such that for each $x\in\partial\Omega_{\m F}$ and each $r>0$,
\[
\sigma_{\star}(\Delta_{\star}(x,2r))\leq C_3\sigma_{\star}(\Delta_{\star}(x,r)).
\]
\item\label{ax.misgood} The measure $m$ is mutually absolutely continuous with respect to the Lebesgue measure; that is, there exists a weight $\tilde w\in L^1_{\loc}(\Omega_{\m F})$ which is positive Lebesgue-a.e. in $\Omega_{\m F}$, and such that for each Borel set $E\subset\Omega_{\m F}$, we may write $m(E)=\dint_Ew(X)\,dX$. In addition, $m$ is doubling in $\overline{\Omega_{\m F}}$, so that there exists a constant $C_4$ such that for each $X\in\overline{\Omega_{\m F}}$ and each $r>0$, we have that
\[
m(B(X,2r)\cap\Omega_{\m F})\leq C_4 m(B(X,r)\cap\Omega_{\m F}).
\]
\item\label{ax.growth} For each $x\in\partial\Omega_{\m F}$ and each $r>0$, the function $\rho$ given by
\begin{equation}\label{eq.growthrho}\notag
\rho(x,r):=\frac{m(B(x,r)\cap\Omega_{\m F})}{r\sigma_{\star}(\Delta_{\star}(x,r))} 
\end{equation}
verifies for some constant $C_5$ and $\ep:=1/C_5$ that 
\begin{equation}\label{eq.growthrho2}\notag
\frac{\rho(x,r)}{\rho(x,s)}\leq C_5\Big(\frac rs\Big)^{1-\ep},\qquad\text{for each }x\in\partial\Omega_{\m F},\quad 0<s<r.
\end{equation}
\item\label{ax.trace} If $D$ is compactly contained in $\Omega_{\m F}$ and $u_i\in C^{\infty}(\overline D)$ is a sequence of functions such that $\int_D|u_i|\,dm\ra0$ and $\int_D|\nabla u_i-v|^2\,dm\ra0$ as $i\ra\infty$, where $v$ is a vector-valued function in $L^2(D,dm)$, then $v\equiv0$.
\end{enumerate}
\end{theorem}

\noindent\emph{Roadmap to the proof of Theorem \ref{thm.sawtooth}.} The proof of the theorem is split into several parts. First, we check the quantitative properties of $\Omega_{\m F}$, hence in Proposition \ref{prop.h1} below we show that the Corkscrew property \ref{ax.h1} holds, and in Proposition \ref{prop.hcsawtooth} we see that the Harnack Chain property \ref{ax.harnack} holds. We then explore the \dADR-``like'' properties of $\sigma_{\star}$ in Propositions \ref{prop.dadr} and \ref{prop.dadr2}, on which we base our verifications of \ref{ax.doublingmeas} in Proposition \ref{prop.doubling} and of \ref{ax.growth} in Proposition \ref{prop.growth}. Finally, we justify in Remark \ref{rm.done} that \ref{ax.misgood} and \ref{ax.trace} are easy consequences of the previously established results in \cite{dfm1} and the existence of interior Corkscrews for $\partial\Omega_{\m F}$.

As stated, let us show that $\Omega_{\m F}$ enjoys the properties \ref{ax.h1} and \ref{ax.harnack}.

\begin{proposition}[Existence of Corkscrew points for the sawtooth domain]\label{prop.h1} The sawtooth domain $\Omega_{\m F}$ has property \ref{ax.h1}. More precisely, for each $x\in\partial\Omega_{\m F}$ and each $r>0$, there exists a point $X\in B(x,r)$ such that $B(X,c_1r)\subset\Omega_{\m F}\cap B(x,r)$, where $c_1>0$ is the uniform small constant given in (\ref{eq.c1}) below.
\end{proposition}

\noindent\emph{Proof.} Fix $x\in\partial\Omega_{\m F}$ and $r>0$.  Suppose first that $x\in\Gamma\cap\partial\Omega_{\m F}$. Fix $k\in\bb Z$ that satisfies 
\[
\frac r{4A_0c_{\bb K}^{-1}}\leq 2^{-k}<\frac r{2A_0c_{\bb K}^{-1}},
\]
and by Lemma \ref{lm.kcontain}, there exists $Q\in\bb D_{\m F}$ with $2^{-k}\leq\ell(Q)\leq c_{\bb K}^{-1}2^{-k}$ and verifying $x\in\overline Q$. We observe that $B(x_Q,a_0\ell(Q))\subset B(x,r)$: let $Y$ be an arbitrary element of the former, and consider the estimate
\[
|Y-x|\leq|Y-x_Q|+|x_Q-x|\leq a_0\ell(Q)+\diam Q\leq(a_0+A_0)\ell(Q)\leq2A_0c_{\bb K}^{-1}2^{-k}<r,
\]
as desired. Now, according to Corollary \ref{cor.centercork}, there exists $I\in\m W_Q^{\operatorname{cs}}$ such that its center $X_I$ is a Corkscrew point for $Q$ with Corkscrew constant $\tilde c=\frac c{1000\sqrt n}$. Moreover, $B(X_I,\tilde ca_0\ell(Q))\subset\opint(\frac12I)$ and therefore $B(X_I,\tilde ca_0\ell(Q))\subset\Omega_{\m F}$. Reckon that $\tilde ca_0\ell(Q)\geq\frac{\tilde cc_{\bb K}a_0}{4A_0}r=:c_{11}r$, whence the ball $B(X_I,c_{11}r)\subset\Omega_{\m F}\cap B(x,r)$ has the desired properties.

We now consider the case that $x\in\Sigma=\partial\Omega_{\m F}\backslash\Gamma$. In this case, $\delta(x)>0$ and consequently there exists a Whitney box $I\subset\Omega_{\m F}$ such that $x\in\partial I^*$, which we now fix. We split into two cases: either $r\leq 10A_2a_2^{-1}\delta(x)$, or not. 

We resolve the former case first. Since $I^*$ is an $n-$cube, $\opint I^*\subset\Omega_{\m F}$, and $x\in\partial\Omega_{\m F}$, it follows that the ray $R$ containing the line segment $[x,X_I]$ has a non-empty intersection with $B(x,r)$. If $r\leq2\diam I$, then take the unique point $Y\in R$ such that $|Y-x|=r/(4\sqrt n)$. This point satisfies $|Y-x|\leq\ell(I)/2<|X_I-x|$ since $x\in\partial I^*$, and therefore $Y\in[x,X_I]\subset I^*$. Note also that $\dist(Y,\partial I^*)\geq r/(4n)$. Hence the ball $B(Y,r/(8n))=B(Y,c_{12}r)$ has the desired properties. If, instead, $r>2\diam I$, then $I^*\subset B(x,r)$ and it follows that $B(X_I,\ell(I)/4)\subset B(x,r)$. On the other hand, owing to (\ref{eq.bdrystar}) we have that 
\[
\ell(I)\geq\tfrac{\delta(x)}{82\sqrt n}\geq\tfrac{a_2}{820A_2\sqrt n}r.
\]
Thus the ball $B(X_I,a_2r/(4000\sqrt nA_2))=B(X_I,c_{13}r)$ has the desired properties.

It remains only to consider the case that $r>10A_2a_2^{-1}\delta(x)$. In this case, let $Q\in\bb D_{\m F}$ be a dyadic cube such that $I\in\m W_Q$, which we now fix. Observe that $\ell(Q)\leq\frac r{20\sqrt n A_2}$, and that for any generation $k\leq k(Q)$, there is a unique $Q_k\in\bb D_{\m F}\cap\bb D^k$ which contains $Q$. Let $k$ be the unique generation such that
\[
\frac r{20\sqrt nA_2}\leq 2^{-k}<\frac r{10\sqrt n A_2},
\]
and choose $Q_k$ as above. According to Corollary \ref{cor.centercork}, there exists a point $X_k\in\Omega$ which is the center of some Whitney box $I_k\in\m W_{Q_k}^{\operatorname{cs}}$ and is also a Corkscrew point for $Q_k$ with Corkscrew constant $\tilde c$. Let us see that with our choice of constants, we have that $B(X_k,\tilde ca_0\ell(Q_k)\subset B(x,r)$. Fix $Y\in B(X_k,\tilde ca_0\ell(Q_k))\subset\frac12 I$, and consider the estimate
\begin{multline*}
|Y-x|\leq\diam I_k+\dist(I_k,Q_k)+\diam Q_k+\dist(I^*,Q)+\diam I^*\\ \leq\sqrt\ell(I_k)+a_0\ell(Q_k)+A_0\ell(Q_k)+4\sqrt nA_2\ell(Q)+2\diam I\\ \leq (5a_0/4+A_0)\ell(Q_k)+6\sqrt nA_2\ell(Q) \leq9\sqrt n A_2\ell(Q_k)<10\sqrt n A_2\ell(Q_k)<r,
\end{multline*}
as desired. Note that $\tilde ca_0\ell(Q_k)\geq\tilde ca_0\frac r{20\sqrt nA_2}$, so that $B(X_k,\frac{\tilde ca_0}{20\sqrt n A_2}r)=:B(X_k,c_{14}r)$ is a ball with the desired properties.

Finally, take
\begin{multline}\label{eq.c1}
c_1=\min\big\{c_{11},c_{12},c_{13},c_{14}\big\}\\ =\min\Big\{\frac{\tilde cc_{\bb K}a_0}{4A_0}~,~\frac1{8n}~,~\frac{a_2}{4000\sqrt n A_2}~,~ \frac{\tilde ca_0}{20\sqrt n A_2}\Big\}
\end{multline}
and reckon that the desired result is proved.\hfill{$\square$}

\begin{proposition}[Harnack Chains in the sawtooth domain]\label{prop.hcsawtooth} The sawtooth domain $\Omega_{\m F}$ has property \ref{ax.harnack}. More precisely, for any $\Lambda\geq1$, there exists  $C_2=C_2(n,d,C_d,\theta,\Lambda)$ such that if $X_1,X_2\in\Omega_{\m F}$ are two points with $\min_i\delta_\star(X_i)\geq r$ for $i=1,2$, and $|X_1-X_2|\leq\Lambda r$, then there is a chain of balls $\{B_m\}_m\subset\Omega_{\m F}$  connecting $X_1$ to $X_2$, and verifying that $\card(\{B_m\}_m)\leq C_2$, and $\dist(B_m,\partial\Omega_{\m F})\approx\diam B_m$, for each  $B_m$.
\end{proposition}

We call a chain of balls $\{B_m\}_m$ as in the above proposition a \emph{Harnack Chain} (for the sawtooth domain).

\noindent\emph{Proof.} Fix $X_1,X_2\in\Omega_{\m F}$ with $\delta_\star(X_i)>r$, $i=1,2$, and $|X_1-X_2|\leq\Lambda r$. We seek to join $X_1$ and $X_2$ via a Harnack Chain that stays far from $\partial\Omega_{\m F}$. First, note that for any $Z\in\Omega_{\m F}$, $\delta(Z)\geq\delta_\star(Z)$, while in the other direction we have that if $Z\in I$ for some $I\in\m W_{\m F}$, then
\[
\delta_\star(Z)\geq\tfrac12\theta\ell(I)\geq\tfrac1{164\sqrt n}\theta\delta(Z).
\]
For each $i=1,2$, fix $Q_i\in\bb D_{\m F}$ and $I_i\in\m W_{Q_i}$ such that $X_i\in\opint I_i^*$. If it can be arranged that $I_1=I_2$, then the result follows immediately from the fact that $\opint I^*$ is an $n-$dimensional open cube. Similarly, if $I_1^*\cap I_2^*\neq\varnothing$, then a Harnack Chain connecting $X_1$ to $X_2$ can be obtained by noting that $I_1^*\cap I_2^*$ is a union of rectangles with no side-length smaller than $\theta\min_i\ell(I_i)$, whence we may use these intersections to ``transfer'' from $X$ to $Y$ in the manner desired. Therefore, if the estimate $\delta_\star(X_1)<\frac{\theta}{200\Lambda}\delta(X_1)=:\frac1M\delta(X_1)$, holds, then 
\[
|X_1-X_2|\leq\Lambda r\leq\tfrac{\theta}{200}\delta(X_1)\leq\tfrac{\theta}2\diam I_1\ll\diam I_1,
\]
so that $I_1$ and $I_2$ touch and by the above discussion, the desired result is achieved with $C_2=C_2(n,\Lambda)$. 

Thus it remains to obtain the conclusion under the supposition that for each $i=1,2$, the estimate $\delta_\star(X_i)\geq\frac1M\delta(X_i)$ holds. In this case, we have that $\delta_\star(X_i)\approx_M\delta(X_i)$, and without loss of generality suppose that $\ell(Q_1)\leq\ell(Q_2)$. In this case, we may connect each $X_i$ to the respective centers of $I_i$, $X_{I_i}$, through Harnack Chains with a uniform number of balls (depending only on $M$). Hence we have reduced the problem to procuring a Harnack Chain between $X_{I_1}$ and $X_{I_2}$. We have that $Q_1$ and $Q_2$ have comparable length, as follows: first, we have the estimate
\begin{equation*}
\delta(X_2)\leq|X_1-X_2|+\delta(X_1)\leq\Lambda r+\delta(X_1)\leq2\Lambda\delta(X_1),
\end{equation*}
which gives that $\ell(I_2)\leq\tfrac{41}2\Lambda\ell(I_1)$, and on the other hand,
\[
n^{-1/2}A_2^{-1}\delta(X_i)/41\leq A_2^{-1}\ell(I_i)\leq\ell(Q_i)\leq a_2^{-1}\ell(I_i)\leq n^{-1/2}a_2^{-1}\delta(X_i)/4,\qquad i=1,2,
\]
which implies that
\[
\delta(X_1)\leq\tfrac{41}4A_2a_2^{-1}\delta(X_2),\qquad\text{and}\qquad \ell(I_1)\leq A_2a_2^{-1}\ell(I_2).
\]
Consequently,
\[
\ell(Q_1)\leq\ell(Q_2)\leq21\Lambda A_2a_2^{-1}\ell(Q_1),
\]
and, furthermore,
\begin{multline}\label{eq.qdist}
\dist(Q_1,Q_2)\leq\dist(I_1,Q_1)+\diam I_1+|X_1-X_2|+\diam I_2+\dist(I_2,Q_2)\\ \leq5\sqrt nA_2(\ell(Q_1)+\ell(Q_2))+\Lambda r \leq 51\sqrt n A_2\Lambda\ell(Q_2).
\end{multline}
Fix $Q_2^a\in\bb D_{\m F}$ as the unique ancestor of $Q_2$ that satisfies
\[
\frac{51\sqrt nA_2\Lambda}{500}\ell(Q_2)\leq\ell(Q_2^a)\leq c_{\bb K}^{-1}\frac{51\sqrt nA_2\Lambda}{500}\ell(Q_2),
\]
and then choose for $Q_1^a\in\bb D_{\m F}$ the ancestor of $Q_1$ that verifies $c_{\bb K}\ell(Q_2^a)\leq\ell(Q_1^a)\leq\ell(Q_2^a)$. By construction, we have that $\ell(Q_1^a)\approx\ell(Q_2^a)\approx\ell(Q_2)\approx\ell(Q_1)$ with uniform constants. By Lemma \ref{lm.whitneychild}, we see that $Q_2$ and its proper parent $Q_2'$ satisfy $\m W_{Q_2}^0\cap\m W_{Q_2'}^0\neq\varnothing$, so that by the construction of $\m W_{Q_2}$ in (\ref{eq.whitneyboxes}), there exists a Harnack Chain of the desired properties connecting $X_{I_2}$ to some point $X_2'$ lying in $I_2'\in\m W_{Q_2'}$. It is easy to see that therefore  we may inductively ``ascend'' through a uniformly finite (since $\ell(Q_2^a)\approx\ell(Q_2)$) sequence of Harnack Chains from $X_{I_2}$ to a point $X_{I_2^a}$ which is the center of a Whitney cube $I_2^a\in\m W_{Q_2^a}$. Now, from (\ref{eq.qdist}), we see that
\[
\dist(Q_1^a,Q_2^a)\leq\dist(Q_1,Q_2)\leq500\ell(Q_2^a),
\]
so that by Lemma \ref{lm.closeregions}, $\m W_{Q_1^a}^0\cap\m W_{Q_2^a}^0\neq\varnothing$. Hence we may pass through a Harnack Chain from $X_{I_2^a}$ to a point $X_{I_1^a}$ which is the center of some $I_1^a\in\m W_{Q_1^a}$. As before but in reverse, we proceed to ``descend'' from $X_{I_1^a}$ to $X_{I_1}$ through a uniformly finite (since $\ell(Q_1^a)\approx\ell(Q_1)$) sequence of Harnack Chains. Hence, in this case the desired result is achieved with a constant $C_2=C_2(n,d,C_d,\Lambda,\theta)$.\hfill{$\square$}

We turn to a study of the properties of the measure $\sigma_{\star}$.

\begin{lemma}[Support of $\sigma_{\star}$]\label{lm.support} The measure $\sigma_{\star}$ is supported on $\partial\Omega_{\m F}$.
\end{lemma}

\noindent\emph{Proof.} It is clear that $\Gamma\cap\partial\Omega_{\m F}\subset\supp\sigma_{\star}$, so we only need to check that $\Sigma$ is in the support of $\sigma_{\star}$. But this is easy: for any bounded open set $U$ intersecting $\Sigma$ and compactly contained in $\bb R^n\backslash\Gamma$, the set $\Sigma\cap U$ is a finite union of non-empty $(n-1)-$dimensional rectangles, so that $\n H^{n-1}(\Sigma\cap U)\in(0,\infty)$, and $\delta(X)\in(0,+\infty)$ for any $X\in\Sigma$. The claim ensues.\hfill{$\square$}

\begin{proposition}[Upper bound for $\sigma_{\star}$]\label{prop.dadr} Let $x\in\partial\Omega_{\m F}$ and $r>0$. Then
\begin{equation}\label{eq.close}
\sigma_{\star}(\Delta_{\star}(x,r))\leq V_1r^d,
\end{equation}
where $V_1=V_1(n,d,C_d,a_0,A_0,\zeta,c,c_{\m H})$. Moreover, if $\delta(x)>r/M_0$ for $M_0=125A_0^2c_{\bb K}^{-2}$, then
\begin{equation}\label{eq.far}
\sigma_{\star}(\Delta_{\star}(x,r))\leq V_2\delta(x)^{d+1-n}r^{n-1}.
\end{equation}
The uniform constant $V_2$ in the last inequality depends only on $n$, $d$, and $V_1$.
\end{proposition}

\noindent\emph{Proof.} Fix $x\in\partial\Omega_{\m F}$ and $r>0$. Let $B:=B(x,r)$ and recall that $\Sigma=\partial\Omega_{\m F}\backslash\Gamma$, $\sigma=\n H^d|_{\Gamma}$. We first prove (\ref{eq.close}) by adapting ideas of the proof for Lemma 3.61 from \cite{hm2}. Observe that
\begin{multline}\label{eq.twoparts}\notag
\sigma_{\star}(\Delta_{\star}(x,r))=\sigma_{\star}(B\cap\Gamma\cap\partial\Omega_{\m F})+\sigma_{\star}(B\cap\Sigma)\\ \leq\sigma_{\star}(B\cap\Gamma)+\sigma_{\star}(B\cap\Sigma)=\sigma(\Delta(x,r))+\sigma_{\star}|_{\Sigma}(B\cap\Sigma)\\ \leq C_dr^d+\sigma_{\star}|_{\Sigma}(B\cap\Sigma).
\end{multline}
Thus we need only show that 
\begin{equation}\label{eq.secondpart}\notag
\sigma_{\star}|_{\Sigma}(B(x,r)\cap\Sigma)\lesssim r^d.
\end{equation}
We will do this by splitting $B\cap\Sigma$ into two parts: one part where the (portions of) faces in $B\cap\Sigma$ correspond to Whitney boxes having small diameter compared to $r$, and the other part  where the (portions of) faces in $B\cap\Sigma$ correspond to Whitney boxes having large diameter compared to $r$. More precisely, it is clear that any $X\in\Sigma$ lies in the face of a fattened Whitney box $J^*$, such that $J\in\m W$, $\opint J^*\subset\Omega_{\m F}$, and $\partial J^*\cap\partial\Omega_{\m F}\neq\varnothing$. Then there exists a Whitney box $I\in\m W$, with $I\notin\m W_Q$ for any $Q\in\bb D_{\m F}$, so that $J^*\cap I\neq\varnothing$ (otherwise, every Whitney box adjacent to $J$ lies in $\m W_{\m F}$, contradicting that $J^*\cap\partial\Omega_{\m F}\neq\varnothing$). Necessarily then, there exists $Q'\in\bb D_{Q_j}$ with $Q_j\in\m F$ and verifying that $I\in\m W_{Q'}$. Denote by $\m F_B$ the sub-collection of those $Q_j\in\m F$ such that there exists $I\in\m R_{Q_j}$ (cf. (\ref{eq.carlesonfam})) intersecting $B\cap\Sigma$.  Let $\m F_B=\m F_1\cup\m F_2$ where $Q_j\in\m F_B$ belongs to $\m F_1$ if $\ell(Q_j)<r$, and $\m F_2=\m F_B\backslash\m F_1$. Then, we may write
\begin{multline}\label{eq.split}
B\cap\Sigma=B\cap\Sigma\cap\big(\bigcup_{Q_j\in\m F_{\m B}}\bigcup_{I\in\m R_{Q_j}}I\big)\\=\Big(B\cap\Sigma\cap\big(\bigcup_{Q_j\in\m F_{1}}\bigcup_{I\in\m R_{Q_j}}I\big)\Big)\bigcup\Big(B\cap\Sigma\cap\big(\bigcup_{Q_j\in\m F_{2}}\bigcup_{I\in\m R_{Q_j}}I\big)\Big)\\=\big(\bigcup_{Q_j\in\m F_1}(B\cap\Sigma_j)\big)\medcup\displaystyle\big(\bigcup_{Q_j\in\m F_2}(B\cap\Sigma_j)\big),
\end{multline}
where $\Sigma_j:=\Sigma\cap(\cup_{I\in\m R_{Q_j}}I)$ for each $Q_j\in\m F$. Our further analysis will be based on the following estimate:
\begin{equation}\label{eq.sigmaj}
\sigma_{\star}(B\cap\Sigma_j)\lesssim\big(\min\big\{r\,,\,\ell(Q_j)\big\}\big)^d,\qquad\text{for each }Q_j\in\m F.
\end{equation}
Suppose momentarily that (\ref{eq.sigmaj}) holds, and we will use it to control $\sigma_{\star}(B\cap\Sigma)$. First, we consider the contribution of $\m F_1$. If $Q_j\in\m F_1$ so that $\ell(Q_j)<r$, we have that $Q_j\in B^*:=B(x,(4a_0K+3A_0)r)$. Indeed, since $Q_j\in\m F_1\subset\m F_B$, it follows that there exists $Q'\in\bb D_{Q_j}$ and $I\in\m W_{Q'}$ such that $B\cap I\neq\varnothing$, and thus for any $q\in Q_j$, we have that
\begin{multline*}
|q-x|\leq\diam Q_j+\dist(Q_j,x)\leq A_0\ell(Q_j)+\dist(I,Q')+\diam I+\dist(I,x)\\ \leq A_0\ell(Q_j)+(2a_0K+A_0)\ell(Q')+\tfrac{2a_0K+A_0}{4\sqrt n}\ell(Q')+r~<(4a_0K+3A_0)r,
\end{multline*}
where we used that $\ell(Q')\leq\ell(Q_j)<r$. Then, 
\begin{multline}\label{eq.f1control}
\sigma_{\star}\big(\bigcup_{Q_j\in\m F_1}(B\cap\Sigma_j)\big)\leq\sum_{Q_j\in\m F_1}\sigma_{\star}(B\cap\Sigma_j)\leq C\sum_{Q_j\in\m F_1}\ell(Q_j)^d\leq CC_da_0^{-d}\sum_{Q_j\in\m F_1}\sigma(Q_j)\\ \leq CC_da_0^{-d}\sigma\big(B^*\cap\Gamma)\leq CC_d^2\big(4K+3\tfrac{A_0}{a_0}\big)^dr^d,
\end{multline}
where $C$ is the uniform constant implicit in (\ref{eq.sigmaj}), and in the second inequality we used (\ref{eq.sigmaj}), in the third inequality we used (\ref{eq.measureofQ}), in the fourth inequality we used that the cubes $Q_j\subset B^*$ are disjoint.

Next we turn to the contribution of $\m F_2$, still supposing that (\ref{eq.sigmaj}) holds. We begin by proving that the cardinality of $\m F_2$ is uniformly bounded. Suppose that $Q_j$ and $Q_k$ belong to $\m F_2$, so that there exist Whitney boxes $I_j\in\m R_{Q_j}, I_k\in\m R_{Q_k}$ intersecting $B\cap\Sigma$, and without loss of generality we may assume that $\ell(Q_k)\leq\ell(Q_j)$. Observe that for $i=j,k$ we have that $\ell(I_i)\leq\frac{2a_0K+A_0}{4\sqrt n}\ell(Q_i)$ since $I_i\in\m R_{Q_i}$. Also note that $\dist(I_j,I_k)\leq\diam B=2r\leq2\ell(Q_k)$, and moreover
\begin{multline*}
4\diam I_j\leq\dist(I_j,\Gamma)\leq\dist(I_j,I_k)+\diam I_k+\dist(I_k,\Gamma)\\ \leq2\ell(Q_k)+\tfrac{41}4(2a_0K+A_0)\ell(Q_k)\leq A(a_0K+A_0)\ell(Q_k),
\end{multline*}
where $A\geq1$ is a large real number with no dependence on any parameter. Now, we see that
\begin{multline*} 
\dist(Q_j,Q_k)\leq\diam Q_k+\dist(I_k,Q_k)+\diam I_k+\dist(I_k,I_j)+\diam I_j+\dist(I_j,Q_j)\\ \leq A(a_0K+A_0)a_1^{-1}\eta^{-\frac{n-1}{n-1-d}}\ell(Q_k),
\end{multline*}
where $a_1=a_1(n,d,C_d,c,c_{\m H},K)$ is the quantity defined in (\ref{eq.xi}). Thus, we have shown that for any $Q_j,Q_k\in\m F_2$, the estimate
\begin{equation}\label{eq.smalldist}
\dist(Q_j,Q_k)\leq A(a_0K+A_0)a_1^{-1}\eta^{-\frac{n-1}{n-1-d}}\min\big\{\ell(Q_j),\ell(Q_k)\big\}=:A_1\min\big\{\ell(Q_j),\ell(Q_k)\big\}
\end{equation}
holds. Let us see that (\ref{eq.smalldist}) implies the uniform boundedness of $\card\m F_2$. Since for all $Q_k\in\m F_2$ we have that $\ell(Q_k)\geq r$ by definition, then we may choose $Q_j\in\m F_2$ so that $\ell(Q_k)\geq\ell(Q_j)$ for all $Q_k\in\m F_2$. Fix such a $Q_j\in\m F_2$, and reckon that by (\ref{eq.smalldist}) and (\ref{eq.centerofQ}), for each $Q_k\in\m F_2$ the set $Q_k\cap\Delta(x_{Q_j},(A_0+A_1)\ell(Q_j))$ is not empty. Accordingly, for each $Q_k\in\m F_2$, there exists a dyadic cube $Q_k'\in\bb D_{Q_k}$ such that $c_{\bb K}\ell(Q_j)\leq\ell(Q_k')\leq\ell(Q_j)$ and $Q_k'\subset\Delta(x_{Q_j},3A_1\ell(Q_j))$. We consider the estimate
\begin{multline*}
c_{\bb K}C_d^{-1}a_0^d\ell(Q_j)^d\card\m F_2\leq\sum_{Q_k\in\m F_2}\sigma\big(\Delta(x_{Q_k'},a_0\ell(Q_k'))\big)\leq\sum_{Q_k\in\m F_2}\sigma(Q_k')\\ \leq\sigma(\Delta(x_{Q_j},3A_1\ell(Q_j))\leq C_d(3A_1)^d\ell(Q_j)^d,
\end{multline*}
and hence obtain that
\[
\card\m F_2\leq c_{\bb K}^{-1}C_d^2\big[3a_0^{-1}A_1\big]^d.
\]
Therefore, we may conclude, using (\ref{eq.sigmaj}), that
\begin{multline}\label{eq.f2control}
\sigma_{\star}\big(\medcup_{Q_j\in\m F_2}\displaystyle(B\cap\Sigma_j)\big)\leq\sum_{Q_j\in\m F_2}\sigma_{\star}(B\cap\Sigma_j)\leq C_d^2\big[3a_0^{-1}A_1\big]^d\sup_{Q_j\in\m F_2}\sigma_{\star}(B\cap\Sigma_j)\\ \leq C_d^2\big[3a_0^{-1} A_1 \big]^d Cr^d.
\end{multline}
Putting (\ref{eq.split}), (\ref{eq.f1control}), and (\ref{eq.f2control}) together, we obtain the desired result modulo the proof of (\ref{eq.sigmaj}). 

We now turn to the proof of (\ref{eq.sigmaj}). Hence take $Q_j\in\m F$ and first suppose that $\ell(Q_j)\leq Mr$ for some $M>0$ to be fixed later. In this case, any $I\in\m R_{Q_j}$ satisfies $\ell(I)\lesssim\ell(Q_j)$, but this estimate is too crude as there may be too many (in fact, infinitely many!) such boxes intersecting $\Sigma_j$. Therefore the idea is to control the number of Whitney boxes $I$ of a given generation that contribute to $\Sigma_j$. To this end, recall $A_2=\frac{2a_0K+A_0}{4\sqrt n}$, and define
\[
\Sigma_j^k:=\Sigma\textstyle\bigcap\big(\textstyle\bigcup_{\tiny\{ I\in\m R_{Q_j}:\,\ell(I)=2^{-k}\}}I\big),\qquad\text{so that}\quad\Sigma_j=\bigcup_{\{k\,:2^{-k}\leq A_2\ell(Q_j)\}}\Sigma_j^k.
\] 
Hence,
\begin{multline}\label{eq.sums}
\sigma_{\star}(B\cap\Sigma_j)=\sum_{k:2^{-k}\leq A_2\ell(Q)}\sigma_{\star}(B\cap\Sigma_j^k)\\ =\sum_{k:a_3\ell(Q_j)<2^{-k}\leq A_2\ell(Q)}\sigma_{\star}(B\cap\Sigma_j^k)+\sum_{k:2^{-k}<a_3\ell(Q_j)}\sigma_{\star}(B\cap\Sigma_j^k)=:T_1+T_2,
\end{multline}
where $a_3>0$ is a small constant to be fixed momentarily. We bound term $T_1$ first. Note that if $X\in\Sigma_j^k$, then there exists $I\in\m R_{Q_j}$ with $\ell(I)>a_3\ell(Q_j)$ and such that $X\in I$. Hence $\delta(X)\geq4\sqrt na_3\ell(Q_j)$. For convenience in the following calculation, set $\m W_{\Sigma_j^1}:=\{I\in\m R_{Q_j}: a_3\ell(Q_j)\leq\ell(I)\leq A_2\ell(Q_j)\}$, and observe that
\begin{multline}\label{eq.adrcalc1}
\sum_{k:a_3\ell(Q_j)<2^{-k}\leq A_2\ell(Q)}\sigma_{\star}(B\cap\Sigma_j^k)\leq\sum_{k:a_3\ell(Q_j)<2^{-k}\leq A_2\ell(Q)}\int_{\Sigma_j^k}\delta(X)^{d+1-n}\dH(X)\\ \leq(4\sqrt na_3)^{d+1-n}\ell(Q_j)^{d+1-n}\n H^{n-1}\big(\textstyle\bigcup_{I\in\m W_{\Sigma_j^1}}(I\cap\Sigma)\big)\\ \leq a_3^{d+1-n} \ell(Q_j)^{d+1-n}\sum_{I\in\m W_{\Sigma_j^1}}\n H^{n-1}(I\cap\Sigma)\\ \leq A_na_3^{d+1-n}\big[\tfrac{a_0K+A_0}{a_3}\big]^n\ell(Q_j)^{d+1-n}\big(\sup_{I\in\m W_{\Sigma_j^1}}\n H^{n-1}(I\cap\Sigma)\big),
\end{multline}
where by $A_n\geq1$ we denote a universal constant depending only on $n$, and in the last line we used the bound (\ref{eq.cardgenw}), since the set $\m W_{\Sigma_j^1}$ can easily be seen to be a subset of a set of the form in Lemma \ref{lm.numbergenw}. Next, let $I\in\m W_{\Sigma_j^1}$, and we seek to bound $\n H^{n-1}(I\cap\Sigma)$. Observe that
\begin{multline}\label{eq.onei}
\n H^{n-1}(I\cap\Sigma)\leq\n H^{n-1}\big(\textstyle\bigcup_{J\in\m W_{\m F}: J^*\cap I\neq\varnothing}\partial J^*\big)\leq A_n\big(\sup_{J\in\m W_{\m F}: J^*\cap I\neq\varnothing}\n H^{n-1}(\partial J^*)\big)\\ \leq A_n(1+\theta)^{n-1}\big(\sup_{J\in\m W_{\m F}: J^*\cap I\neq\varnothing}\ell(J)^{n-1}\big) \leq A_n\ell(I)^{n-1}\leq A_nA_2^{n-1}\ell(Q_j)^{n-1}.
\end{multline}
We may combine (\ref{eq.adrcalc1}) with (\ref{eq.onei}) to see that
\begin{equation}\label{eq.sum1}
T_1\leq A_na_3^{d+1-2n}\big[a_0K+A_0\big]^{2n-1}\ell(Q_j)^d,
\end{equation}
which is the desired bound for $T_1$. We remark that the estimates in (\ref{eq.onei}) and (\ref{eq.sum1}) also allow us to say that for any $I\in\m W$, it holds that
\begin{equation}\label{eq.forallI}
\sigma_{\star}(I\cap\Sigma)\leq A_n\ell(I)^d.
\end{equation}

Now we seek to bound $T_2$. Set $\m W_{\Sigma_j^k}:=\{I\in\m R_{Q_j}:\ell(I)=2^{-k},~ I\cap\Sigma_j^k\neq\varnothing\}$,  and observe that
\begin{multline}\label{eq.t2control}
T_2\leq\sum_{k:2^{-k}<a_3\ell(Q_j)}~\sum_{I\in\m R_{Q_j}:\ell(I)=2^{-k}}\sigma_{\star}(I\cap\Sigma)\\ \leq\sum_{k:2^{-k}<a_3\ell(Q_j)}~\card(\m W_{\Sigma_j^k})\big(\sup_{I\in\m W_{\Sigma_j^k}}\sigma_{\star}(I\cap\Sigma)\big)\leq A_n\sum_{k:2^{-k}<a_3\ell(Q_j)}~\card(\m W_{\Sigma_j^k})2^{-kd},
\end{multline}
where we used (\ref{eq.forallI}) in the last inequality. Hence, it will suffice to show that for some large uniform constant $C$ and some $\zeta\in(0,1)$, the estimate
\begin{equation}\label{eq.t21}
\card(\m W_{\Sigma_j^k})\leq C(2^k\ell(Q_j)^{d-\zeta}
\end{equation}
holds. Let us establish (\ref{eq.t21}) then. If  $I\in\m W_{\Sigma_j^k}$, then there exists $Q_I\in\bb D_{Q_j}$ such that $I\in\m W_{Q_I}$, and moreover there exists $J\in\m W_{\m F}$ so that $I\cap J^*\neq\varnothing$. In particular, there exists $Q'\in\bb D_{\m F}$ such that $J\in\m W_{Q'}$. Observe that
\[
\ell(Q')\leq a_2^{-1}\ell(J)\leq4a_2^{-1}\ell(I)\leq 4a_2^{-1}a_3\ell(Q_j),
\]
whence $\ell(Q')<\ell(Q_j)$ provided that $a_3\leq a_2/8$, which we assume from now on. Since $Q'\in\bb D_{\m F}$ and $Q_j\in\m F$, it follows that $Q'$ and $Q_j$ are disjoint. Consequently,
\begin{multline}\label{eq.t22}
\dist(Q_I,\Gamma\backslash Q_j)\leq\dist(Q_I,Q')\leq\diam Q_I+\dist(I,Q_I)+\diam I+\dist(J,Q')\\ \leq (A_0a_2^{-1}+6\sqrt nA_2a_2^{-1})\ell(I).
\end{multline}
Since for any $q\in Q_I$ we have that $\dist(x,\Gamma\backslash Q_j)\leq\diam Q_I+\dist(Q_I,\Gamma\backslash Q_j)$, then by using (\ref{eq.t22}), it follows that
\begin{equation*}
Q_I\subset\big\{x\in Q_j~:~\dist(x,\Gamma\backslash Q_j)\leq\big[10\sqrt na_2^{-1}A_22^{-k}\ell(Q_j)^{-1}\big]\ell(Q_j)\big\}=:V.
\end{equation*}
Observe that $V$ is the set considered in property \ref{item.bdrythin} of Lemma \ref{lm.dyadiccubes}, with \break $\rho=10\sqrt na_2^{-1}A_22^{-k}\ell(Q_j)^{-1}$. We may apply the inequality in \ref{item.bdrythin} so long as $\rho<a_0$, which in our case will be true as long as $a_3\leq\frac{a_0a_2}{100\sqrt nA_2}$. Henceforth we fix $a_3$ to be given by the right-hand side of the last inequality, and note that it also satisfies $a_3\leq a_2/8$. Then, we have that
\begin{equation}\label{eq.t23}
\sigma\big(\textstyle\bigcup_{I\in\m W_{\Sigma_j^k}}Q_I\big)\leq\sigma(V)\leq A_0(10\sqrt na_2^{-1}A_2)^{\zeta}2^{-k\zeta}\ell(Q_j)^{-\zeta}\sigma(Q_j).
\end{equation}
Now, we reckon the estimate
\begin{multline}\label{eq.t24}
\card(\m W_{\Sigma_j^k})2^{-kd}=\sum_{I\in\m W_{\Sigma_j^k}}\ell(I)^d\leq A_2^d\sum_{I\in\m W_{\Sigma_j^k}}\ell(Q_I)^d\leq C_d(a_0^{-1}A_2)^d\sum_{I\in\m W_{\Sigma_j^k}}\sigma(Q_I)\\ \leq C_d(a_0^{-1}A_2)^d\sum_{\tiny\begin{matrix}\tilde Q\in\bb D_{Q_j}\text{ s.t. }\\\exists I\in\m W_{\Sigma_j^k}\cap\m W_{\tilde Q}\end{matrix}}\sum_{I\in\m W_{\tilde Q}}\sigma(\tilde Q)\leq C_d(a_0^{-1}A_2)^dN_0\sum_{\tiny\begin{matrix}\tilde Q\in\bb D_{Q_j}\text{ s.t. }\\\exists I\in\m W_{\Sigma_j^k}\cap\m W_{\tilde Q}\end{matrix}}\sigma(\tilde Q)\\ \leq C_d(a_0^{-1}A_2)^dN_0\big[C_d^2(2A_0a_0^{-1})^d\big]^{\frac{A_2}{a_2}}\sigma\big(\textstyle\bigcup\limits_{\tiny\begin{matrix}\tilde Q\in\bb D_{Q_j}\text{ s.t. }\\\exists I\in\m W_{\Sigma_j^k}\cap\m W_{\tilde Q}\end{matrix}}\tilde Q\big)\leq A_32^{-k\zeta}\ell(Q_j)^{d-\zeta},
\end{multline}
where $A_3$ is a uniform constant. In the fifth inequality, we used Corollary \ref{cor.numberwq}, in the sixth inequality we used (\ref{eq.cardchildren}) and the fact that $\ell(Q_I)/\ell(Q_{I'})\leq A_2/a_2$ for any $I,I'\in\m W_{\Sigma_j^k}$, and in the last inequality we used (\ref{eq.t23}). It is clear that (\ref{eq.t24}) gives (\ref{eq.t21}). Going back to (\ref{eq.t2control}), we can now conclude that
\begin{equation*}
T_2\leq A_nA_3\ell(Q_j)^{d-\zeta}\sum_{k:2^{-k}<a_3\ell(Q_j)}2^{-k\zeta}\leq\tfrac{2^{\zeta}}{2^{\zeta}-1}a_3^{\zeta}A_nA_3\ell(Q_j)^d.
\end{equation*}
Putting this last estimate together with (\ref{eq.sum1}) and (\ref{eq.sums}) gives (\ref{eq.sigmaj}) in the case that $\ell(Q_j)\leq Mr$.

It remains only to show that (\ref{eq.sigmaj}) holds in the case that $\ell(Q_j)>Mr$. Observe that if $x\in\Sigma$, then the ball $B$ is centered on an $(n-1)-$dimensional face of some Whitney box $J^*$, with $\opint J^*\subset\Omega_{\m F}$. Suppose that $\ell(J)\geq A_n' r$, where $A_n'\geq1$ is chosen so that $B\cap\Sigma$ is a subset of the boundary faces of the Whitney cubes adjacent to $J$ (including also $J$). It is clear that this is a constraint solely depending on $n$. In this case, if $X\in\partial I^*$ for any $I$ touching $J$, then $\delta(X)\geq\dist(I^*,\Gamma)\geq2\sqrt n\ell(I)\geq\sqrt n\ell(J)/2\geq r/2$, and we have the following bound:
\begin{multline*}
\sigma_{\star}(B\cap\Sigma_j)\leq\sigma_{\star}(B\cap\Sigma)\leq\sigma_{\star}\big(\medcup_{I\text{ touching }J}(\partial I^*\cap B)\big)\leq\textstyle\sum_{I\text{ touching }J}\sigma_{\star}(\partial I^*\cap B)\\ \leq A_n\sup_{I\text{ touching }J}\sigma_{\star}(\partial I^*\cap B)\leq A_n\sup_{I\text{ touching }J}\int_{\partial I^*\cap B}\delta(X)^{d+1-n}\dH(X)\\ \leq A_n2^{n-1-d}r^{d+1-n}\sup_{I\text{ touching }J}\n H^{n-1}(B\cap\partial I^*)\leq A_n2^{n-1-d}r^d,
\end{multline*}
where $A_n$ is a universal constant depending only on $n$, in the fourth inequality we used that the number of Whitney boxes adjacent to $J$ is uniformly bounded (depending only on $n$), and in the last line we used the facts that at least one of $\n H^{n-1}(B\cap\partial I^*)>0$, that for any such $I$ there exists $x'\in\partial I^*$ satisfying $B\cap\partial I^*\subset B(x_I,2r)\cap\partial I^*$, and that  each $\partial I^*$ is an $(n-1)-$Ahlfors-David regular set.

Now suppose that either $x\in\Gamma$, or $x\in\Sigma$ with $\ell(J)\leq A_n'r$. The bound $\delta(x)\leq42\sqrt nA_n'r$ holds trivially in the former case, and in the latter it holds because of the estimate
\[
\delta(x)\leq\dist(J^*,\Gamma)+\diam J^*\leq42\sqrt n\ell(J).
\]
For each $I\in\m R_{Q_j}$ intersecting $B$, there exists $Q_I\in\bb D_{Q_j}$ such that $I\in\m W_{Q_I}$, and we have that
\[
4\sqrt n\ell(I)\leq\dist(I,\Gamma)\leq\dist(I,x)+\delta(x)\leq 43\sqrt nA_n'r.
\]
Hence $\ell(I)\leq11\sqrt nA_n'r$, and for any $q_I\in Q_I$, we reckon that
\begin{multline}\notag
|q_I-x|\leq\diam Q_I+\dist(Q_I,I)+\diam I+\dist(I,x)\leq A_0\ell(Q_I)+A_2\ell(Q_I)+\sqrt n\ell(I)+r\\ \leq (2A_2a_2^{-1}+\sqrt n)\ell(I)+r\leq100\sqrt nA_n'A_2a_2^{-1}r=:A_4r.
\end{multline}
Thus $Q_I\subset B(x,A_4r)$. Now let $\{Q^i\}\subset\bb D_{Q_j}$ be a covering of $B(x,A_4r)\cap Q_j$ such that $Mr/2\leq\ell(Q^i)\leq Mr$ (which is possible since $\ell(Q_j)>Mr$) and such that the $Q^i$ are pairwise disjoint. It is easy to see then that $\dist(Q^{i_1},Q^{i_2})\leq(4A_4+A_0)\min\{\ell(Q^{i_1}),\ell(Q^{i_2})\}$, whence we deduce as in the paragraph following (\ref{eq.smalldist}) that $\card\{Q^i\}\leq N_1$ where $N_1=N_1(d,C_d,a_0,A_0,A_4)$. Now take
\[
M=1000nA_n'A_2a_2^{-1}a_0^{-1}.
\]
When $M$ is given as above, our present scenario is very similar to the one for $T_2$ above. More precisely, suppose that $I$ intersects $\Sigma_j$; we have that $Q_I\subset Q^i$ for some $Q^i$ as above, and $\ell(I)\lesssim r\ll\ell(Q^i)\approx Mr$. We may find $J\in\m W_{\m Q'}$ such that $J^*\cap I\neq\varnothing$ and $Q'\in\bb D_{\m F}$, so that $\ell(Q')\approx\ell(J)\approx\ell(I)\ll Mr$. With our choice of $M$, we have that $\ell(Q')<\ell(Q^i)$, so that $Q'\cap Q^i=\varnothing$. This observation gives us the estimate $\dist(Q_I,(Q_j)^c)\lesssim\ell(I)$, and we may once again use Lemma \ref{lm.dyadiccubes} \ref{item.bdrythin} (owing to our choice of $M$) to control the cardinality of the Whitney boxes $I$ intersecting $\Sigma_j^k$ by $C(2^k\ell(Q^i))^{d-\zeta}$. Thus it is easy to see that we obtain the desired result in a similar way as we did for $T_2$, by formally replacing $\ell(Q_j)$ with $\ell(Q^i)\approx Mr$. Thus ends the proof of (\ref{eq.sigmaj}).

We turn to the proof of (\ref{eq.far}). In this case, we have that $r<M_0\delta(x)$. If $r\geq\delta(x)/4$, then the desired result follows immediately from using (\ref{eq.close}); more precisely, we have that
\[
\sigma_{\star}(\Delta_\star(x,r))\leq Cr^d\leq4^{n-1-d}C\delta(x)^{d+1-n}r^{n-1},
\]
where $C$ is the constant from (\ref{eq.close}). Now suppose that $r<\delta(x)/4$. In this case, $B(x,r)\cap\Gamma=\varnothing$, and for any $X\in\Delta_\star(x,r)$, we have that $\delta(X)\geq3\delta(x)/4$, whence we deduce that
\[
\sigma_{\star}(\Delta_\star(x,r))\leq\big(\tfrac43\big)^{n-1-d}\delta(x)^{d+1-n}\n H^{n-1}(\Delta_\star(x,r))\leq A_n\big(\tfrac43\big)^{n-1-d}\delta(x)^{d+1-n}r^{n-1},
\]
and in the last inequality we used the $(n-1)-$Ahlfors-David regularity of each face $\partial J^*$ intersecting $B$, and that the number of fattened Whitney boxes $J^*$ which intersect $\Delta_\star(x,r)$ is uniformly bounded (depending only on $n$).\hfill{$\square$}

\begin{remark}\label{rm.fp} Observe that Proposition \ref{prop.dadr} implies in particular that $\n H^{n-1}(\Sigma\cap K)<\infty$ for any compact set $K\subset\bb R^n$. This is easily seen by fixing a compact set $K\subset\bb R^n$ and noticing that therefore $\delta(X)\leq M_K$ for any $X\in\Sigma\cap K$, which implies that
\[
\n H^{n-1}(\Sigma\cap K)\leq M_K^{n+1-d}\sigma_{\star}(\Sigma\cap K)<\infty.
\]
Moreover, since $\n H^d(\Gamma\cap K)\in(0,+\infty)$ for any compact set intersecting $\Gamma$, it follows that $\n H^{n-1}(\Gamma)=0$, and therefore $\partial\Omega_{\m F}$ satisfies $\n H^{n-1}(\partial\Omega_{\m F}\cap K)<+\infty$ for any compact set $K\subset\bb R^n$.
\end{remark}

We now concentrate on a lower bound for $\sigma_{\star}$.

\begin{proposition}[Lower bound for $\sigma_{\star}$]\label{prop.dadr2} Let $x\in\partial\Omega_{\m F}$ and $r>0$. Suppose that $M_0$ is given by (\ref{eq.m0}) below. If $\delta(x)\geq r/M_0$, then
\begin{equation}\label{eq.far2}
\sigma_{\star}(\Delta_{\star}(x,r))\geq v_1\delta(x)^{d+1-n}r^{n-1},
\end{equation}
where $v_1=v_1(n,d,M_0,\theta)$. If $\delta(x)< r/M_0$, then
\begin{equation}\label{eq.close2}
\sigma_{\star}(\Delta_{\star}(x,r))\geq v_2r^d,
\end{equation}
where $v_2=v_2(n,d,C_d,\theta,c_{\bb K},A_0,a_0,c,c_{\m H},M_0)$.
\end{proposition}

\noindent\emph{Proof.} We consider (\ref{eq.far2}) first, so that $\delta(x)>0$ which implies that $x\in\partial J^*$ for some $J\in\m W_{\m F}$. Hence $r/M_0\leq\delta(x)\leq42\sqrt n\ell(J)$. Observe the estimate
\begin{multline*}
\sigma_{\star}(\Delta_\star(x,r))\geq\sigma_{\star}(\partial J^*\cap\Sigma\cap B(x,r))=\int_{\partial J^*\cap\Sigma\cap B(x,r)}\delta(X)^{d+1-n}\dH(X)\\ \geq(2M_0)^{d+1-n}\delta(x)^{d+1-n}\n H^{n-1}(\partial J^*\cap\Sigma)\geq c_n(2M_0)^{d+1-n}\delta(x)^{d+1-n}\theta^{n-1}\ell(J)^{n-1}\\ \geq c_n2^dM_0^{d+2-2n}\delta(x)^{d+1-n}\theta^{n-1}r^{n-1}
\end{multline*}
where in the fourth inequality we made use of Lemma \ref{lm.faces}.

We proceed to prove (\ref{eq.close2}), using ideas of the proof for Lemma 3.61 in \cite{hm2}. First, observe that by Remark \ref{rm.fp} and Theorem \ref{thm.criterionfp}, we have that $\Omega_{\m F}$ is a set of locally finite perimeter. Hence, by the structure theorem for sets of finite perimeter, Theorem \ref{thm.structure}, it follows that $\Vert\partial\Omega_{\m F}\Vert=\n H^{n-1}\mres\partial^\star\Omega_{\m F}$. We will use these facts below.

Suppose that $\delta(x)<r/M_0$, so that there exists $\hat x\in\Gamma$ with $|x-\hat x|\leq r/M_0$. Now fix $\hat Q\in\bb D$ with $\hat x\in\hat Q$ and such that $c_{\bb K}r/M_0\leq\ell(\hat Q)\leq r/M_0$. If $M_0>\max\{9,A_0^2\}$, then we may guarantee that $\hat Q\subset B(\hat x,r/\sqrt{M_0})\subset B(x,r)=:B$. We consider two cases.

{\bf Case 1.} The ball $B(\hat x,r/\sqrt{M_0})$ meets some $Q_j\in\m F$ with $\ell(Q_j)\geq r/M_0$. Then we may procure a dyadic cube $Q\in\bb D_{Q_j}$ with $r/{2M_0}\leq\ell(Q)\leq r/M_0$ and $Q\subset B(\hat x,2r/\sqrt{M_0})$. By Lemma \ref{lm.Bbelowsawtooth}, the ball $B(x_Q,r')=B(x_Q,a_0\ell(Q)/(5A_2a_2^{-1}))$ lies in $\bb R^n\backslash\Omega_{\m F}$, while if we further assume that $M_0\geq16$, then for any $Y\in B(x_Q,r')$ we have that
\begin{equation*}
|Y-x|\leq|Y-x_Q|+|x_Q-\hat x|+|\hat x-x|<\tfrac{a_0a_2}{5A_2}\ell(Q)+\tfrac{2r}{\sqrt{M_0}}+\tfrac r{M_0} \leq 4r/\sqrt{M_0}<r,
\end{equation*}
whence it is known that $B(x_Q,r')\subset B\backslash\Omega_{\m F}$. On the other hand, we have shown in Proposition \ref{prop.h1} that $\partial\Omega_{\m F}$ has interior Corkscrew points, so that there exists $X\in\Omega_{\m F}$ verifying that $B(X,c_1r)\subset\Omega_{\m F}\cap B$, and $c_1$ is given in (\ref{eq.c1}). We can now appeal to the relative isoperimetric inequality (\ref{eq.iso}) to conclude that 
\begin{multline}\label{eq.isouse}
\Vert\partial\Omega_{\m F}\Vert(B(x,r))\geq a_n\min\big\{\n L^n(B(X,c_1r))~,~\n L^n(B(x_Q,r'))\big\}^{\frac{n-1}n}\\ \geq a_n\min\big\{c_1,\tfrac{a_0a_2}{10A_2M_0}\big\}^{n-1}r^{n-1},
\end{multline}
where $a_n$ is a uniform constant depending only on $n$.
Consequently,
\begin{multline}\label{eq.lowerbound1}
\sigma_{\star}(\Delta_\star(x,r))\geq\sigma_{\star}(\Delta_\star(x,r)\cap\Sigma)\geq M_0^{n-1-d}r^{d+1-n}\n H^{n-1}(\Delta_\star\cap\Sigma)\\=M_0^{n-1-d}r^{d+1-n}\n H^{n-1}(B(x,r)\cap\partial\Omega_{\m F})\geq M_0^{n-1-d}r^{d+1-n}\n H^{n-1}(B(x,r)\cap\partial^\star\Omega_{\m F})\\=M_0^{n-1-d}r^{d+1-n}\Vert\partial\Omega_{\m F}\Vert(B(x,r))\geq a_n(a_0a_2A_2^{-1})^{n-1}M_0^{-d}r^d,
\end{multline}
where we have used the structure theorem for sets of finite perimeter, Theorem \ref{thm.structure}.

{\bf Case 2.} There is no $Q_j$ as in case 1. It follows that if $Q_j\in\m F$ meets $B(\hat x,r/\sqrt{M_0})$, then $\ell(Q_j)\leq r/M_0$. Let $\hat{\m F}$ denote the collection of those $Q_j\in\m F$ which intersect $\hat\Delta=\Delta(\hat x,r/\sqrt{M_0})$. Then we may split $\frac12\hat\Delta$ as
\[
\tfrac12\hat\Delta=\big(\tfrac12\hat\Delta\backslash\big(\medcup_{\hat{\m F}}Q_j\big)\big)\cup\big(\tfrac12\hat\Delta\cap(\medcup_{\hat{\m F}}Q_j)\big).
\]
We now consider two sub-cases: either the estimate
\begin{equation}\label{eq.eat}
\sigma\big(\tfrac12\hat\Delta\backslash\big(\medcup_{\hat{\m F}}Q_j\big)\big)\geq\frac12\sigma\big(\tfrac12\hat\Delta\big)
\end{equation}
holds, or it does not. If it does, then we deduce that
\begin{multline*}
\sigma_{\star}(B(x,r)\cap\partial\Omega_{\m F})\geq\sigma\big(B(x,r)\cap\Gamma\backslash\medcup_{\m F}Q_j\big)\geq\sigma\big(\hat\Delta\backslash\medcup_{\hat{\m F}}Q_j\big)\geq \frac12\sigma\big(\tfrac12\hat\Delta\big)\\ \geq C_d^{-1}2^{-d-1}M_0^{-d/2}r^d=:a_5r^d,
\end{multline*}
which yields the desired result. We are left to consider the case that (\ref{eq.eat}) does not hold. Then, instead, we have that
\begin{equation}\label{eq.noeat}
\medsum_{\m F'}\sigma(Q_j)\geq\sigma\big(\tfrac12\hat\Delta\cap\big(\medcup_{\hat{\m F}}Q_j\big)\big)\geq\frac12\sigma\big(\frac12\hat\Delta\big)\geq a_5r^d,
\end{equation}
where $\m F'$ consists of those $Q_j\in\hat{\m F}$ which intersect $\frac12\Delta$. For each $Q_j\in\m F'$, fix any one of the $(n-1)-$dimensional cubes $P_j\subset\partial\Omega_{\m F}$ constructed in Proposition \ref{prop.cubeinface}, and denote its center by $x_j^\star$. We now claim that for each $Q_j\in\m F'$, the ball $B_j^*=B(x_j^\star,16A_0c_{\bb K}^{-1}\ell(Q_j))$ contains both an interior and an exterior Corkscrew ball for $\Omega_{\m F}$, with respect to the surface ball $B_{Q_j}^*\cap\partial\Omega_{\m F}$ (with Corkscrew constants that could depend on $K_0$). 

Indeed, by virtue of Lemma \ref{lm.Bbelowsawtooth}, the ball $B_j=B(x_{Q_j}, a_0a_2\ell(Q_j)/(5A_2))$ is contained in $\bb R^n\backslash\Omega_{\m F}$, and we note that for any $Y\in B_j$,
\begin{multline*}
|Y-x_j^\star|<\frac{a_0a_2}{5A_2}\ell(Q_j)+\diam Q_j+\dist(P_j,Q_j)+\diam P_j\\ \leq2A_0\ell(Q_j)+11A_0c_{\bb K}^{-1}\ell(Q_j)+\tfrac{13}8\theta a_0c_{\bb K}^{-1}\ell(Q_j)\leq15A_0c_{\bb K}^{-1}\ell(Q_j),
\end{multline*}
so that $B_j\subset B_j^*\backslash\Omega_{\m F}$. By a similar reasoning, we also have that $Q_j\subset B_j^*$. For the interior Corkscrew ball, let $\hat Q\in\bb D_{\m F}$ be the proper parent of $Q_j$ and fix $I\in\m W_{\hat Q}^{\operatorname{cs}}$. Note that $B(X_I,\ell(I)/2)\subset I\subset\opint I^*\subset\Omega_{\m F}$, and by our choice of radius for $B_j^*$, we also have that $B(X_I,\ell(I)/2)\subset B_j^*$ similarly as in the estimate above. Then the ball $B(X_I,a_0c\ell(Q_j)/(41\sqrt n))$ is contained in $B_j^*\cap\Omega_{\m F}$ by (\ref{eq.whitneyboxesc}). Thus by using the relative isoperimetric inequality in a manner analogous to (\ref{eq.isouse}), we deduce that
\[
\Vert\partial\Omega_{\m F}\Vert(B_j^*)\geq a_n(a_0a_2/(5A_2))^{n-1}\ell(Q_j)^{n-1}.
\]
We also have that for any $Y\in B_j^*$,
\[
\delta(Y)\leq16A_0c_{\bb K}^{-1}\ell(Q_j)+\diam P_j+\dist(P_j,\Gamma)\leq400A_0c_{\bb K}^{-1}\ell(Q_j),
\]
and therefore, analogous to (\ref{eq.lowerbound1}), we obtain that
\begin{equation}\label{eq.eachqj}
\sigma_{\star}(B_j^*\cap\partial\Omega_{\m F})\geq a_n(a_0a_2A_2^{-1})^{n-1}A_0^{d+1-n}c_{\bb K}^{n-1-d}\ell(Q_j)^d.
\end{equation}

We show that for $M_0\geq1250A_0^2c_{\bb K}^{-2}$, we have that $B_j^*\subset B$. Fix any $Y\in B_j^*$ and observe that
\begin{multline*}
|Y-x|\leq\diam(B_j^*)+\diam Q_j+\tfrac12\diam B\big(\hat x,r/\sqrt{M_0}\big)+|\hat x-x|\\ \leq32A_0c_{\bb K}^{-1}\ell(Q_j)+A_0\ell(Q_j)+\frac r{\sqrt{M_0}}+\frac r{M_0}\leq\frac{35A_0c_{\bb K}^{-1}}{\sqrt{M_0}}r<r,
\end{multline*}
as claimed. Henceforth we will take 
\begin{equation}\label{eq.m0}
M_0:=125A_0^2c_{\bb K}^{-2}.
\end{equation}

Let us now show that we can muster a sub-collection $\m F''\subset\m F'$ of cubes $Q_j\in\m F'$ such that the balls in $\{B_j^*\}_{Q_j\in\m F''}$ are pairwise disjoint and
\begin{equation}\label{eq.aftercover}
\sum_{Q_j\in\m F''}\ell(Q_j)\gtrsim r^d.
\end{equation}
Since for any $Q_j\in\m F'\subset\hat{\m F}$ we have that $\ell(Q_j)<r/M_0$, it follows that there exists $k_0\in\bb Z$ such that $\ell(Q_j)\leq2^{-k_0}$ for all $Q_j\in\m F'$, and $\ell(Q)=2^{-k_0}$ for some $Q\in\m F'$. For any $k\geq k_0$, let $\m F'_k=\{Q_j\in\m F':\ell(Q_j)=2^k\}$. Fix a sub-collection $\m F_{k_0}''$ of $\m F_{k_0}'$ which is \emph{$B_j^*-$maximal} in the sense that the balls $\{B_j^*\}_{Q_j\in\m F_{k_0}''}$ are pairwise disjoint, but where adjoining any other cube in $\m F'\backslash\m F_{k_0}''$ makes some of these balls overlap. Next, define inductively for each $k>k_0$ the sub-collection $\m F_k''$ which is the union of all $\m F_{\tilde k}''$, $k_0\leq\tilde k<k$, and adjoined with a sub-collection of $\m F_k'$ such that $\m F''_k$ is $B_j^*-$maximal. We then set $\m F''=\cup_{k\geq k_0}\m F_k''$ and observe that it satisfies that the balls in $\{B_j^*\}_{Q_j\in\m F''}$ are pairwise disjoint, and that each $Q_m\in\m F'\backslash\m F''$ is such that $B_m^*$ intersects $B_j^*\in\m F''$ for some $Q_j\in\m F''$ with $\ell(Q_m)\leq\ell(Q_j)$ (otherwise, $Q_m$ would have had to belong to $\m F_k''$ for some $k$). 
 
Recall that for any $Q_j\in\m F'$, $Q_j\subset B_j^*$. If $Q_m\in\m F'$ intersects $Q_j$ with $\ell(Q_m)\leq\ell(Q_j)$, then
\[
\dist(Q_m,Q_j)\leq\diam B_m^*+\diam B_j^*\leq32A_0c_{\bb K}^{-1}\ell(Q_j),
\]
and thus
\[
Q_m\subset B(x_{Q_j},34A_0c_{\bb K}^{-1}\ell(Q_j))\cap\Gamma,
\]
which, together with our observations in the last paragraph, implies that
\[
\medcup_{\m F'}Q_j\subset\medcup_{Q_j\in\m F''}B\big(x_{Q_j},34A_0c_{\bb K}^{-1}\ell(Q_j)\big)\cap\Gamma.
\]
Hence, we estimate
\begin{equation*}
\sum_{\m F'}\sigma(Q_j)\leq\sigma\big(\medcup_{\m F''}B\big(x_{Q_j},34A_0c_{\bb K}^{-1}\ell(Q_j)\big)\cap\Gamma\big)\leq C_d(34A_0c_{\bb K}^{-1})^d\displaystyle\sum_{\m F''}\ell(Q_j)^d
\end{equation*}
By combining this last estimate with (\ref{eq.noeat}), we obtain (\ref{eq.aftercover}) with implicit constant $a_5C_d^{-1}(34^{-1}A_0^{-1}c_{\bb K})^{d}$.  Finally, we combine (\ref{eq.eachqj}) and (\ref{eq.aftercover}) to conclude that 
\begin{multline*}
\sigma_{\star}(B\cap\partial\Omega_{\m F})\geq\sum_{Q_j\in\m F''}\sigma_{\star}(B_j^*\cap\partial\Omega_{\m F})\geq a_n(a_2^2A_2^{-2})^{n-1}A_0^dc_{\bb K}^{n-1-d}\sum_{Q_j\in\m F''}\ell(Q_j)^d\\ \geq35^{-d}C_d^{-2}M_0^{-d/2}a_n(a_2^2c_{\bb K}A_2^{-2})^{n-1}r^d,
\end{multline*}
which does complete our argument for (\ref{eq.close2}).\hfill{$\square$}

We are ready to show

\begin{proposition}[The measure $\sigma_{\star}$ is doubling on the sawtooth domain]\label{prop.doubling} The measure $\sigma_{\star}$ satisfies \ref{ax.doublingmeas}.
\end{proposition}

\noindent\emph{Proof.} Fix $x\in\partial\Omega_{\m F}$ and $r>0$. We split the proof of the proposition into three cases.

{\bf Case 1.} $2r<M_0\delta(x)$. Then we also have that $r<M_0\delta(x)$, and therefore,
\[
\sigma_{\star}\big(B(x,2r)\cap\Gamma\big)\leq 2^{n-1}V_2r^{n-1}\delta(x)^{d+1-n}\leq2^{n-1}V_2v_1^{-1}\sigma_{\star}\big(B(x,r)\cap\Gamma\big).
\]

{\bf Case 2.} $r<M_0\delta(x)\leq2r$. Here, we obtain that
\begin{multline*}
\sigma_{\star}\big(B(x,2r)\cap\Gamma\big)\leq2^dV_1r^d=2^dV_1r^{d+1-n}r^{n-1}\leq2^{-2d+1-n}M_0^{d+1-n}V_1\delta(x)^{d+1-n}r^{n-1}\\ \leq2^{-2d+1-n}M_0^{d+1-n}V_1v_1^{-1}\sigma_{\star}\big(B(x,r)\cap\Gamma\big).
\end{multline*}

{\bf Case 3.} $r\geq M_0\delta(x)$. We easily estimate that
\[
\sigma_{\star}\big(B(x,2r)\cap\Gamma\big)\leq2^dV_1r^d\leq2^dV_1v_2^{-1}\sigma_{\star}\big(B(x,r)\cap\Gamma\big),
\]
as desired.\hfill{$\square$}

Next, we turn to verifying the growth condition \ref{ax.growth}. In preparation, we record the following useful estimate from \cite{dfm1}.

\begin{lemma}[Behavior of $m$, \cite{dfm1} Lemma 2.3, Remark 2.4]\label{lm.m} For any $\alpha>0$, there exists a constant $M(\alpha)>0$, depending only on $n$, $d$, $C_d$, and $\alpha$, such that the following statements hold for any $X\in\bb R^n$ and any $r>0$.
\begin{enumerate}[(i)]
	\item If $\delta(X)\geq\alpha r$, then
	\begin{equation}\label{eq.m1}\notag
	M(\alpha)^{-1}r^n\delta(X)^{d+1-n}\leq m\big(B(X,r)\big)\leq M(\alpha)r^n\delta(X)^{d+1-n}.
	\end{equation}
	\item If $\delta(X)\leq\alpha r$, then
	\begin{equation}\label{eq.m2}\notag
	M(\alpha)^{-1}r^{d+1}\leq m\big(B(X,r)\big)\leq M(\alpha)r^{d+1}.
	\end{equation}
\end{enumerate}
\end{lemma}

\begin{proposition}[Growth condition]\label{prop.growth} The measures $\sigma_{\star}$ and $m$ satisfy \ref{ax.growth}. More precisely, there exist constants $V_5\geq1$ and $\ep\in(0,1)$ so that for each $x\in\partial\Omega_{\m F}$ and all $r,s$ with $0<s<r$, the estimate
\begin{equation}\label{eq.othergrowth}\notag
\frac{m\big(B(x,r)\cap\Omega_{\m F}\big)}{m\big(B(x,s)\cap\Omega_{\m F}\big)}\leq V_5\Big(\frac rs\Big)\frac{\sigma_{\star}\big(B(x,r)\cap\Gamma\big)}{\sigma_{\star}\big(B(x,s)\cap\Gamma\big)}
\end{equation}
holds.
\end{proposition}

\noindent\emph{Proof.} Recall that $M_0$ is the uniform constant given in (\ref{eq.m0}). Note that $m\big(B(x,r)\cap\Omega_{\m F}\big)\leq m\big(B(x,r)\big)$ for any $x\in\partial\Omega_{\m F}$ and $r>0$, while in the other direction, Proposition \ref{prop.h1} implies the existence of a Corkscrew point $X=X_{x,r}\in\Omega_{\m F}$ such that $B(X,c_1r)\subset\Omega_{\m F}\cap B(x,r)$. Observe that such a Corkscrew point satisfies
\[
c_1r\leq\delta(X_{x,r})\leq r+\delta(x).
\]
Now fix $x\in\partial\Omega_{\m F}$ and $r,s>0$ with $0<s<r$. We consider three cases.

{\bf Case 1.} $\delta(x)\geq r/M_0\geq s/M_0$. Here, note that $\delta(X_{x,s})\leq2M_0\delta(x)$. We have that
\begin{multline*}
\frac{m\big(B(x,r)\cap\Omega_{\m F}\big)}{m\big(B(x,s)\cap\Omega_{\m F}\big)}\leq\frac{m\big(B(x,r)\big)}{m\big(B(X_{x,s},c_1s)\big)}\leq\frac{M(\tfrac1{M_0})r^n\delta(x)^{d+1-n}}{M(1)^{-1}(c_1s)^n\delta(X_{x,s})^{d+1-n}}\\[2mm] \leq\Big[(2M_0)^{n-1-d}\frac{M(1)M(\tfrac1{M_0})}{c_1^n}\Big]\Big(\frac rs\Big)\frac{r^{n-1}\delta(x)^{d+1-n}}{s^{n-1}\delta(x)^{d+1-n}}\\ \leq\Big[(2M_0)^{n-1-d}\frac{M(1)M(\tfrac1{M_0})V_2}{c_1^nv_1}\Big]\Big(\frac rs\Big)\frac{\sigma_{\star}\big(B(x,r)\cap\Gamma\big)}{\sigma_{\star}\big(B(x,s)\cap\Gamma\big)},
\end{multline*}
where we have used (\ref{eq.far}) and (\ref{eq.far2}), and thus established the desired estimate in this case.

{\bf Case 2.} $\delta(x)\leq s/M_0\leq r/M_0$. In this case, we have that $\delta(X_{x,s})\leq2s$. Reckon that
\begin{multline*}
\frac{m\big(B(x,r)\cap\Omega_{\m F}\big)}{m\big(B(x,s)\cap\Omega_{\m F}\big)}\leq\frac{m\big(B(x,r)\big)}{m\big(B(X_{x,s},c_1s)\big)}\leq\frac{M(\tfrac1{M_0})r^{d+1}}{M(\tfrac2{c_1})^{-1}(c_1s)^{d+1}}\\[2mm] \leq\Big[\frac{M(\tfrac2{c_1})M(\tfrac1{M_0})}{c_1^{d+1}}\Big]\Big(\frac rs\Big)\frac{r^d}{s^d} \leq\Big[\frac{M(\tfrac2{c_1})M(\tfrac1{M_0})V_1}{c_1^{d+1}v_2}\Big]\Big(\frac rs\Big)\frac{\sigma_{\star}\big(B(x,r)\cap\Gamma\big)}{\sigma_{\star}\big(B(x,s)\cap\Gamma\big)},
\end{multline*}
where this time we made use of (\ref{eq.close}) and (\ref{eq.close2}). 

{\bf Case 3.} $s/M_0<\delta(x)<r/M_0$. Now we see that $\delta(X_{x,s})\leq2M_0\delta(x)$, and estimate
\begin{multline*}
\frac{m\big(B(x,r)\cap\Omega_{\m F}\big)}{m\big(B(x,s)\cap\Omega_{\m F}\big)}\leq\frac{m\big(B(x,r)\big)}{m\big(B(X_{x,s},c_1s)\big)}\leq\frac{M(\tfrac1{M_0})r^{d+1}}{M(1)^{-1}(c_1s)^n\delta(X_{x,s})^{d+1-n}}\\[2mm] \leq\Big[(2M_0)^{n-1-d}\frac{M(1)M(\tfrac1{M_0})}{c_1^n}\Big]\Big(\frac rs\Big)\frac{r^d}{s^{n-1}\delta(x)^{d+1-n}} \\ \leq\Big[(2M_0)^{n-1-d}\frac{M(1)M(\tfrac1{M_0})}{c_1^n}\Big]\Big(\frac rs\Big)\frac{\sigma_{\star}\big(B(x,r)\cap\Gamma\big)}{\sigma_{\star}\big(B(x,s)\cap\Gamma\big)},
\end{multline*}
using (\ref{eq.close2}) and (\ref{eq.far}). The desired result is established in any case.\hfill{$\square$}

\begin{remark}\label{rm.done} Incidentally, Lemma \ref{lm.m} together with Proposition \ref{prop.h1}   give that \ref{ax.misgood} holds (the doubling property on $\partial\Omega_{\m F}$ follows from the existence of interior Corkscrew balls; for a similar analysis see the first paragraph of the proof of Proposition \ref{prop.growth}). Finally, it is a trivial application of the fundamental results in \cite{dfm1} that $m$ satisfies the axiom \ref{ax.trace}, since $m|_{\Omega_{\m F}}$ is merely the restriction of the function $m$ on $\Omega$ which satisfies this property.
\end{remark}

\section{Carleson measures, discrete Carleson measures, and extrapolation}\label{sec.extrapolation}

\begin{definition}[Carleson measures]\label{def.carleson} We say that a non-negative Borel measure $\lambda$ on $\Omega$ is a \emph{Carleson measure} if
	\begin{equation}\label{eq.carlesonnorm}\notag
	\vertiii{\lambda}_{\m C}:=\sup\limits_{\Delta\subset\Gamma}\frac{1}{\sigma(\Delta)}\lambda(T_{\Delta})<\infty,
	\end{equation}
	where, if $\Delta=B(x,r)\cap\Gamma$, then $T_{\Delta}=B(x,r)\cap\Omega$ is the \emph{tent over $\Delta$}. The supremum runs over all surface balls $\Delta\subset\Gamma$. We call $\vertiii{\lambda}_{\m C}$ the \emph{Carleson norm} of $\lambda$, and we write $\m C$ for the set of all Carleson measures in $\Omega$. 
\end{definition}

A main tool in our proof is the extrapolation of Carleson measures, which we use in the dyadic setting. We borrow the definitions and results from \cite{hm2}, where this result has been considered in a co-dimension 1 setting; see also \cite{chm}. In the setting of higher co-dimension, this framework has appeared in \cite{dm}.

\begin{definition}[Discrete Carleson measures] Let $\{\alpha_Q\}_{Q\in\bb D}$ be a sequence of non-negative numbers indexed by $Q\in\bb D$, and for any sub-collection $\bb D'\subset\bb D$, we define
	\[
	\f m(\bb D'):=\sum_{Q\in\bb D'}\alpha_Q.
	\]
	We say that $\f m$ is a \emph{discrete Carleson measure} on $\bb D$ with respect to $\sigma$ (written $\f m\in\m C$) if
	\[
	\Vert\f m\Vert_{\m C}:=\sup_{Q\in\bb D}\frac{\f m(\bb D_Q)}{\sigma(Q)}<\infty.
	\]
	Similarly, we have a local version: For a fixed $Q_0\in\bb D$, we say that $\f m$ is a discrete Carleson measure on $\bb D_{Q_0}$ with respect to $\sigma$ (written $\f m\in\m C(Q_0)$) if
	\[
	\Vert\f m\Vert_{\m C(Q_0)}:=\sup_{Q\in\bb D_{Q_0}}\frac{\f m(\bb D_Q)}{\sigma(Q)}<\infty.
	\]
	Moreover, set $\bb D_Q^{\short}:=\bb D_Q\backslash\{Q\}$, and given a disjoint family $\m F\subset\bb D$, we define the restriction of $\f m$ to the sawtooth $\bb D_{\m F}$ by
	\[
	\f m_{\m F}(\bb D'):=\f m(\bb D'\cap\bb D_{\m F})=\sum_{Q\in\bb D'\backslash(\cup_{\m F}\bb D_{Q_j})}\alpha_Q,\qquad\text{for }\bb D'\subset\bb D_{Q_0}.
	\]
\end{definition}

The following result concerns the extrapolation of Carleson measures, and is a main tool in our proof.

\begin{theorem}[Extrapolation of Carleson measures; \cite{dm} \cite{hm2}]\label{thm.extrapolation} Let $\Gamma$ be a closed \dADR set, and recall that $\sigma=\n H^d|_{\Gamma}$. Fix $Q_0\in\bb D$ and a dyadically doubling Borel measure $\mu$ on $Q_0$. Assume that there is some sequence of non-negative numbers $\{\alpha_Q\}_{Q\in\bb D(Q_0)}$ such that the corresponding $\f m$ satisfies
	\[
	\Vert\f m\Vert_{\m C(Q_0)}\leq M_0,
	\]
	for some $M_0<\infty$. Suppose that there exists $\xi>0$ such that for every $Q\in\bb D_{Q_0}$, and every disjoint family $\m F\subset\bb D_Q$ verifying
	\[
	\Vert\f m_{\m F}\Vert_{\m C(Q)}=\sup_{Q'\in\bb D_Q}\frac{\f m(\bb D_{\m F,Q'})}{\sigma(Q')}\leq\xi,
	\]
	we have that $\m P_{\m F}\mu$ satisfies the following property:
	\begin{equation}\label{eq.propextra}
	\forall\ep\in(0,1),~\exists C_{\ep}>1~\text{such that }\Big(F\subset Q,~\frac{\sigma(F)}{\sigma(Q)}\geq\ep\implies\frac{\m P_{\m F}\mu(F)}{\m P_{\m F}\mu(Q)}\geq\frac1{C_{\ep}}\Big).
	\end{equation}
	Then, there exist $\eta_0\in(0,1)$ and $C_0<\infty$ such that, for every $Q\in\bb D_{Q_0}$,
	\begin{equation}\label{eq.concluextra}\notag
	F\subset Q,~\frac{\sigma(F)}{\sigma(Q)}\geq1-\eta_0\implies\frac{\mu(F)}{\mu(Q)}\geq\frac1{C_0}.
	\end{equation}
	In other words, $\mu\in A_{\infty}^{\dyadic}(Q_0)$.
\end{theorem}

Let us elucidate how Theorem \ref{thm.extrapolation} will be used to prove Theorem \ref{thm.main}. In the hypothesis of the latter, we have that the measure $d\lambda(X):=\delta(X)^{d-n}\f a^2\,dX$ is a (continuous) Carleson measure, where $\f a$ is defined in (\ref{eq.disagreement}). In the following lemma, we define the natural discrete version of this measure, and show that it is indeed a discrete Carleson measure as one might expect.

\begin{lemma}\label{lm.todyadic} Suppose that $\m A_0$, $\m A$ are two uniformly elliptic matrices, such that their disagreement $\f a$ defined in (\ref{eq.disagreement}) satisfies that $d\lambda(X)=\delta(X)^{d-n}\f a^2\,dX$ is a Carleson measure. Then, for every $Q_0\in\bb D$, the collection $\f m=\{\alpha_Q\}_{Q\in\bb D_{Q_0}}$ with
	\begin{equation}\label{eq.discretecarleson}
	\alpha_Q:=\sum_{I\in\m W_Q}\frac{\sup_{Y\in I^*}|\f E(Y)|^2}{\ell(I)^{n-d}}|I|,\qquad Q\in\bb D,
	\end{equation}
	is a discrete Carleson measure, and in fact, 
	\[
	\Vert\f m\Vert_{\m C(Q_0)}\leq(41\sqrt n)^{n-d}(7\sqrt nA_2a_0^{-1})^dC_d^2\vertiii{\lambda}_{\m C}.
	\]
\end{lemma}

\noindent\emph{Proof.} Let $Q\in\bb D_{Q_0}$, write $t_Q=7\sqrt nA_2\ell(Q)$, and consider the estimates
\begin{multline*}
\frac{\f m(\bb D_Q)}{\sigma(Q)}=\frac1{\sigma(Q)}\sum_{Q'\in\bb D_Q}\sum_{I\in\m W_{Q'}}\frac{\sup_{Y\in I^*}|\f E(Y)|^2}{\ell(I)^{n-d}}|I|\\ \leq(41\sqrt n)^{n-d}\frac1{\sigma(Q)}\sum_{I\in\m R_Q}\dint_I\frac{\sup_{Y\in I^*}|\f E(Y)|^2}{\delta(X)^{n-d}}\,dX\\ \leq (41\sqrt n)^{n-d}\frac1{\sigma(Q)}\sum_{I\in\m R_Q}\dint_I\frac{\f a^2(X)}{\delta(X)^{n-d}}\,dX\leq (41\sqrt n)^{n-d}\frac1{\sigma(Q)}\dint_{R_Q}\frac{\f a^2(X)}{\delta(X)^{n-d}}\,dX\\ \leq (41\sqrt n)^{n-d}(7\sqrt nA_2a_0^{-1})^dC_d^2\frac1{\sigma(\Delta(x_Q,t_Q))}\dint_{B(x_Q,t_Q)\cap\Omega}\frac{\f a^2(X)}{\delta(X)^{n-d}}\,dX\\[2mm] \leq(41\sqrt n)^{n-d}(7\sqrt nA_2a_0^{-1})^dC_d^2\vertiii{\lambda}_{\m C},
\end{multline*}
where we have used in the third line that $B(X,\delta(X)/2)\supset I^*$ for any $X\in I$, and later we used  that $R_Q\subset B(x_Q,7\sqrt nA_2\ell(Q))$ (by the same argument as (\ref{eq.sawtoothbdd})), and (\ref{eq.ahlforsreg}).\hfill{$\square$}

\section{Review of the elliptic theory for sets with boundaries of high co-dimension}\label{sec.theory}

Let us review the necessary background and theory of the David-Feneuil-Mayboroda operators \cite{dfm1}. Before starting, we remark that many of the results in this section have direct analogues for our sawtooth domains by virtue of Theorem \ref{thm.sawtooth} and the elliptic theory for sets of mixed dimension carried out in \cite{dfm20}.
 
Formally, we write
\[
L=-\dv A\nabla,
\]
with $A:\Omega\ra\bb M_n(\bb R)$, where $\bb M_n(\bb R)$ is the set of $n\times n$ real-valued matrices, and we require that $A$ satisfies the following weighted boundedness and ellipticity conditions:
\begin{gather}
\delta(X)^{n-d-1}A(X)\xi\cdot\zeta\leq C_A|\xi||\zeta|,\quad\text{for each }X\in\Omega\text{ and every }\xi,\zeta\in\bb R^n,\notag\\ \delta(X)^{n-d-1}A(X)\xi\cdot\xi\geq C_A^{-1}|\xi|^2,\quad\text{for each }X\in\Omega\text{ and every }\xi\in\bb R^n.\label{eq.elliptic}
\end{gather}
Recall that we denote $w(X)=\delta(X)^{-n+d+1}$ and $m(E)=\dint_Ew(X)\,dX$. By $\m A$ we denote the matrix $w^{-1}A$, so that 
\[
\dint_{\Omega}A\nabla u\nabla v=\dint_{\Omega}\m A\nabla u\nabla v\,dm.
\]
The matrix $\m A$ satisfies unweighted ellipticity and boundedness conditions
\begin{gather}
\m A(X)\xi\cdot\zeta\leq C|\xi||\zeta|,\quad\text{for each }X\in\Omega\text{ and every }\xi,\zeta\in\bb R^n,\notag\\ \m A(X)\xi\cdot\xi\geq C^{-1}|\xi|^2,\quad\text{for each }X\in\Omega\text{ and every }\xi\in\bb R^n.\label{eq.elliptic2}
\end{gather}

In order to rigorously define the operator $L$, we need a suitable domain and corresponding range. As in \cite{dfm1}, we consider the following weighted Sobolev space,
\begin{equation*}\label{eq.sobspace}
W=\dot W^{1,2}_w(\Omega):=\big\{u\in L^1_{loc}(\Omega)\,:\,\nabla u\in L^2(\Omega,dm)\big\},
\end{equation*}
and set
\begin{equation*}\label{eq.wnorm}
\Vert u\Vert_W=\Big(\dint_{\Omega}|\nabla u|^2\,dm\Big)^{\frac12},\qquad u\in W.
\end{equation*}
Actually, it is proven in \cite{dfm1} that
\begin{equation*}\label{eq.bigspace}
W=\big\{u\in L^1_{loc}(\bb R^n):\nabla u\in L^2(\bb R^n,dm)\big\}.
\end{equation*}

If $E\subset\bb R^n$ is a Borel set, we let $C_c^{\infty}(E)$ denote the space of compactly supported, smooth functions on $E$. We call $W_0$ the completion of $C_c^{\infty}(\Omega)$ in the norm $\Vert\cdot\Vert_W$. Finally, denote by $\m M(\Gamma)$ the set of $\sigma-$measurable functions on $\Gamma$, and then set
\begin{equation*}\label{eq.fracsob}
H=\dot H^{1/2}(\Gamma):=\Big\{g\in\m M(\Gamma)\,:\,\int_{\Gamma}\int_{\Gamma}\frac{|g(x)-g(y)|^2}{|x-y|^{d+1}}\,d\sigma(x)\,d\sigma(y)<\infty\Big\}.
\end{equation*}
The significance of $H$ is that it plays a role for $W$ analogous in many ways to the role that the fractional Sobolev space $H^{\frac12}$ plays for the classical Sobolev space $W^{1,2}$. 

In addition to $W$ which is a space of functions defined globally, we introduce a local version. Let $E\subset\bb R^n$ be an open set. The set of functions $W_r(E)$ is defined as
\begin{equation*}\label{eq.localspace}
W_r(E)=\big\{f\in L^1_{loc}(E)\,:\,\varphi f\in W~\text{for all}~\varphi\in C_c^{\infty}(E)\big\}
\end{equation*}
where $\varphi f$ is seen as a function on $\bb R^n$.

The following two results establish that we can make sense of traces on $\Gamma$ of functions in this weighted Sobolev space.

\begin{theorem}[Trace operator, Theorem 3.4 of \cite{dfm1}]\label{thm.traceop} There exists a bounded linear operator $T:W\ra H$ (a trace operator) with the following properties. The trace of $u\in W$ is such that, for $\sigma-$almost every $x\in\Gamma$,
	\begin{equation*}\label{eq.tracelim}
	Tu(x)=\lim\limits_{r\ra0}\,\fiint_{B(x,r)}u(X)\,dX,
	\end{equation*}
	and, analogously to the Lebesgue density property,
	\begin{equation*}\label{eq.lebesguedensity}
	\lim\limits_{r\ra0}\,\fiint_{B(x,r)}|u(X)-Tu(x)|\,dX=0.
	\end{equation*}
\end{theorem}

\begin{lemma}[Local traces, Lemma 8.1 of \cite{dfm1}] Let $E\subset\bb R^n$ be an open set. For every function $u\in W_r(E)$, we can define the trace of $u$ on $\Gamma\cap E$ by
	\begin{equation*}\label{eq.localtrace}
	Tu(x)=\lim\limits_{r\ra0}\,\fiint_{B(x,r)}u(X)\,dX\qquad\text{for }\sigma-\text{almost every }x\in\Gamma\cap E,
	\end{equation*}
and $Tu\in L^1_{loc}(\Gamma\cap E,\sigma)$. Moreover, for every choice of $f\in W_r(E)$ and $\varphi\in C_c^{\infty}(E)$,
\begin{equation*}\label{eq.tracemult}
\big(T(\varphi u)\big)(x)=\varphi(x)Tu(x)\qquad\text{for }\sigma-\text{almost every }x\in\Gamma\cap E.
\end{equation*}
\end{lemma}

Next, we give a meaning to a local solution of the problem $Lu=0$.

\begin{definition}[Local weak solutions]\label{def.weaksol} Let $E\subseteq\Omega$ be an open set. We say that $u\in W_r(E)$ is a \emph{solution} of $Lu=0$ in $E$ if for any $\varphi\in C_c^{\infty}(E)$,
\begin{equation}\label{eq.weakid}
\dint_{\Omega}A\nabla u\nabla\varphi\,dX=\dint_{\Omega}\m A\nabla u\nabla\varphi\,dm=0.
\end{equation}
\end{definition}

The following lemma allows us in certain cases to extend the space of test functions we have available for (\ref{eq.weakid}).

\begin{lemma}[Extensions of test space, Lemma 8.3 of \cite{dfm1}]\label{lm.moretestfns}  Let $E\subset\Omega$ be an open set and let $u\in W_r(E)$ be a solution of $Lu=0$ in $E$. Also denote by $E^{\Gamma}$ the interior of $E\cup\Gamma$. The identity (\ref{eq.weakid}) holds:
	\begin{itemize}
		\item when $\varphi\in W_0$ is compactly supported in $E$;
		\item when $\varphi\in W_0$ is compactly supported in $E^{\Gamma}$ and $u\in W_r(E^{\Gamma}).$
		\item when $E=\Omega$, $\varphi\in W_0$, and $u\in W$.
	\end{itemize}
\end{lemma}

We now consider the mapping properties of the operator $L$. The proof of the next result is a routine application of the Lax-Milgram theorem.

\begin{lemma}[The bilinear form and the operator $L$]\label{lm.laxmilgram} Let $A$ be a matrix satisfying (\ref{eq.elliptic}). Then there exists a unique, linear, invertible operator $L: W\ra W^*$ such that if $u\in W$, then the identity
\begin{equation}\label{eq.lmweak}\notag
B[u,\phi]:=\dint_{\Omega}A\nabla u\nabla\phi\,dX=\dint_{\Omega}\phi Lu\,dX,
\end{equation}
holds for each $\phi\in W$.
\end{lemma} 

We have an analogous version of the Harnack inequality.

\begin{lemma}[Harnack inequality; Lemma 8.9 of \cite{dfm1}]\label{lm.harnack} Let $B$ be a ball such that $3B\subseteq\Omega$, and let $u\in W_r(3B)$ be a non-negative solution in $3B$. Then
	\begin{equation}\label{eq.harnack}\notag
	\sup_Bu\leq C\inf_Bu,
	\end{equation}
	where $C$ depends only on $n$, $d$, $C_d$, and $C_A$.
\end{lemma}

Now, we exhibit results concerning the \emph{Green function}. 

\begin{lemma}[Green's function, Lemma 10.1 of \cite{dfm1}]\label{lm.greenfn}  There exists a non-negative function $g:\Omega\times\Omega\ra\bb R\cup\{+\infty\}$ with the following properties.
\begin{enumerate}[(i)]
	\item For any $Y\in\Omega$ and any $\alpha\in C_c^{\infty}(\bb R^n)$ such that $\alpha\equiv1$ in a neighborhood of $y$,
	\[
	(1-\alpha)g(\cdot,Y)\in W_0.
	\]
	In particular, $g(\cdot,Y)\in W_r(\bb R^n\backslash\{Y\})$ and $T[g(\cdot,Y)]=0$.
	\item For every choice of $Y\in\Omega$, $R>0$, and $q\in[1,\frac n{n-1})$,
	\[
	g(\cdot,Y)\in W^{1,q}(B(Y,R)):=\big\{u\in L^q(B(Y,R)),\nabla u\in L^q(B(Y,R))\big\}.
	\]
	\item For $Y\in\Omega$ and $\varphi\in C_c^{\infty}(\Omega)$,
	\begin{equation*}\label{eq.greenid}
	\dint_{\Omega}A\nabla_Xg(X,Y)\nabla\varphi(X)\,dX=\varphi(Y).
	\end{equation*}
	In particular, $g(\cdot,Y)$ is a solution of $Lu=0$ in $\Omega\backslash\{Y\}$.
	\item For $r>0$, $Y\in\Omega$, and $\ep>0$,
	\begin{equation*}\label{eq.partialgreenbounds}
	\dint_{\Omega\backslash B(Y,r)}|\nabla_Xg(X,Y)|^2\,dm(X)\quad\leq\quad\left\{\begin{matrix}Cr^{1-d},&\text{if }4r\geq\delta(Y),\\ \frac{Cr^{2-n}}{w(Y)},&\text{if }2r\leq\delta(Y),~n\geq3,\\ \frac{C_{\ep}}{w(Y)}\Big(\frac{\delta(Y)}r\Big)^{\ep},&\text{if }2r\leq\delta(Y),~n=2,\end{matrix}\right.
	\end{equation*}
	where $C>0$ depends on $d,n, C_d,C_A$, and $C_{\ep}>0$ depends on $d, C_d, C_A,\ep$.
	\item For $X,Y\in\Omega$ such that $X\neq Y$ and $\ep>0$,
	\begin{equation*}\label{eq.greenbounds}
	0\leq g(X,Y)\quad\leq\quad\left\{\begin{matrix}C|X-Y|^{1-d},&\text{if }4|X-Y|\geq\delta(Y),\\ \frac{C|X-Y|^{2-n}}{w(Y)},&\text{if }2|X-Y|\leq\delta(Y),~n\geq3,\\ \frac{C_{\ep}}{w(Y)}\Big(\frac{\delta(Y)}{|X-Y|}\Big)^{\ep},&\text{if }2|X-Y|\leq\delta(Y),~n=2,\end{matrix}\right.
	\end{equation*}
	where $C>0$ depends on $d,n, C_d,C_A$, and $C_{\ep}>0$ depends on $d, C_d, C_A,\ep$.
	\item For $q\in[1,\frac n{n-1})$ and $R\geq\delta(Y)$,
	\begin{equation*}\label{eq.partialgreenbounds2}
	\dint_{B(Y,R)}|\nabla_Xg(X,Y)|^q\,dm(X)\leq C_qR^{d(1-q)+1},
	\end{equation*}
	where $C_q>0$ depends on $d, n, C_d, C_A,$ and $q$.
	\item For $Y\in\Omega$, $R\geq\delta(Y)$, $t>0$, and $p\in[1,\frac{2n}{n-2}]$ (if $n\geq3$) or $p\in[1,+\infty)$ (if $n=2$),
	\begin{equation*}\label{eq.somebound}
	\frac{m(\{X\in B(Y,R),\,g(X,Y)>t\})}{m(B(Y,R))}\leq C_p\Big(\frac{R^{1-d}}t\Big)^{\frac p2},
	\end{equation*}
	where $C_p>0$ depends on $d,n, C_d, C_A$ and $p$.
\end{enumerate}
\end{lemma}

The next result is the representation formula given by the Green's function.

\begin{lemma}[Green representation formula, Lemma 10.7 of \cite{dfm1}] Let $g:\Omega\times\Omega\ra\bb R\cup\{+\infty\}$ be the non-negative function constructed in Lemma \ref{lm.greenfn}. Then, for any $f\in C_c^{\infty}(\Omega)$, the function $u$ defined by
	\begin{equation*}\label{eq.greenrep}
	u(X)=\int_{\Omega}g(X,Y)f(Y)\,dY
	\end{equation*}
belongs to $W_0$ and is a solution of $Lu=f$ in the sense that the identity
\begin{equation*}\label{eq.weaksense}
\int_{\Omega}A\nabla u\cdot\nabla\varphi=\int_{\Omega}\m A\nabla u\cdot\nabla\varphi\,dm=\int_{\Omega}f\varphi
\end{equation*}
holds for every $\varphi\in W_0$.
\end{lemma}

If $A$ is a matrix satisfying    (\ref{eq.elliptic}), then its transpose $A^T$ also satisfies   (\ref{eq.elliptic}). We denote $L^T=-\dv A^T\nabla$, and $g^T$ is the Green's function of Lemma \ref{lm.greenfn} for the operator $L^T$. Lemma 10.6 of \cite{dfm1} tells us that
\begin{equation}\label{eq.gisgt}
g(X,Y)=g^T(Y,X),\qquad\text{for all }X,Y\in\Omega, X\neq Y.
\end{equation}

We now use Green's functions for a representation formula concerning the difference of two solutions. A proof of it in the setting of co-dimension $1$ chord-arc domains may be found in \cite{chm} (more specifically, see their Lemma 3.12-Lemma 3.20), and its proof in our setting is essentially the same (see Remark \ref{rm.closetohm}), given that our Green's function in Lemma \ref{lm.greenfn} satisfies the properties analogous to those of the pioneering construction in \cite{gw}. Thus we omit the details of the proof, but we do provide a heuristic that formally justifies the desired identity.

\begin{lemma}[Difference of solutions, \cite{chm} Lemma 3.18]\label{lm.diffid} Suppose that $A_0, A_1$ are two matrices satisfying   (\ref{eq.elliptic}). Let $L_0=-\dv A_0\nabla$, $L_1=-\dv A_1\nabla$, and let $E$ be a Borel set in $\Gamma$. Suppose that $u_i\in W$ solves $L_iu_i=0$ in $\Omega$ and that $Tu_0=Tu_1=f\in H^{\frac12}(\Gamma)\cap C_c(\Gamma)$. Then the identity
	\begin{equation}\label{eq.diffid}
	u_1(X)-u_0(X)=\dint_{\Omega}(A_0-A_1)^T(Y)\nabla_Yg_{L_1}^T(Y,X)\nabla u_0(Y)\,dY,
	\end{equation}
	holds for almost every $X\in\Omega$, and for all $X\in\Omega\backslash\overline{\supp{(A_0-A_1)}}$.
\end{lemma}

\noindent\emph{Heuristic for the proof.} Let $F:=u_1-u_0$, and observe that $L_1F=L_1u_1-L_1u_0=-L_1u_0$. On the other hand, $L_1u_0=-\dv A_1\nabla u_0=-\dv\big((A_1-A_0)\nabla u_0)$. Therefore,
\[
\dint_{\Omega}A_1\nabla F\nabla\varphi=\dint_{\Omega}(A_0-A_1)\nabla u_0\nabla\varphi,\qquad\varphi\in C_c^{\infty}(\Omega).
\]
Equivalently,
\begin{equation}\label{eq.diffid3}
\dint_{\Omega}A_1^T\nabla\varphi\nabla F=\dint_{\Omega}(A_0-A_1)^T\nabla\varphi\nabla u_0,\qquad\varphi\in C_c^{\infty}(\Omega).
\end{equation}
Then, \emph{formally} plugging in $\varphi=g_{L_1}^T(\cdot,X)$ and using the fact that $L_1^Tg_{L_1}^T(\cdot,X)=\delta_X$, we obtain the desired result. The main issue with our logic is that in general we are not justified in plugging in $g_{L_1}^T(\cdot,X)$ for $\varphi$, because $g_{L_1}^T(\cdot,X)$ may not belong to $W_0$. However, we do point out that if $X\in\Omega\backslash\overline{\supp{(A_0-A_1)}}$ (which, incidentally, is always the situation in this paper), then we can make sense of $\varphi=g_{L^1}^T$ in the right-hand side of (\ref{eq.diffid3}), and this realization implies the claimed identity over any such $X$.\hfill{$\square$}

\subsection{The harmonic measure in a domain with boundary of high co-dimension}\label{sec.harmmeas}

In \cite{dfm1}, the Dirichlet problem
\begin{equation}\label{eq.dirichlet}
\left\{\begin{matrix}Lu=0&\text{in }\Omega,\\u=f&\text{on }\Gamma,\end{matrix}\right.
\end{equation}
was seen to have a suitably interpreted weak solution. Moreover, it was shown that there exists a family of positive regular Borel measures $\omega^X$ on $\Gamma$ indexed over $X\in\Omega$, called the \emph{harmonic measure}, such that for any boundary function $f\in C_c^0(\Gamma)$, the solution to (\ref{eq.dirichlet}) can be written as
\begin{equation}\label{eq.solnharmonicmeas}
u(X):=\int_{\Gamma}f\,d\omega^X.
\end{equation}
Here, $C_c(\Gamma)$ is the space of continuous functions on $\Gamma$ with compact support. Let us write the precise statement below. Let $C(\bb R^n)$ be the space of continuous functions on $\bb R^n$.

\begin{lemma}[Setup for harmonic measure; Lemma 9.4 of \cite{dfm1}] There exists a bounded linear operator
\begin{equation*}
U:C_c(\Gamma)\ra C(\bb R^n)
\end{equation*}
such that, for every $f\in C_c(\Gamma)$,
\begin{enumerate}[(i)]
	\item the restriction of $Uf$ to $\Gamma$ is $f$;
	\item we have that $\sup_{\bb R^n}Uf=\sup_{\Gamma}f$ and $\inf_{\bb R^n}Uf=\inf_\Gamma f$;
	\item we have that $Uf\in W_r(\Omega)$ and $Uf$ solves $L(Uf)=0$ in $\Omega$;
	\item if $B$ is a ball centered on $\Gamma$ and $f\equiv0$ in $B$, then $Uf$ lies in $W_r(B)$;
	\item if $f\in C_c(\Gamma)\cap H$, then $Uf\in W$, and $Uf$ is the unique solution of the Dirichlet problem with data $f$.
\end{enumerate}
\end{lemma}

\begin{lemma}[Harmonic measure; Lemmas 9.5 and 9.6 of \cite{dfm1}] There exists a unique positive regular Borel measure $\omega^X$ on $\Gamma$ such that
\[
Uf(X)=\int_{\Gamma}f(y)\,d\omega^X(y)
\]
for any $f\in C_c(\Gamma)$. Besides, for any Borel set $E\subset\Gamma$,
\[
\omega^X(E)=\sup\{\omega^X(K)\,:\,E\supset K, K\textup{ compact}\}=\inf\{\omega^X(V)\,:\,E\subset V, V\text{ open}\}.
\]
Moreover, for each $X\in\Omega$, $\omega^X$ is a probability measure. That is, $\omega^X(\Gamma)=1$.
\end{lemma}

We  now record some results on the harmonic measure, proved mostly in \cite{dfm1}. The first lemma below tells us qualitatively how the family of harmonic measures behaves over $X\in\Omega$.

\begin{lemma}[Lemma 9.7 of \cite{dfm1}]\label{lm.harmsoln} Let $E\subseteq\Gamma$ be a Borel set and define the function $u_E$ on $\Omega$ by $u_E(X)=\omega^X(E)$. Then
	\begin{enumerate}[(i)]
		\item \emph{(Qualitative non-degeneracy)} if there exists $X\in\Omega$ such that $u_E(X)=0$, then $u_E\equiv0$;
		\item \emph{(Harmonic meeasure is a local weak solution)}the function $u_E$ lies in $W_r(\Omega)$ and is a solution in $\Omega$;
		\item \emph{(Local trace of harmonic measure)} if $B\subseteq\bb R^n$ is a ball such that $E\cap B=\varnothing$, then $u_E\in W_r(B)$ and $Tu_E=0$ on $\Gamma\cap B$.
	\end{enumerate}
\end{lemma}

The next lemma allows us to control from below the harmonic measure on a surface ball by the Green function in certain settings.

\begin{lemma}[Green's function and the harmonic measure, Lemma 11.9 of \cite{dfm1}]\label{lm.greenlow} Let $x_0\in\Gamma$ and $r>0$ be given, and set $X_0\in\Omega$ to be a Corkscrew point for $\Delta(x_0,r)$ given by Lemma \ref{lm.corkscrew}. Then for all $X\in\Omega\backslash B(X_0,\delta(X_0)/4)$,
\begin{equation}\label{eq.greenlow1}
r^{d-1}g(X,X_0)\leq C\omega^X(B(x_0,r)\cap\Gamma),
\end{equation}
and
\begin{equation*}\label{eq.greenlow2}\notag
r^{d-1}g(X,X_0)\leq C\omega^X(\Gamma\backslash B(x_0,2r)),
\end{equation*}
where $C>0$ depends only on $d,n, C_d$ and $C_A$.
\end{lemma}

The following lemma gives non-degeneracy of the harmonic measure.

\begin{lemma}[Quantitative non-vanishing, Lemma 11.10 of \cite{dfm1}]\label{lm.nondeg} Let $\alpha>1$, $x_0\in\Gamma$, and $r>0$ be given, and let $X_0\in\Omega$ be a Corkscrew point for $\Delta(x_0,r)$. Then
\begin{gather}
\omega^X(B(x_0,r)\cap\Gamma)\geq C_{\alpha}^{-1}\quad\text{for }X\in B(x_0,r/\alpha),\label{eq.nondeg1}\notag\\  \omega^X(B(x_0,r)\cap\Gamma)\geq C_{\alpha}^{-1}\quad\text{for }X\in B(X_0,\delta(X_0)/\alpha),\label{eq.nondeg2}\\ \omega^X(\Gamma\backslash B(x_0,r))\geq C_{\alpha}^{-1}\quad\text{for }X\in\Omega\backslash B(x_0,\alpha r),\notag\label{eq.nondeg3}\\\intertext{and}\omega^X(\Gamma\backslash B(x_0,r))\geq C_{\alpha}^{-1}\quad\text{for }X\in B(X_0,\delta(X_0)/\alpha),\label{eq.nondeg4}\notag
\end{gather}
where $C_{\alpha}>0$ depends only upon $d, n, C_d, C_A$, and $\alpha$.
\end{lemma}

We may control the harmonic measure from above by the Green's function in certain settings, as the next lemma tells us.

\begin{lemma}[Lemma 11.11 of \cite{dfm1}]\label{lm.greenhigh} Let $x_0\in\Gamma$ and $r>0$ be given, and set $X_0\in\Omega$ to be a Corkscrew point for $\Delta(x_0,r)$. Then
\begin{equation}\label{eq.greenhigh1}
\omega^X(B(x_0,r)\cap\Gamma)\leq Cr^{d-1}g(X,X_0)\quad\text{for }X\in\Omega\backslash B(x_0,2r),
\end{equation}
and
\begin{equation}\label{eq.greenhigh2}\notag
\omega^X(\Gamma\backslash B(x_0,2r))\leq Cr^{d-1}g(X,X_0)\quad\text{for }X\in B(x_0,r)\backslash B(X_0,\delta(X_0)/4),
\end{equation}
where $C>0$ depends only upon $d,n, C_d,$ and $C_A$.
\end{lemma}

Next, we have a doubling property of the harmonic measure on surface balls.

\begin{lemma}[Harmonic measure is doubling, Lemma 11.12 of \cite{dfm1}]\label{lm.doublinglemma}  For $x_0\in\Gamma, r>0$, and $\alpha>1$, we have that
	\begin{gather}
	\omega^X(B(x_0,2r)\cap\Gamma)\leq C_{\doubling}\omega^X(B(x_0,r)\cap\Gamma)\quad\text{for }X\in\Omega\backslash B(x_0,2\alpha r),\label{eq.doubling1}\notag\\ \intertext{and that}\omega^X(\Gamma\backslash B(x_0,r))\leq C_{\doubling}\omega^X(\Gamma\backslash B(x_0,2r))\quad\text{for }X\in\Omega\cap B(x_0,r/\alpha),\label{eq.doubling2}\notag
	\end{gather}
	where $C_{\doubling}>0$ depends only on $n,d, C_d, C_A$, and $\alpha$.
\end{lemma}

The doubling property of the harmonic measure and the elementary properties of the dyadic cubes gives us the following corollaries.

\begin{corollary}\label{cor.likesurfaceball} Let $Q\in\bb D$, and recall that $\Delta_Q=\Delta(x_Q,a_0\ell(Q))$ is the surface ball which $Q$ contains. Then,
	\begin{equation}\label{eq.harmcube}\notag
	\omega^X(\Delta_Q)\leq\omega^X(Q)\leq C_{\doubling}^{1+\log_2(\frac{A_0}{a_0})}\omega^X(\Delta_Q),\qquad\text{for each }X\in\Omega\backslash B(x_Q,2A_0\ell(Q)).
	\end{equation}
\end{corollary}

\begin{corollary}[Harmonic measure is dyadically doubling]\label{cor.dd} Fix $Q_0\in\bb D$ and $X_0\in\Omega\backslash B(x_{Q_0},3A_0\ell(Q_0))$. Then $\omega^{X_0}$ is a dyadically doubling measure in $Q_0$.
\end{corollary}

\noindent\emph{Proof.} This result follows easily from Lemma \ref{lm.dimpliesdd} and Lemma \ref{lm.doublinglemma}.

The following notion is fundamental in our analysis of the absolute continuity of the harmonic measure.

\begin{definition}[Poisson kernel]\label{def.poisson} Fix $X\in\Omega$, and suppose that $\omega^X_L\ll\sigma$. Then we denote by $k^X_L=\frac{d\omega^X_L}{d\sigma}$ the Radon-Nikodym derivative of $\omega^X_L$ with respect to $\sigma$, and refer to it as the \emph{Poisson kernel}.
\end{definition}

We will concern ourselves with the quantitative absolute continuity of the harmonic measure, but first we have to adapt the definitions of $A_{\infty}$ and $RH$ for it to be meaningful for the harmonic measures (as families of probability measures) that we consider in this article.

\begin{definition}[$A_{\infty}$, $RH_p$, and dyadic analogues for harmonic measure]\label{def.ainftyw} We say that the harmonic measure $\{\omega^X\}_{X\in\Omega}$ is \emph{of class $A_{\infty}$} with respect to the surface measure $\sigma$, or simply $\omega\in A_{\infty}(\sigma)$, if for every $\ep>0$, there exists $\xi=\xi(\ep)>0$ such that for any surface ball $\Delta$,  every surface ball $\Delta'\subseteq\Delta$, and every Borel set $E\subset\Delta'$, we have that
	\begin{equation}\label{eq.ainftycondw}\notag
	\frac{\sigma(E)}{\sigma(\Delta')}<\xi\implies\frac{\omega^{X_\Delta}(E)}{\omega^{X_\Delta}(\Delta')}<\ep,
	\end{equation}
	where $X_\Delta$ is a Corkscrew point for $\Delta$ as in Definition \ref{def.corkscrew}. Analogously, we say that $\omega\in A_{\infty}^{\dyadic}$ if for each $Q_0\in\bb D$ and $X_{Q_0}$ a Corkscrew point for $Q_0$, we have that $\omega^{X_{Q_0}}\in A_{\infty}^{\dyadic}(Q_0)$ with uniform constants.
	
	 Given $p\in(1,\infty)$, if $\omega\ll\sigma$, then we say that  $\{\frac{d\omega^X}{d\sigma}\}_{X\in\Omega}$ is \emph{of class $RH_p$}, or simply $k=\frac{d\omega}{d\sigma}\in RH_p$, if there exists a constant $C_0\geq1$ such that for each surface ball $\Delta=\Gamma\cap B(x,r)$ with a Corkscrew point $X_{\Delta}\in\Omega$, we have the estimate
	 \begin{equation}\label{eq.poissonrh}
	 \Big(\frac1{\sigma(\Delta')}\int_{\Delta'}(k^{X_{\Delta}})^p\,d\sigma\Big)^{1/p}\leq C_0\frac1{\sigma(\Delta')}\int_{\Delta'}k^{X_{\Delta}}\,d\sigma,\qquad\text{for each surface ball }\Delta'\subseteq\Delta.
	 \end{equation}
	 We call $C_0$ the $RH_p$ characteristic of  $\frac{d\omega}{d\sigma}$. Analogously, if $\omega\ll\sigma$, we say that $k=\frac{d\omega}{d\sigma}\in RH_p^{\dyadic}$ if for each $Q_0\in\bb D$ and $X_{Q_0}$ a Corkscrew point for $Q_0$, we have that $k^{X_{Q_0}}\in RH_p^{\dyadic}(Q_0)$ with uniform $RH_p$ characteristic (see Definition \ref{def.rh}).
\end{definition}

Next we state a global comparison principle for the harmonic measure.

\begin{lemma}[Change of poles, Lemma 11.16 of \cite{dfm1}]\label{thm.cphm}  Let $x_0\in\Gamma$ and $r>0$ be given, and let $X_0\in\Omega$   be a Corkscrew point for $\Delta(x_0,r)$. Let $E,F\subseteq\Delta_0:=B(x_0,r)\cap\Gamma$ be two Borel subsets of $\Gamma$ such that both $\omega^{X_0}(E)$ and $\omega^{X_0}(F)$ are positive. Then
\begin{equation}\label{eq.cphm1}\notag
C^{-1}\frac{\omega^{X_0}(E)}{\omega^{X_0}(F)}\leq\frac{\omega^{X}(E)}{\omega^{X}(F)}\leq C\frac{\omega^{X_0}(E)}{\omega^{X_0}(F)},\qquad\text{for }X\in\Omega\backslash B(x_0,2r),
\end{equation}
where $C>0$ depends only on $n, d, C_d,$ and $C_A$. In particular, with the choice $F=B(x_0,r)\cap\Gamma$,
\begin{equation}\label{eq.cphm2}
C^{-1}\omega^{X_0}(E)\leq\frac{\omega^X(E)}{\omega^X(\Delta_0)}\leq C\omega^{X_0}(E)\qquad\text{for }X\in\Omega\backslash B(x_0,2r),
\end{equation}
where again $C>0$ depends only on $n, d, C_d$, and $C_A$.
\end{lemma}

We will also need to use a comparison principle for locally-defined solutions. The precise statement is as follows.

\begin{theorem}[Local comparison principle, Theorem 11.17 of \cite{dfm1}]\label{thm.localcp}  Let $x_0\in\Gamma$ and $r>0$ and let $X_0\in\Omega$ be a Corkscrew point for $\Delta(x_0,r)$. Let $u,v\in W_r(B(x_0,2r))$ be two non-negative, not identically zero, solutions of $Lu=Lv=0$ in $B(x_0,2r)$, such that $Tu=Tv=0$ on $\Gamma\cap B(x_0,2r)$. Then
	\begin{equation}\label{eq.localcp1}\notag
	C^{-1}\frac{u(X_0)}{v(X_0)}\leq\frac{u(X)}{v(X)}\leq C\frac{u(X_0)}{v(X_0)}\quad\text{for }X\in\Omega\cap B(x_0,r),
	\end{equation}
	where $C>0$ depends only on $n,d, C_d$ and $C_A$.
\end{theorem}

Let us show that, much like for the harmonic measure, we also have a change of poles for the Poisson kernel.

\begin{lemma}[Change of Poles for Poisson kernel]\label{lm.poissonpole} Let $\Delta\subset\Delta_0\subset\Gamma$ be surface balls in $\Gamma$, and set $X_0, X$ to be Corskcrew points of $\Delta_0,\Delta$ respectively. If $\omega\ll\sigma$ then
	\begin{equation}\label{eq.changeofpoleequiv}\notag
	k^{X}(y)\approx\frac{k^{X_0}(y)}{\omega^{X_0}(\Delta)},\qquad\text{for } \sigma-\text{a.e. }y\in\Delta.
	\end{equation}
\end{lemma}

\noindent\emph{Proof.} Write $\Delta_0=B(x_0,r_0)\cap\Gamma$ and $\Delta=B(x,r)\cap\Gamma$. Let $X_0'$ be a Corkscrew point for $4\Delta_0$. Then $X_0'\notin B(x_0,2r_0)$, and hence $X_0'\notin B(x,2r)$. By the Harnack chains we have that
\[
\omega^{X_0'}(E)\approx\omega^{X_0}(E),\qquad\text{for any Borel set }E\subset\Delta.
\]
We now apply (\ref{eq.cphm2}) to obtain that for any Borel set $E\subset\Delta$,
\[
\omega^X(E)\approx\frac{\omega^{X_0'}(E)}{\omega^{X_0'}(\Delta)}\approx\frac{\omega^{X_0}(E)}{\omega^{X_0}(\Delta)}.
\]
The desired result now follows by the Lebesgue differentiation theorem and letting $E\searrow y\in Q$.\hfill{$\square$}

The next lemma collects the results which allow us to compare harmonic measures (and Green functions, incidentally) for operators which agree locally near a surface ball. We remark in passing that the proof shown is based on the local comparison principle stated above, but it is also possible to obtain the results  in the following lemma without appealing to the local comparison principle, by means of the identity (\ref{eq.diffid}), the properties of the Green function, and the Caccioppoli inequality.

\begin{lemma}[Comparison of harmonic measures near the boundary]\label{lm.poissonequiv} Fix $x\in\Gamma$, $r>0$, let $X_0$ be a Corkscrew point (with Corkscrew constant $c<1$) for the surface ball $\Delta_0:=\Delta(x,r)$, and suppose that $A_0$ and $A_1$ are two matrices satisfying (\ref{eq.elliptic}) and $A_0\equiv A_1$ in $B(x,4c^{-1}r)\cap\Omega$. Let $L_0=-\dv A_0\nabla$ and $L_1=-\dv A_1\nabla$.  The following statements hold.
\begin{enumerate}[(i)]
	\item For each surface ball $\Delta'\subset\Delta_0$, we have that
\begin{equation}\label{eq.equiv1}
\tfrac1{C}\omega_1^{X_0}(\Delta')\leq\omega_0^{X_0}(\Delta')\leq C\omega_1^{X_0}(\Delta').
\end{equation}
\item The measures $\omega_1^{X_0}$ and $\omega_0^{X_0}$ are mutually absolutely continuous on $\Delta_0$.
\item If $\omega_0^{X_0}|_{\Delta_0}\ll\sigma|_{\Delta_0}$, then $\omega_1^{X_0}|_{\Delta_0}\ll\sigma|_{\Delta_0}$, and $k_0^{X_0}(y)\approx k_1^{X_0}(y)$, for  $\sigma-$a.e. $y\in\Delta_0$.
\end{enumerate}
\end{lemma}

\noindent\emph{Proof.} (i). Let $\tilde X_0\in\Omega$ be a Corkscrew point for $\Delta(x,4c^{-1}r)$, so that $\tilde X_0\in\Omega\backslash B(x,4r)$. Note that since $L_0\equiv L_1$ in $B(x,4r)\cap\Omega$, then $A_0^T\equiv A_1^T$ in $B(x,4r)\cap\Omega$. As such, we may apply Theorem \ref{thm.localcp} to the Green functions $g_0^T(\cdot,\tilde X_0)$ and $g_1^T(\cdot,\tilde X_0)$, to deduce that
\begin{equation}\label{eq.localcpgreens}
\frac{g_0^T(X_0 ,\tilde X_0)}{g_1^T(X_0  ,\tilde X_0)}\approx\frac{g_0^T(Y,\tilde X_0)}{g_1^T(Y,\tilde X_0)},\qquad\text{for every }Y\in B(x,r)\cap\Omega.
\end{equation}
We now use (\ref{eq.gisgt}), (\ref{eq.greenlow1}), and (\ref{eq.greenhigh1})  to obtain that
\begin{equation}\label{eq.greeneq}\notag
g_i^T(X_0,\tilde X_0)=g_i(\tilde X_0,X_0)\approx\frac{\omega_i^{\tilde X_0}(\Delta_0)}{r^{d-1}},\qquad i=0,1.
\end{equation}
By the Harnack inequality and Harnack chains we have that
\[
\omega_i^{\tilde X_0}(\Delta_0)\approx\omega_i^{X_0}(\Delta_0)\approx1,\qquad i=0,1,
\]
and thus using these last results in (\ref{eq.localcpgreens}), we see that
\begin{equation}\label{eq.localcpgreens2}
g_0^T(Y,\tilde X_0)\approx g_1^T(Y,\tilde X_0),\qquad\text{for every }Y\in B(x,r)\cap\Omega.
\end{equation}

Now fix a surface ball $\Delta'=\Delta(y,r')\subset\Delta_0$ and let $Y'$ be a Corkscrew point for $\Delta'$. Then $Y'\in B(x,r)\cap\Omega$. Since we may write
\[
g_i^T(Y',\tilde X_0)=g_i(\tilde X_0,Y')\approx\frac{\omega_i^{\tilde X_0}(\Delta')}{{r'}^{d-1}},\qquad i=0,1,
\]
then, using (\ref{eq.localcpgreens2}), we observe that $\omega_0^{\tilde X_0}(\Delta')\approx\omega_1^{\tilde X_0}(\Delta')$, and (\ref{eq.equiv1}) immediately follows.

(ii). Let us see that (\ref{eq.equiv1}) implies the mutual absolute continuity of $\omega_1^{X_{\Delta}}$ and $\omega_0^{X_{\Delta}}$, it suffices to use the (outer and inner) regularity of the measures, the Besicovitch covering theorem \cite{besicovitch} applied to a bounded open set, and (\ref{eq.equiv1}). More precisely, let $E\subseteq\Delta_0$ be a Borel set such that $\omega_0^{X_0}(E)>0$. Then by the inner regularity of $\omega_0^{X_0}$, there exists a compact set $K\subseteq E$ such that $\omega_0^{X_0}(K)\approx\omega_0^{X_0}(E)$. To prove that $\omega_1^{X_0}(E)>0$, it suffices to show that $\omega_1^{X_0}(K)>0$. Let $V\subseteq\Delta_0$ be an open set in the subspace topology of $\Delta_0$ such that $V\supset K$. Then we can write $V=\Gamma\cap\big(\medcup_{x\in V}B(x,r_x)\big)$, for suitable finite $r_x>0$. Let $\m B:=\{B(x,r_x)\}$. The latter is a Besicovitch covering of the bounded set $V$, and thus applying the Besicovitch covering theorem, we can write $V=\Gamma\cap\Big(\cup_{B\in\m B'}B\Big)$ where $\m B'$ is a subcollection of $\m B$ such that the balls intersect an at most uniformly finite (depending only on $n$) number of times. Then we may use (\ref{eq.equiv1}) to estimate $\omega_0^{X_0}(K)$ from above by $\omega_1^{X_0}(V)$ times a constant independent of $V$. Since this is true for any $V$, it follows by the outer regularity of $\omega_1^{X_0}$ that $\omega_1^{X_0}(K)>0$, which completes the proof of (ii).

(iii). Since $\omega_0^{X_0}|_{\Delta_0}\ll\sigma|_{\Delta_0}$ and we have seen that (ii) holds, then $\omega_1^{X_0}|_{\Delta_0}\ll\sigma|_{\Delta_0}$ follows. Now fix $y\in\Delta_0$, and for each $k\in\bb N$, let $\Delta_k=\Delta(y,r_k)\subset\Delta_0$ with $r_k\searrow0$ as $k\ra\infty$. According to (i), we may then write
\[
\frac{\omega_0^{X_0}(\Delta_k)}{\sigma(\Delta_k)}\approx\frac{\omega_1^{X_0}(\Delta_k)}{\sigma(\Delta_k)},\qquad\text{for each }k\in\bb N\text{ and every }y\in_0.
\]
Finally, we send $k\ra\infty$ in the above estimate, and due to the Lebesgue Differentiation Theorem we arrive at the desired result.\hfill{$\square$}

The next comparison principle between Poisson kernels is similar in spirit to the one shown in the previous lemma, but adapted to our dyadic grid.

\begin{lemma}[Comparison of Poisson Kernels in a cube, Lemma 3.24 in \cite{chm}]\label{lm.poissoncube} \ \\ Fix $Q_0\in\bb D$, let $X_0\in\Omega$ be a Corkscrew point (with Corkscrew constant $c$) for the surface ball $\Delta(x_{Q_0},10c^{-1}\sqrt nA_2\ell(Q_0))$, let $\m F\subset Q_0$ be a disjoint family, and suppose that $A_0$ and $A_1$ are two matrices satisfying (\ref{eq.elliptic}) and $A_0\equiv A_1$ in $R_{Q_0}\backslash(\Omega_{\m F}\backslash\Omega_{\m F,Q_0})$. Let $L_0=-\dv A_0\nabla$ and $L_1=-\dv A_1\nabla$.  If the corresponding harmonic measures $\omega_0^{X_0},\omega_1^{X_0}$ are absolutely continuous with respect to $\sigma$, then for each $t\in(0,1)$ we have that
\[
\frac1{C_t}k_1^{X_0}(y)\leq k_0^{X_0}(y)\leq C_{t}k_1^{X_0}(y),\qquad\text{for }\sigma-\text{almost every }y\in Q_0\backslash\Sigma_{Q_0,t},
\]
and $\Sigma_{Q_0,t}:=\big\{x\in Q_0:\dist(x,\Gamma\backslash Q_0)<t\ell(Q_0)\}$.
\end{lemma}

\noindent\emph{Sketch of proof.} The proof is essentially the same as in Lemma 3.24 of \cite{chm} (see Remark \ref{rm.closetohm}), except that our claim is slightly sharper by requiring that $A_0\equiv A_1$ only in $R_{Q_0}\backslash(\Omega_{\m F}\backslash\Omega_{\m F,Q_0})$ as opposed to in all of $R_{Q_0}$. It turns out that this is enough though: in \cite{chm}, the authors cover $Q_0\backslash\Sigma_{Q_0,t}$ with a uniformly finite (cardinality depending on $t$) collection of surface balls $\{\Delta_k=\Delta(x_k,r_k)\}_k$, and chosen in such a way that $x_k$ and $r_k\approx_{t}\ell(Q_0)$ verify the containment $B_k\cap\Omega=B(x_k,Cr_k)\cap\Omega\subseteq R_{Q_0}$, for some uniform large constant $C\approx1$. Actually, their method of proof gives that if $X\in B_k\cap\Omega$ and $X\in I\in\m W$, then $I\in\m W_Q$ with $Q\in\bb D_{Q_0}$. Hence, if $X\in\Omega_{\m F}\cap B_k$, then $X\in\Omega_{\m F,Q_0}$. It follows that $(\Omega_{\m F}\backslash\Omega_{\m F,Q_0})\cap B_k=\varnothing$, whence we have that $B_k\subset R_{Q_0}\backslash(\Omega_{\m F}\backslash\Omega_{\m F,Q_0})$. The rest of the proof is elementary; one employs Lemma \ref{lm.poissonequiv} and the fact that $L_0=L_1$ in $B_k$ to get the desired result on each surface ball $\Delta_k$, and via Harnack Chains and the Harnack Inequality, one can teleport (with constants depending on $t$) from a Corkscrew point for $\Delta_k$ to the fixed point $X_0$.\hfill{$\square$}

\begin{remark}\label{rm.unbounded} The main reason why we require the slightly sharper version of this result as opposed to in \cite{hm1}, \cite{chm}, is because we will  decide to use our analogue of the Dahlberg-Jerison-Kenig sawtooth lemma, Lemma \ref{lm.djkg}, on the unbounded sawtooth domain $\Omega_{\m F}$ as opposed to the bounded sawtooth $\Omega_{\m F,Q_0}$ whose mixed-dimension elliptic theory we have not fully developed (although it would not be hard to make it work given our theory in this paper and in \cite{dfm20}; it would just be tedious rather than difficult).
\end{remark}

We next see how to relate  the solvability of the Dirichlet problem with the quantitative absolute continuity of the harmonic measure.

\begin{theorem}[Relationship between $A_{\infty}$ and the Dirichlet problem]\label{thm.equivalence} Assume that $\Gamma$ is a closed \dADR set with $d\in[1,n-1)$ not necessarily an integer. Suppose that the matrix $A$ satisfies (\ref{eq.elliptic}), let $L=-\dv A\nabla$, let $\omega$ be the harmonic measure associated to $L$, and let $p,p'\in(1,\infty)$, $\frac1p+\frac1{p'}=1$. Then, the following statements are equivalent:
	\begin{enumerate}[(a)]
		\item For each $f\in C_c(\Gamma)$, the solution to the Dirichlet problem $u$ satisfies the non-tangential bound
		\begin{equation}\label{eq.nontangentialbound}
		\Vert Nu\Vert_{L^{p'}(\Gamma)}\leq C\Vert f\Vert_{L^{p'}(\Gamma)},
		\end{equation}
		where $Nu$ is the non-tangential maximal function and $C$ is a uniform constant.
		\item We have that $\omega\ll\sigma$ and $\frac{d\omega}{d\sigma}\in RH_p$. 
		\item We have that $\omega\ll\sigma$, and there is a uniform constant $C_0$ such that for every surface ball $\Delta=\Gamma\cap B(x,r)$, there exists $X_{\Delta}\in\Omega$, which is a Corkscrew point for $\Delta$, verifying the  following scale-invariant $L^p$ estimate:
		\begin{equation}\label{eq.poissonlpest}
		\int_{\Delta}(k^{X_{\Delta}})^p\,d\sigma\leq C_0\sigma(\Delta)^{1-p}.
		\end{equation}
	\end{enumerate}
\end{theorem}

\noindent\emph{Proof.} (a)$\implies$(b). Fix a surface ball $\Delta=\Gamma\cap B(x_0,r)$ and $X_{\Delta}\in\Omega$ a Corkscrew point for the surface ball $\Delta$. Let $X_0\in\Omega\backslash B(x_0,2r)$ be a Corkscrew point for the surface ball $4\Delta$, and immediately by Harnack Chains and the Harnack Inequality we see that $\omega^{X_0}\approx\omega^{X_{\Delta}}$, whence we need only prove the desired result with pole $X_0$ as opposed to $X_{\Delta}$.  We will show that $\omega^{X_0}\ll\sigma$ on $\Delta$ and that $\frac{d\omega^{X_0}}{d\sigma}\in RH_p(\sigma,\Delta)$ via the characterization in Theorem \ref{thm.ainftyprop} (vi).  Owing to Lemma \ref{lm.harmsoln} (i), we have that for any $X\in\Omega$, $\omega^X\ll\omega^{X_0}$.  Accordingly, for each $X\in\Omega$ let $\n K(X,\cdot)=\frac{d\omega^X}{d\omega^{X_0}}$ be the Radon-Nikodym derivative (see \cite{folland}). Since $\{\omega^X\}$ is a family of probability measures, we trivially have that $\n K(X,\cdot)\in L^1(\Gamma,\omega^{X_0})$, for each $X\in\Omega$.

Now fix a non-negative $f\in C_c(\Delta)$, and let
\[
u(Y):=\int_{\Gamma}f\,d\omega^Y,\qquad\text{for each }Y\in\Omega.
\]
We claim that there exists a uniform constant $C>0$ such that
\begin{multline}\label{eq.majorized}
(\m M_{\omega^{X_0}}f)(x):=\sup_{\Delta(x,s)\subseteq\Delta}~\frac1{\omega^{X_0}(\Delta(x,s))}\int_{\Delta(x,s)}f\,d\omega^{X_0}\\ \leq C(Nu)(x),\qquad\text{for each }x\in\Delta.
\end{multline}
Assume the claim for a moment. Then, since $u$ is the solution to the Dirichlet problem with data $f$ and (\ref{eq.nontangentialbound}) holds, we have the estimate
\[
\Vert\m M_{\omega^{X_0}}f\Vert_{L^{p'}(\Delta,\sigma)}\lesssim\Vert Nu\Vert_{L^{p'}(\Delta,\sigma)}\leq\Vert Nu\Vert_{L^{p'}(\Gamma,\sigma)}\lesssim\Vert f\Vert_{L^{p'}(\Gamma,\sigma)}=\Vert f\Vert_{L^{p'}(\Delta,\sigma)},
\]
valid for each non-negative $f\in C_c(\Delta)$. If $f\in L^{p'}(\Delta,\sigma)$ is non-negative, we may approximate it by $C_c(\Delta)$ non-negative functions in a standard way, so that the estimate
\[
\Vert\m M_{\omega^{X_0}}f\Vert_{L^{p'}(\Delta,\sigma)}\lesssim\Vert f\Vert_{L^{p'}(\Delta,\sigma)}
\]
is valid for all non-negative $f\in L^{p'}(\Delta,\sigma)$. Consequently, according to Theorem \ref{thm.ainftyprop} (vi), we deduce that $\frac{d\omega^{X_0}}{d\sigma}\in RH_p(\Delta)$ with $RH_p$ characteristic independent of $\Delta$, as desired.

Thus we proceed to prove (\ref{eq.majorized}). Fix $x\in\Delta$, $\alpha>0$, and $X\in\gamma^{\alpha}(x)$ such that $s:=|X-x|$ satisfies $\Delta(x,s)\subseteq\Delta$. By definition of $\gamma^{\alpha}(x)$, we have that $\delta(X)\leq s\leq(1+\alpha)\delta(X)$. Observe that by the non-negativity of $f$ and the properties of the Radon-Nikodym derivative, it is the case that
\begin{equation}\label{eq.aimpliesb1}
u(X)\geq\int_{\Delta(x,s)}f\,d\omega^X=\int_{\Delta(x,s)}f\n K(X,\cdot)\,d\omega^{X_0}.
\end{equation}
Since $\omega^{X_0}$ is a doubling measure on $\Delta$ (see Lemma \ref{lm.doublinglemma}), we may use the Lebesgue Differentiation Theorem for doubling measures \cite{folland} to obtain that for $\omega^{X_0}-\text{a.e. }y\in\Delta(x,s)$,
\[
\n K(X,y)=\lim_{\Delta'\searrow y}\frac1{\omega^{X_0}(\Delta')}\int_{\Delta'}\n K(X,\cdot)\,d\omega^{X_0}=\lim_{\Delta'\searrow y}\frac{\omega^X(\Delta')}{\omega^{X_0}(\Delta')},
\]
(here, $\Delta'$ is a surface ball centered at $y$ and contained in $\Delta(x,s)$). Note that necessarily we have $X_0\in\Omega\backslash B(x,2s)$.  Denote by $A_{\Delta(x,s)}$ a Corkscrew point for $\Delta(x,s)$, and so from the Comparison Principle (\ref{eq.cphm2}) we may conclude that
\[
\frac{\omega^{A_{\Delta(x,s)}}(\Delta')}{\omega^{X_0}(\Delta')}\approx\frac1{\omega^{X_0}(\Delta(x,s))},\qquad\text{for all }\Delta'\searrow y~\text{and all}~\Delta(x,s)\subseteq\Delta.
\]
On the other hand, for any $y\in\Delta(x,s)$, $\frac1{1+\alpha}s\leq|X-y|\leq2s$, which implies by the Harnack chains that
\[
\omega^X(\Delta')\approx\omega^{A_{\Delta(x,s)}}(\Delta').
\]
Putting all these observations together and back into (\ref{eq.aimpliesb1}), we deduce that
\[
(Nu)(x)\gtrsim\frac1{\omega^{X_0}(\Delta(x,s))}\int_{\Delta(x,s)}f\,d\omega^{X_0},
\]
for each $x\in\Delta$ and each $s>0$ such that $\Delta(x,s)\subseteq\Delta$. Since $\Delta(x,s)\subseteq\Delta$ is arbitrary, the claim (\ref{eq.majorized}) follows.

(b)$\implies$(a). This is a consequence of Theorem 4.1 in \cite{mz} (formally, they have symmetric $A$, but this assumption can be dropped).

(b)$\implies$(c). By assumption we already have that $\omega\ll\sigma$. Now fix $\Delta\subset\Gamma$ and $X_{\Delta}$ the Corkscrew point for $\Delta$ given by property (b) above. We apply (\ref{eq.poissonrh}) with $\Delta'=\Delta$ and raise the inequality to the $p-$th power, to get that
\[
\int_{\Delta}(k^{X_\Delta})^p\,d\sigma\leq C_0^p\sigma(\Delta)^{1-p}\omega^{X_\Delta}(\Delta)^p\leq C\sigma(\Delta)^{1-p},
\]
where we used the non-degeneracy of the harmonic measure in the last inequality. Hence (\ref{eq.poissonlpest}) is established.

(c)$\implies$(b). By assumption we already have that $\omega\ll\sigma$. Now fix a surface ball $\Delta\subset\Gamma$ and another surface ball $\Delta'\subseteq\Delta$. According to (\ref{eq.poissonlpest}) we have the estimate
\[
\Big(\frac1{\sigma(\Delta')}\int_{\Delta'}(k^{X_{\Delta'}})^p\,d\sigma\Big)^{1/p}\leq C_0^{\frac1p}\frac1{\sigma(\Delta')}.
\]
Next, we use Lemma \ref{lm.poissonpole} applied with surface balls $\Delta,\Delta'$ to obtain that there exists a uniform (in $\Delta,\Delta'$) constant $\tilde c$ such that
\[
k^{X_{\Delta'}}\geq\tilde c\frac{k^{X_\Delta}}{\omega^{X_\Delta}(\Delta')},\qquad\sigma-\text{a.e. on }\Delta'.
\]
Putting these observations together, we deduce that
\[
\Big(\frac1{\sigma(\Delta')}\int_{\Delta'}(k^{X_{\Delta}})^p\,d\sigma\Big)^{1/p}\leq\frac{C_0^{\frac1p}}{\tilde c}\frac1{\sigma(\Delta')}\omega^{X_\Delta}(\Delta')=\frac{C_0^{\frac1p}}{\tilde c}\frac1{\sigma(\Delta')}\int_{\Delta'}k^{X_\Delta}\,d\sigma,
\]
as desired.\hfill{$\square$}

When (a) of the above theorem occurs we say that $(\operatorname{D})_{p'}$ is solvable for $L$ or that $L$ is solvable in $L^{p'}$. In such case, for every $f\in L^{p'}(\bb R^d)$ there exists a unique $u$ such that $Lu=0$ in $\bb R^n$, (\ref{eq.nontangentialbound}) holds and $u$ converges non-tangentially to $f$ for $\sigma-$almost every $x\in\Gamma$.

\begin{remark}[Equivalence of $RH_p$ and $RH_p^{\dyadic}$]\label{rm.equivrh} Suppose throughout this remark that $\omega\ll\sigma$. In Theorem \ref{thm.equivalence}, we saw that the condition (\ref{eq.poissonlpest}) is equivalent to $\frac{d\omega}{d\sigma}\in RH_p$. Consider the following dyadic analogue of condition (\ref{eq.poissonlpest}): For each $Q\in\bb D$ and for $X_Q$ a Corkscrew point for $Q$, the estimate
\begin{equation}\label{eq.dyadicanalogue}
\int_Q(k^{X_Q})^p\,d\sigma\lesssim\sigma(Q)^{1-p}
\end{equation}
holds. Using Corollary \ref{cor.likesurfaceball} and the Harnack inequality to flexibly move the poles, it is not difficult to see that the condition (\ref{eq.dyadicanalogue}) is equivalent to (\ref{eq.poissonlpest}). Of course, this new condition (\ref{eq.dyadicanalogue}) is also equivalent to the condition that $\frac{d\omega}{d\sigma}\in RH_p^{\dyadic}$. It follows that $\frac{d\omega}{d\sigma}\in RH_p$ is equivalent to $\frac{d\omega}{d\sigma}\in RH_p^{\dyadic}$.
\end{remark}

To end this section, we record the $L^p$-control of the square function by the nontangential maximal function under the assumption that the harmonic measure lies in $A_{\infty}$.

\begin{theorem}[$S<N$; Theorem 3.1 of \cite{mz}]\label{thm.djk} Suppose that $d\in[1,n-1)$ is not necessarily an integer, $\Gamma$ is \dADR, and that $A$ is a (not necessarily symmetric, see Remark \ref{rm.nonsym} below) matrix satisfying (\ref{eq.elliptic}). Assume that for some $p'\in(1,\infty), (D)_{p'}$ is solvable for $L$. Write $u$ for the solution to the Dirichlet problem $\operatorname{(D)}_{p'}$ with data $f\in L^{p'}(\Gamma)$. Then, for all apertures $\alpha>0$ we have the estimate
\[
\Vert Su\Vert_{L^{p'}(\Gamma)}\lesssim\Vert f\Vert_{L^{p'}(\Gamma)},
\]
where the implicit constants depend on $n$, $d$, $C_d$, $C_A$, $\alpha$, and the $RH_p$ constant of $\omega$.
\end{theorem}

\begin{remark}\label{rm.nonsym} We remark that the previous theorem is stated in \cite{mz} for symmetric matrices and $d\in\bb N$ only, but in fact their method of proof generalizes to non-symmetric matrices, mainly using (\ref{eq.gisgt}), and to all $d\in\bb R$, $d\in[1,n-1)$. For concreteness, the fact that $d\in\bb N$ was never explicitly used in the proof (recall that the construction of the dyadic cubes Lemma \ref{lm.dyadiccubes} works for all real $d\in(0,n-1)$), and  the symmetry of the matrix $A$ is used only in Step 3 of the proof of their Proposition 1.16, and explicitly arising only in their calculation (3.78), where $A=A^T$ is used to maneuver the integration by parts. However, we note that their function $G=G(X)$ in (3.78) really is $G(X_Q,X)$ (see their estimate (3.73)), while $G(X_Q,X)=G^T(X,X_Q)$ by (\ref{eq.gisgt}), and the latter is the  ``correct'' Green's function for which the needed cancellation $\dv A^TG^T(\cdot,X_Q)=0$ will hold in their (3.78) (of course, in the symmetric setting, there is no difference between these Green's functions). Thus, in the non-symmetric setting, the last two lines of their calculation (3.78) read 
\begin{multline*}
\ldots=\frac12\dint_{\bb R^n}u^2(X)A^T(X)\nabla_XG(X_Q,X)\nabla\psi_N(X)\,dX\\-\dint_{\bb R^n}u(X)G(X_Q,X)A(X)\nabla u(X)\nabla\psi_N(X)\,dX =:\frac12I-II.
\end{multline*}
One then uses the representation $G(X_Q,X)=G^T(X,X_Q)$ again while bounding $|I|$ and $|II|$ (see their (3.79)-(3.81)) to exploit the fact that $G^T(\cdot,X_Q)$ solves $-\dv A^TG^(\cdot,X_Q)=0$ in the fat sawtooth domain, so that one may use the Caccioppoli inequality and the Harnack inequality as required. At the last step when bounding $|I|$ and $|II|$, we switch from $G^T(X_I,X_Q)$  to $G(X_Q,X_I)$ and invoke Lemma \ref{lm.greenlow} to obtain the required control \emph{with $\omega_L$} (and no dependence on   $\omega_{L^T}$). As explained here, one recovers their full Proposition 1.16 in the non-symmetric case; the rest of the proof of their Theorem 3.1 sees no obstacle from the non-symmetric point of view.
\end{remark}

	\begin{remark}\label{rm.dreal} We note that Theorem \ref{thm.djk} has content for all $d\in[1,n-2)$ with $d$ not necessarily an integer, and for all closed \dADR unbounded $\Gamma$. Indeed, for any such $\Gamma$, we consider the special operator $L_{\operatorname{DEM}}$ of \cite{dem} Theorem 6.7, and recall that $\omega_{\operatorname{DEM}}\ll\sigma$. By Theorem \ref{thm.equivalence}, it follows that there exists $p>1$ such that $(\operatorname{D})_p$ is solvable for $L_{\operatorname{DEM}}$. We thus see that the hypotheses of Theorem \ref{thm.djk} are verified for this operator, and hence the control of the square function by the non-tangential maximal function holds.
	\end{remark}

\section{The  projection lemma for the dyadically-generated sawtooth}\label{sec.djk}

Having shown that the triple $(\Omega_{\m F},m,\sigma_{\star})$ satisfies the axioms \ref{ax.h1}-\ref{ax.trace} in Section \ref{sec.measure}, we may   appeal to the elliptic theory set forth in \cite{dfm20} to conclude that there exists a harmonic measure $\omega_{\star}$ on $\partial\Omega_{\m F}$ associated to the operator $L=-\dv A\nabla$ whose matrix $A$ satisfies  (\ref{eq.elliptic}). This harmonic measure $\omega_\star$ on the sawtooth boundary $\Omega_{\m F}$ enjoys many similar properties to the harmonic measure on $\Gamma$ which were reviewed in Subsection \ref{sec.harmmeas}. In particular, we note that $\omega_{\star}$ has the doubling property; that is, we have Lemma \ref{lm.doublinglemma} with $\Omega$, $\Gamma$, and $\omega$ replaced by $\Omega_{\m F}$, $\partial\Omega_{\m F}$, and $\omega_{\star}$, respectively (see Lemma 15.43 of \cite{dfm20}).

The following lemma is an analogue of the Dahlberg-Jerison-Kenig sawtooth lemma \cite{djk}; it has already been shown in \cite{dm} on a similar setting, but we include its proof here too for completeness.

\begin{lemma}[Dyadic sawtooth lemma, global version]\label{lm.djkg} Suppose that $\Gamma$ is a \dADR set with $d\in(0,n-1)$. Fix $Q_0\in\bb D$, let $\m F=\{Q_j\}_j\subset\bb D_{Q_0}$ be a family of pairwise disjoint dyadic cubes, and let $\m P_{\m F}$ be the corresponding projection operator, as in (\ref{eq.projection}). Let $X_0$ be the Corkscrew point for $Q_0$ with respect to both $\Omega$ and $\Omega_{\m F}$, whose existence is shown in Proposition \ref{prop.simultcs} and Corollary \ref{cor.bigcorkscrew}. Let $A$ be a matrix of essentially bounded, real coefficients satisfying the weighted ellipticity condition  (\ref{eq.elliptic}). Let $\omega=\omega^{X_0}$ and $\omega_{\star}=\omega_{\star}^{X_0}$ denote the respective harmonic measures for the domains $\Omega$ and $\Omega_{\m F}$, with the fixed pole $X_0$ as above. Let $\mu=\mu^{X_0}$ be the measure defined on $Q_0$ as
	\begin{equation}\label{eq.mug}
	\mu(F)=\omega_{\star}\big(F\backslash(\cup_jQ_j)\big)+\sum_{Q_j\in\m F}\frac{\omega(F\cap Q_j)}{\omega(Q_j)}\omega_{\star}(P_j),\qquad F\subset Q_0,
	\end{equation}
	where $P_j$ is the $n-$dimensional cube constructed in Proposition \ref{prop.cubeinface}. Then $\m P_{\m F}\mu$ depends only on $\omega_{\star}$ and not on $\omega$. More precisely,
	\begin{equation}\label{eq.pfmug}
	\m P_{\m F}\mu(F)=\omega_{\star}\big(F\backslash(\cup_jQ_j)\big)+\sum_{Q_j\in\m F}\frac{\sigma(F\cap Q_j)}{\sigma(Q_j)}\omega_{\star}(P_j),\qquad F\subset Q_0.
	\end{equation}
	Moreover, there exists $\theta>0$ such that for all $Q\in\bb D_{Q_0}$ and Borel sets $F\subset Q$, we have that
	\begin{equation}\label{eq.projthetag}
	\Big(\frac{\m P_{\m F}\omega(F)}{\m P_{\m F}\omega(Q)}\Big)^{\theta}\lesssim\frac{\m P_{\m F}\mu(F)}{\m P_{\m F}\mu(Q)}\lesssim\frac{\m P_{\m F}\omega(F)}{\m P_{\m F}\omega(Q)}.
	\end{equation}
\end{lemma}

\noindent\emph{Proof.}  We first verify that (\ref{eq.pfmug}) holds. Let $F\subset Q_0$ be a Borel set, and observe that
\begin{multline*}
\m P_{\m F}\mu(F)=\mu\big(F\backslash(\cup_jQ_j)\big)+\sum_{Q_j\in\m F}\frac{\sigma(F\cap Q_j)}{\sigma(Q_j)}\mu(Q_j)\\ =\omega_{\star}\big(F\backslash(\cup_jQ_j)\big)+\sum_{Q_j\in\m F}\frac{\sigma(F\cap Q_j)}{\sigma(Q_j)}\sum_{Q_k\in\m F}\frac{\omega(Q_j\cap Q_k)}{\omega(Q_k)}\omega_{\star}(P_k),
\end{multline*}
which implies (\ref{eq.pfmug}) since $\m F$ is a pairwise disjoint family.

We now show the right-hand side inequality in (\ref{eq.projthetag}), so we fix $Q\in\bb D_{Q_0}$ and $F\subset Q$. First, suppose that there exists $Q_j\in\m F$ such that $Q\subset Q_j$. In this case, we have that
\begin{equation*}
\frac{\m P_{\m F}\mu(F)}{\m P_{\m F}\mu(Q)}=\frac{\frac{\sigma(F\cap Q_j)}{\sigma(Q_j)}\omega_{\star}(P_j)}{\frac{\sigma(Q\cap Q_j)}{\sigma(Q_j)}\omega_{\star}(P_j)}=\frac{\frac{\sigma(F\cap Q_j)}{\sigma(Q_j)}\omega(Q_j)}{\frac{\sigma(Q\cap Q_j)}{\sigma(Q_j)}\omega(Q_j)}=\frac{\m P_{\m F}\omega(F)}{\m P_{\m F}\omega(Q)}.
\end{equation*}

Therefore, it remains to consider the case that $Q$ is not contained in any $Q_j$, $Q_j\in\m F$; in this case, we have that $Q\in\bb D_{\m F,Q_0}\subset\bb D_{\m F}$, and that, if $Q_j\cap Q\neq\varnothing$, then $Q_j\subsetneq Q$. Let $x_j^{\star}$ be the center of $P_j$ and let $r_j\approx\ell(Q_j)$ be as in Notation \ref{not.centerofpj}, so that $P_j\subseteq\Delta_{\star}(x_j^{\star},r_j)$. Since $\ell(P_j)\approx\ell(Q_j)\approx r_j$, we have that
\begin{equation}\label{eq.djk1}
\omega^{X_0}_\star(P_j)\geq\omega^{X_0}_\star\big(\Delta_\star(x_j^\star,\ell(P_j)/2)\big)\gtrsim\omega^{X_0}_\star\big(\Delta_\star(x_j^\star,r_j)\big),
\end{equation}
where we used the doubling property of $\omega^{X_0}_\star$ in the last estimate (if $\ell(Q_j)\approx\ell(Q)$, then we consider a point $X'\in\Omega_{\m F}$, where $X'$ is a Corkscrew point for the surface ball $\Delta(x_{Q_0},A\ell(Q_0))$, with $A>0$ large enough, and such that $\omega^{X'}_\star$ is doubling on the surface balls centered at $x_j^\star$ and of radius less than $r_j$; then, by the Harnack Inequality Lemma 11.35 and Lemma 12.19 of \cite{dfm20}, and Harnack Chains, we have that $\omega^{X_0}_\star\approx\omega^{X'}_\star$). Using (\ref{eq.pfmug}), we have that
\begin{multline}\label{eq.djk2}
\m P_{\m F}\mu(Q)=\omega_{\star}\big(Q\backslash(\cup_jQ_j)\big)+\sum_{Q_j\in\m F, Q_j\subsetneq Q}\frac{\sigma(Q\cap Q_j)}{\sigma(Q_j)}\omega_{\star}(P_j)\\=\omega_{\star}\big(Q\backslash(\cup_jQ_j)\big)+\sum_{Q_j\in\m F, Q_j\subsetneq Q}\omega_{\star}(P_j)\\ \gtrsim\omega_{\star}\big(Q\backslash(\cup_jQ_j)\big)+\sum_{Q_j\in\m F, Q_j\subsetneq Q}\omega_{\star}\big(\Delta_\star(x_j^\star,r_j)\big)\\ \geq\omega_\star\Big(\big(Q\cap\partial\Omega_{\m F}\big)\medcup\big(\medcup_{Q_j\in\m F, Q_j\subsetneq Q}\Delta_\star(x_j^\star,r_j)\big)\Big)\geq \omega_\star\big(\Delta_\star^Q\big),
\end{multline}
where $\Delta_\star^Q$ is the surface ball in Proposition \ref{prop.starball}, in the third line we used (\ref{eq.djk1}), and in the last line we used Proposition \ref{prop.starball}, and the doubling property of $\omega_\star$ and Propositions \ref{prop.containment} and \ref{prop.mumeasure0} to see that $\omega_\star(Q\backslash(\cup_jQ_j))=\omega_\star(Q\cap\partial\Omega_{\m F})$. 

Now let $X_Q$ be the Corkscrew point (simultaneous for $\Omega$ and $\Omega_{\m F}$) for the cube $Q$. By the change of poles (Lemma 15.61 of \cite{dfm20}) and the doubling property of $\omega_\star$, for any Borel set $H_\star\subset\Delta_\star(y_Q,\hat r_Q)$ (see Notation \ref{not.centerofpj}), we have that
\[
\omega_\star^{X_Q}(H_\star)\approx\frac{\omega_\star^{X_0}(H_\star)}{\omega^{X_0}_\star(\Delta_\star(y_Q,\hat r_Q))}\approx\frac{\omega_\star^{X_0}(H_\star)}{\omega^{X_0}_\star(\Delta_\star^Q)},
\]
where in the last step we have used the fact that both $\Delta_\star(y_Q,\hat r_Q)$ and $\Delta_\star^Q$ have radius comparable to $\ell(Q)$, and that $\dist(\Delta_\star^Q,Q)\lesssim\ell(Q)$, so that we can compare both of these surface balls to a bigger (with radius still equivalent to $\ell(Q)$) surface ball containing both of them, and hence achieve the stated equivalence between these surface balls by using the doubling property of $\omega_\star$ (and, if needed, Harnack Chains and Harnack Inequality).

Using this last result, (\ref{eq.bigqball}), and (\ref{eq.djk2}), we see that
\begin{multline}\label{eq.djk3}
\frac{\m P_{\m F}\mu(F)}{\m P_{\m F}\mu(Q)}\lesssim\frac{\omega_\star^{X_0}(F\backslash(\cup_jQ_j))}{\omega_\star^{X_0}(\Delta_\star^Q)}+\sum_{Q_j\in\m F, Q_j\subseteq Q}\frac{\sigma(F\cap Q_j)}{\sigma(Q_j)}\frac{\omega_\star^{X_0}(P_j)}{\omega_\star^{X_0}(\Delta_\star^Q)}\\ \approx\omega_\star^{X_Q}(F\backslash(\cup_jQ_j))+\sum_{Q_j\in\m F, Q_j\subsetneq Q}\frac{\sigma(F\cap Q_j)}{\sigma(Q_j)}\omega_\star^{X_Q}(P_j).
\end{multline}
Next, we claim that the estimates
\begin{equation}\label{eq.compare}
\omega_\star^{X_Q}(F\backslash\cup_jQ_j)\lesssim\omega^{X_Q}(F\backslash\cup_jQ_j),\qquad\omega_\star^{X_Q}(P_j)\lesssim\omega^{X_Q}(Q_j)
\end{equation}
hold. Indeed, the first one follows immediately by the maximum principle (Lemma 12.8 in \cite{dfm20}) since $\Omega_{\m F}\subset\Omega$. For the second estimate, let $u(X):=\omega^X(Q_j)$ and $u_\star(X):=\omega_\star^X(P_j)$, and we first note that $\omega^X(Q_j)\approx1=u_\star(X)$ for each $X\in P_j$, by (\ref{eq.nondeg2}) and the fact that $\dist(P_j,Q_j)\approx\diam P_j\approx\diam Q_j\approx\ell(Q_j)$. Since also we have that $u_\star(X)=0\leq u(X)$ for each $X\in\partial\Omega_{\m F}\backslash P_j$, we may thus apply again the maximum principle to conclude that $u_\star(X_Q)\lesssim u(X_Q)$, as desired.

Finally, from (\ref{eq.djk3}) and (\ref{eq.compare}) we deduce that
\begin{multline*}
\frac{\m P_{\m F}\mu(F)}{\m P_{\m F}\mu(Q)}\lesssim\omega^{X_Q}(F\backslash(\cup_jQ_j))+\sum_{Q_j\in\m F, Q_j\subsetneq Q}\frac{\sigma(F\cap Q_j)}{\sigma(Q_j)}\omega^{X_Q}(Q_j)\\ \lesssim\frac{\omega^{X_0}(F\backslash(\cup_jQ_j))}{\omega^{X_0}(Q)}+\sum_{Q_j\in\m F, Q_j\subsetneq Q}\frac{\sigma(F\cap Q_j)}{\sigma(Q_j)}\frac{\omega^{X_0}(Q_j)}{\omega^{X_0}(Q)}\\=\frac{\m P_{\m F}\omega(F)}{\omega(Q)}=\frac{\m P_{\m F}\omega(F)}{\m P_{\m F}\omega(Q)},
\end{multline*}
where in the second line we used the change of poles for $\omega$ (\ref{eq.cphm2}), and in the last line we used that $\m P_{\m F}\omega(Q)=\omega(Q)$. This ends the proof of the right-hand side of (\ref{eq.projthetag}). The left-hand side is obtained because $\m P_{\m F}\omega$ is dyadically doubling by Lemma \ref{lm.pisdoubling}, $\m P_{\m F}\mu$ is dyadically doubling by Lemma 9.51 of \cite{dm} (whose argument is very similar to that of Lemma B.2 from \cite{hm2}; see our Remark \ref{rm.closetohm}), and $A_{\infty}^{dyadic}$ is an equivalence relationship among dyadically doubling measures (see Lemma \ref{lm.ainftrans}).\hfill{$\square$}

\section{Proof of Theorem \ref{thm.main}}\label{sec.proof}

In this section, we give the proof of Theorem \ref{thm.main}; we mainly follow the outline in \cite{hm1}; see also \cite{chm}.

Let $A_0$, $A$ be two matrices that satisfy (\ref{eq.elliptic}), write $\m A_0=w^{-1}A_0$, $\m A=w^{-1}A$, and suppose that $d\lambda(X)=\frac{\f a(X)^2}{\delta(X)^{n-d}}dX$ is a (continuous) Carleson measure, where $\f a$ is defined in (\ref{eq.disagreement}). As in Lemma \ref{lm.todyadic}, the natural discretization of the Carleson measure $\lambda$ is the collection $\f m=\{\alpha_Q\}_{Q\in\bb D}$ with
\[
\alpha_Q=\sum_{I\in\m W_Q}\frac{\sup_{Y\in I^*}|\f E(Y)|^2}{\ell(I)^{n-d}}|I|,\qquad Q\in\bb D,
\]
where $\f E(Y)=\m A(Y)-\m A_0(Y)$. Let $L_0=-\dv A_0\nabla$,  $L=-\dv A\nabla$, and let $\omega_0$, $\omega$ be the harmonic measures of $L_0$, $L$ respectively. 

Our program is to apply Theorem \ref{thm.extrapolation} to eventually obtain that $\omega_L\in A_{\infty}^{\dyadic}$. Of course, this will imply that $\omega_L\ll\sigma$ and that $\frac{d\omega_L}{d\sigma}\in RH_q^{\dyadic}$ for some $q>0$, and by Remark \ref{rm.equivrh}, the latter is equivalent to $\frac{d\omega_L}{d\sigma}\in RH_q$, which in turn means that $\omega_L\in A_{\infty}$.

Thus we ought to verify the hypotheses of Theorem \ref{thm.extrapolation}.  Fix $Q_0\in\bb D$, and observe that $\Vert\f m\Vert_{\m C(Q_0)}\lesssim\vertiii{\lambda}_{\m C}$ by Lemma \ref{lm.todyadic}. Given $\xi>0$ small enough and to be chosen later, we fix a disjoint family $\m F=\{Q_j\}_j\in\bb D_{Q_0}$ that verifies the estimate
\begin{equation}\label{eq.cond1}
\Vert\f m_{\m F}\Vert_{\m C(Q_0)}=\sup_{Q\in\bb D_{Q_0}}\frac{\f m(\bb D_{\m F,Q})}{\sigma(Q)}\leq\xi.
\end{equation}
Recall that $R_{Q_0}\subset B(x_{Q_0},7\sqrt nA_0\ell(Q_0))$ (see (\ref{eq.sawtoothbdd})). Let $X_0\in\Omega$ be a Corkscrew point (with Corkscrew constant $c$) for the surface ball $\Delta(x_{Q_0},10c^{-1}\sqrt n A_2\ell(Q_0))$. Then
\[
|X_0-x_{Q_0}|\geq\delta(X_0)\geq c10c^{-1}\sqrt n A_2\ell(Q_0)=10\sqrt nA_2\ell(Q_0),
\]
which implies that $X_0\in\Omega\backslash B(x_{Q_0},10\sqrt nA_2\ell(Q_0))\subset\Omega\backslash R_{Q_0}^{**}\subset\Omega\backslash R_{Q_0}$, where 
\[
R_{Q_0}^{**}:=\opint\Big(\bigcup_{I\in\m R_{Q_0}}I^{**}\Big),\qquad I^{**}=(1+2\theta)I.
\]
Moreover, according to Corollary \ref{cor.dd}, we have that $\omega^{X_0}$ is dyadically doubling in $Q_0$, while we also have that $\delta(X_0)\approx\dist(X_0,Q_0)\approx\ell(Q_0)$.

We want to show that $\m P_{\m F}\omega_L^{X_0}$ satisfies (\ref{eq.propextra}), with uniform constants and with $\omega_L^{X_0}$ in place of $\mu$. Since $L_0$ is solvable in some $L^{p'}$, then by Theorem \ref{thm.equivalence}, $\omega_{L_0}^{X_0}\ll\sigma$ and $k_0^{X_0}\in RH_p(\Delta(x_{Q_0},A_0\ell(Q_0))$ (with reverse H\"older characteristic independent of $Q_0$).

\subsection{Step 0: A qualitative reduction} We first make a reduction that allows us to conjure \emph{qualitative} absolute continuity properties of the harmonic measure $\omega_L$.

\begin{definition}[Tubes encasing the boundary]\label{def.tube} Fix $\tau>0$. By a \emph{$\tau-$tube around $\Gamma$}, we mean the open set $\Gamma_{\tau}:=\big\{X\in\Omega\,:\,\text{dist}(X,\Gamma)<\tau\big\}$.
\end{definition}

We define $A_{\tau}$ as $A_{\tau}=A_0$ in the $\tau-$tube $\Gamma_{\tau}$, and $A_{\tau}=A$ in $\bb R^n\backslash\Gamma_{\tau}$. In all of the following steps we work with $L_{\tau}=-\dv A_{\tau}\nabla$ in place of $L$. We note that the ellipticity constants of $A_{\tau}$ are controlled by those of $A$ and $A_0$, uniformly in $\Omega$. The same is true of the condition on the disagreement $\f a$.
 
Let us now exploit Lemma \ref{lm.poissonequiv} to deduce absolute continuity properties of $\omega_{L_{\tau}}$, with dependence on $\tau$.

\begin{corollary}[Comparability of harmonic measures in tubes]\label{lm.match} Retain the notation above.  Then $\omega_{\tau}^{X_0}\ll\sigma$, and if $\tau$ is small enough depending on $n$, $d$, $C_d$ only, then $k_\tau^{X_0}\in RH_p(\Delta(x_{Q_0},A_0\ell(Q_0)))$, with the $RH_p$ characteristic depending on $\tau$ and $\ell(Q_0)$.
\end{corollary} 

We emphasize that this is a qualitative result (the dependence on $\tau$ and $\ell(Q_0)$ is non-optimal) - see also related comments after the proof.

\noindent\emph{Proof.} Fix a surface ball $\Delta=\Delta(x_0,r)$ with $r\in(0,\frac{c\tau}4)$, and let $X_{\Delta}$ be a Corkscrew point for $\Delta$. By Lemma \ref{lm.poissonequiv}, we have that $\omega_0^{X_\Delta}$ is mutually absolutely continuous with $\omega_\tau^{X_\Delta}$ on $\Delta$. Recall that for each $i=0,\tau$,  $\omega_i^{X_{\Delta}}$ is mutually absolutely continuous with   $\omega_i^{X_0}$. It follows that $\omega_0^{X_\Delta}\ll\sigma$ on $\Gamma$, and therefore that $\omega_\tau^{X_0}\ll\omega_\tau^{X_\Delta}\ll\omega_0^{X_\Delta}\ll\sigma$ on $\Delta$. Since $\Delta\subset\Gamma$ was arbitrary, we have that $\omega_\tau^{X_0}\ll\sigma$ on $\Gamma$.

Next, fix a surface ball $\Delta_{\tau}:=\Delta(x,c\tau/4)$ with $x\in\Delta(x_{Q_0},10c^{-1}\sqrt n A_2\ell(Q_0))$, and let $X_{\Delta_{\tau}}$ be a Corkscrew point for $\Delta_{\tau}$. Then Lemma \ref{lm.poissonequiv} (iii) gives that $k_0^{X_{\Delta_{\tau}}}(y)\approx k_\tau^{X_{\Delta_{\tau}}}(y)$ for $\sigma-$a.e. $y\in\Delta_{\tau}$. Then the Harnack Chains and Harnack inequality guarantee the estimates
\[
k_\tau^{X_0}(y)\approx_{\tau}k_\tau^{X_{\Delta_{\tau}}}(y)\approx k_0^{X_{\Delta_{\tau}}}(y)\approx_{\tau}k_0^{X_0}(y),\qquad\text{for }\sigma-\text{a.e. }y\in\Delta_{\tau}.
\]
Since we may cover $\Delta(x_{Q_0},A_0\ell(Q_0))$ by surface balls $\{\Delta_{\tau}\}$ as above and with uniformly bounded overlap, it follows that 
\[
k_\tau^{X_0}(y)\approx_{\tau}k_0^{X_0}(y),\qquad\text{for }\sigma-\text{a.e. }y\in\Delta(x_{Q_0},A_0\ell(Q_0))\supset Q_0.
\]
The desired result ensues.\hfill{$\square$}

It follows that we may assume that all the harmonic measures $\omega_{\tau}=\omega_{L_{\tau}}$ are absolutely continuous with respect to $\sigma$, and $k_{\tau}^{X_0}=k_{L_{\tau}}^{X_0}\in RH_p(Q_0)$ with the $RH_p$ characteristic depending on $\tau$ and $\ell(Q_0)$. The dependence on $\tau$ and $\ell(Q_0)$ will not be an issue because these facts are used only qualitatively. 

Therefore, in Step 1 below, we will have \emph{a priori} that $\omega^{X_0}_{\tau}\ll\sigma$ and that $k^{X_0}_{\tau}\in L^p(Q_0,\sigma)$. We eventually establish a reverse H\"older inequality for $k_\tau$ with $RH$ exponent and characteristic independent of $\tau$ and $\ell(Q_0)$. We will finally pass to the limit using Lemma \ref{lm.opconv} below to conclude that $\omega_L\ll\sigma$ and $\frac{d\omega_L}{d\sigma}\in RH_q$. This will in turn imply as desired that $L$ is solvable in $L^{q'}$ by Theorem \ref{thm.equivalence}.

\subsection{Step 1: Exploit smallness of $\Vert\f m_{\m F}\Vert_{\m C(Q_0)}$} Introduce the operator $L_1$ defined as $L_1=L_{\tau}$ in $\Omega_{\m F,Q_0}$, and $L_1=L_0$ in $\Omega\backslash\Omega_{\m F,Q_0}$. We write $\omega_1$ for the harmonic measure associated to the operator $L_1$, and $g_1$ for the Green function associated to $L_1$. We have that $k_0^{X_0}\in RH_p(\Delta(x_{Q_0},A_0\ell(Q_0)))$, and in particular by Theorem \ref{thm.equivalence}, we have that
\begin{equation*} 
\int_{Q_0}\big(k_0^{X_0}\big)^p\,d\sigma\leq\int_{\Delta(x_{Q_0},A_0\ell(Q_0))}\big(k_0^{X_0}\big)^p\,d\sigma\lesssim C_0\sigma(\Delta(x_{Q_0},A_0\ell(Q_0)))^{1-p}\approx \sigma(Q_0)^{1-p};
\end{equation*}
thus, in summary,
\begin{equation}\label{eq.poissonrh0}
\int_{Q_0}\big(k_0^{X_0}\big)^p\,d\sigma\lesssim\sigma(Q_0)^{1-p}.
\end{equation}

Our immediate goal in Step 1 is to show that (\ref{eq.poissonrh0}) remains true when $k_0^{ X_0}$ is replaced by $k_1^{X_0}$, the Poisson kernel for the operator $L_1$ defined above.

Let $f\geq0$ be a continuous function supported on $Q_0$, such that $\Vert f\Vert_{L^{p'}(Q_0,\sigma)}=1$, and let $u_0$ and $u_1$ be the corresponding solutions to the Dirichlet problems for $L_0$ and $L_1$ with boundary data $f$. Set $\m E_1(Y)=A_1(Y)-A_0(Y)=\m E(Y){\bf 1}_{\Omega_{\m F,Q_0}}(Y)$, where $\m E(Y)=A_{\tau}(Y)-A_0(Y)$. Then, we may write
\begin{multline}\label{eq.calc1}
F_1(X_0):=|u_1(X_0)-u_0(X_0)|=\Big|\dint_{\bb R^n}\nabla_Yg_1^T(Y,X_0)\m E_1^T(Y)\nabla u_0(Y)\,dY\Big|\\ \leq\dint_{\Omega_{\m F,Q_0}}|\nabla_Yg_1(X_0,Y)||\m E(Y)||\nabla u_0(Y)|\,dY \\ \leq\sum\limits_{Q\in\bb D_{\m F,Q_0}}\sum_{I\in\m W_Q}\dint_{I^*}|\nabla_Yg_1(X_0,Y)||\m E(Y)||\nabla u_0(Y)|\,dY \\ \leq\sum\limits_{Q\in\bb D_{\m F,Q_0}}\sum_{I\in\m W_Q}\big(\sup_{Y\in I^*}|\m E(Y)|\big)\Big(\dint_{I^*}|\nabla_Yg_1(X_0,Y)|^2\,dY\Big)^{\frac12}\Big(\dint_{I^*}|\nabla u_0(Y)|^2\,dY\Big)^{\frac12},
\end{multline}
where we have used (\ref{eq.diffid}) in the first line, and later H\"older's inequality. By definition of $X_0$, we have that $v(Y)=g_1(X_0,Y)=g_1^T(Y,X_0)$ is a non-negative solution of $L_1^Tv=0$ in $R_{Q_0}^{**}$ (as $X_0\notin R_{Q_0}^{**}$). Hence, we can apply Caccioppoli's inequality (see Lemma 8.6 of \cite{dfm1}) to obtain that
\begin{multline}\label{eq.cacc1}
\dint_{I^*}|\nabla_Yg_1(X_0,Y)|^2\,dY\lesssim\frac1{\ell(I)^{-n+d+1}}\dint_{I^*}|\nabla_Yg_1(X_0,Y)|^2w(Y)\,dY\\ \lesssim\ell(I)^{-2}\frac1{\ell(I)^{-n+d+1}}\dint_{I^{**}}|g_1(X_0,Y)|^2\,dm\approx \dint_{I^{**}}\frac{|g_1(X_0,Y)|^2}{\delta(Y)^2}\,dY,
\end{multline}
for any $I\in\m W_Q$, $Q\in\bb D_{\m F,Q_0}$. Fix such an $I$ and $Q$. We have by the Harnack inequality and Harnack chains   that $g_1(X_0,Y)\approx g_1(X_0,X_Q)$ for all $Y\in I^{**}$ and where $X_Q$ is a Corkscrew point for $Q$. Then by  Lemma \ref{lm.greenlow}, Lemma \ref{lm.greenhigh}, and Corollary \ref{cor.likesurfaceball}, for every $Y\in I^{**}$ we have that
\begin{equation}\label{eq.usecfms}
\frac{g_1(X_0,Y)}{\delta(Y)}\approx\frac{g_1(X_0,X_Q)}{\delta(X_Q)}\approx\frac{\omega_1^{X_0}\big(\Delta(x_Q,a_0\ell(Q))\big)}{\ell(Q)\ell(Q)^{d-1}}\approx\frac{\omega_1^{X_0}(Q)}{\sigma(Q)}.
\end{equation}
Putting together (\ref{eq.cacc1}) and (\ref{eq.usecfms}),   we see that
\begin{equation}\label{eq.calc2}
\dint_{I^*}|\nabla_Yg_1(X_0,Y)|^2\,dY\lesssim\Big(\frac{\omega_1^{X_0}(Q)}{\sigma(Q)}\Big)^2|I|.
\end{equation}
Plugging (\ref{eq.calc2}) into (\ref{eq.calc1}), using the fact that
\[
\sup_{Y\in I^*}|\m E(Y)|\approx\ell(I)^{1-(n-d)}\sup_{Y\in I^*}|\f E(Y)|,
\]
and using the Cauchy-Schwartz inequality, we obtain that
\begin{multline}\label{eq.calc3}
F_1(X_0)\lesssim\sum_{Q\in\bb D_{\m F,Q_0}}\sum_{I\in\m W_Q}\ell(I)^{1-\frac{n-d}2}\frac{\omega_1^{X_0}(Q)}{\sigma(Q)}\Big(\frac{\sup_{Y\in I^*}|\f E(Y)|^2}{\ell(I)^{n-d}}|I|\Big)^{\frac12}\Big(\dint_{I^*}|\nabla u_0|^2\Big)^{\frac12}\\ \lesssim\sum_{Q\in\bb D_{\m F,Q_0}}\Big(\frac{\alpha_Q}{\sigma(Q)}\Big)^{\frac12}\sigma(Q)\frac{\omega_1^{X_0}(Q)}{\sigma(Q)}\Big(\dint_{U_Q}|\nabla u_0(X)|^2\delta(X)^{2-n}\,dX\Big)^{\frac12}\\ \lesssim\sum_{Q\in\bb D_{\m F,Q_0}}\Big(\frac{\f m(\bb D_{\m F,Q})}{\sigma(Q)}\Big)^{\frac12}\int_Q\big(M(k_1^{X_0}{\bf 1}_{Q_0})\big)(x)\Big(\dint_{\gamma_d^Q(x)}|\nabla u_0(X)|^2\delta(X)^{2-n}\,dX\Big)^{\frac12}\,dx\\ \lesssim\Vert\f m_{\m F}\Vert_{\m C(Q_0)}^{\frac12}\sum_{Q\in\bb D_{\m F,Q_0}}\int_Q\big(M(k_1^{X_0}{\bf 1}_{Q_0})\big)(x)\Big(\dint_{\gamma^{\alpha_1}(x)}|\nabla u_0(X)|^2\delta(X)^{2-n}\,dX\Big)^{\frac12}\,dx\\ \lesssim\Vert\f m_{\m F}\Vert_{\m C(Q_0)}^{\frac12} \int_\Gamma\big(M(k_1^{X_0}{\bf 1}_{Q_0})\big)(x)\big(S_{\alpha_1}u_0\big)(x)\,dx
\end{multline}
where in the second line we used that $\ell(I)\approx\ell(Q)$ for each $I\in\m W_Q$, the definition of $\alpha_Q$, and the bounded overlap of the dylated Whitney cubes $I^*$; in the third line we used for each $x\in Q$ the estimate  
\begin{multline}\label{eq.justifyhl}\notag
\frac{\omega^{X_0}_1(Q)}{\sigma(Q)}=\frac1{\sigma(Q)}\int_Qk_1^{X_0}\,d\sigma\lesssim\frac1{\Delta(x,A_0\ell(Q))}\int_{\Delta(x, A_0\ell(Q))}k_1^{X_0}{\bf 1}_{Q_0}\,d\sigma\\ \leq\big(M(k_1^{X_0}{\bf 1}_{Q_0})\big)(x),\qquad\text{for each }x\in Q,
\end{multline} 
and in the fourth line we chose $\alpha_1$ so that $\gamma_d(x)\subset\gamma^{\alpha_1}(x)$ for all $x\in\Gamma$ (see the paragraph following (\ref{eq.localdyadicntcone})). Using (\ref{eq.cond1}),(\ref{eq.calc3}) and the H\"older's inequality, we furnish the estimate
\begin{multline*}
F_1(X_0)\lesssim\xi^{\frac12}\Vert M(k_1^{X_0}{\bf 1}_{Q_0})\Vert_{L^p(\Gamma,\sigma)}\Vert S_{\alpha_1}u_0\Vert_{L^{p'}(\Gamma,\sigma)}\lesssim\xi^{\frac12}\Vert k_1^{X_0}\Vert_{L^p(Q_0,\sigma)}\Vert f\Vert_{L^{p'}(Q_0,\sigma)}\\ =\xi^{\frac12}\Vert k_1^{X_0}\Vert_{L^p(Q_0,\sigma)},
\end{multline*}
where we have used Theorem \ref{thm.djk} and the fact that the trace of $u_0$ is $f$. As such, we have that
\begin{equation}\label{eq.calc4}
|u_1(X_0)-u_0(X_0)|=F_1(X_0)\lesssim\xi^{\frac12}\Vert k_1^{X_0}\Vert_{L^p(Q_0,\sigma)}.
\end{equation}
We will now see that (\ref{eq.calc4}) implies the desired result. By the definitions of $u_0$, $u_1$, $f$, and H\"older's inequality, (\ref{eq.calc4}) gives that
\begin{multline*}
\int_{Q_0}fk_1^{X_0}\,d\sigma\lesssim\xi^{\frac12}\Vert k_1^{X_0}\Vert_{L^p(Q_0,\sigma)}+\int_{Q_0}fk_0^{X_0}\,d\sigma\\ \leq\xi^{\frac12}\Vert k_1^{X_0}\Vert_{L^p(Q_0,\sigma)}+\Vert f\Vert_{L^{p'}(\Gamma,\sigma)}\Vert k_0^{X_0}\Vert_{L^p(Q_0,\sigma)} \leq\xi^{\frac12}\Vert k_1^{X_0}\Vert_{L^p(Q_0,\sigma)}+\Vert k_0^{X_0}\Vert_{L^p(Q_0,\sigma)}.
\end{multline*}
Since the continuous functions are dense in $L^{p'}(\Gamma,\sigma)$, by taking supremum over all possible $f$ as described earlier, we obtain that
\begin{equation}\label{eq.hide}
\Vert k_1^{X_0}\Vert_{L^p(Q_0,\sigma)}\lesssim\xi^{\frac12}\Vert k_1^{X_0}\Vert_{L^p(Q_0,\sigma)}+\Vert k_0^{X_0}\Vert_{L^p(Q_0,\sigma)}.
\end{equation}
It follows that as long as $\xi^{\frac12}$ is small enough (depending only on the permissible constants), we may hide the first term on the right-hand side of the above inequality to the left-hand side; hence we get that
\[
\Vert k_1^{X_0}\Vert_{L^p(Q_0,\sigma)}\lesssim\Vert k_0^{X_0}\Vert_{L^p(Q_0,\sigma)}.
\]
And hence, since $k_0^{X_0}$ satisfies (\ref{eq.poissonrh0}), we obtain that $k_1^{X_0}$ satisfies (\ref{eq.poissonrh0}) as well, with the implicit constant independent of $\tau$ and $Q_0$.

\subsection{Self-improvement of Step 1.} We currently have (\ref{eq.poissonrh0}) only for the largest cube $Q_0$; but let us see that we can extend (\ref{eq.poissonrh0}) to obtain a reverse H\"older estimate on every dyadic subcube of $Q_0$.

Fix $Q\in\bb D_{Q_0}$. Let $X^Q$ be a Corkscrew point for the surface ball $\Delta(x_Q,10c^{-1}\sqrt nA_2\ell(Q))$, so that $X^Q\in\Omega\backslash R_Q^{**}$ (see the remarks about $X_0$ following (\ref{eq.cond1})). Define a new operator $L_1^Q=L_{\tau}$ in $\Omega_{\m F,Q}$ and $L_1^Q=L_0$ otherwise in $\Omega\backslash\Omega_{\m F,Q}$, and let $k_{L_1^Q}^{X^Q}$ denote the Poisson kernel for $L_1^Q$ with pole at $X^Q$. Given our proof above, it is easy to see that
\begin{equation}\label{eq.poissonrhq}
\int_Q(k_{L_1^Q}^{X^Q})^p\,d\sigma\leq C_1\sigma(Q)^{1-p},
\end{equation}
for some $C_1$ independent of $Q$ (and $\tau$). Indeed, if $Q\subset Q_k$ for some $Q_k\in\m F$ then we obtain that $\bb D_{\m F,Q}=\varnothing$ so that $\Omega_{\m F,Q}=\varnothing$ and so $L_1^Q\equiv L_0$ in $\Omega\backslash\Gamma$. In that case, (\ref{eq.poissonrhq}) holds by hypothesis. Otherwise,   trivially we have that $\Vert\f m_{\m F}\Vert_{\m C(Q)}\leq\Vert\f m_{\m F}\Vert_{\m C(Q_0)}\leq\xi$, and consequently, if $Q$ is not contained in any $Q_k\in\m F$, then we may simply repeat the previous argument with respect to $Q$, and we obtain (\ref{eq.poissonrhq}) exactly as before. This proves the claim.

Now, by the non-degeneracy of the harmonic measure we have that $\int_Qk_{L_1^Q}^{X^Q}\,d\sigma\gtrsim1$, and combining this estimate with (\ref{eq.poissonrhq}) we obtain that
\begin{equation}\label{eq.poissonrh1}
\Big(\frac1{\sigma(Q)}\int_Q(k_{L_1^Q}^{X^Q})^p\,d\sigma\Big)^{\frac1p}\lesssim\frac1{\sigma(Q)}\int_Qk_{L_1^Q}^{X^Q}\,d\sigma.
\end{equation}
Next, we want to pass from $k_{L_1^Q}^{X^Q}$ to $k_{L_1}^{X^Q}$. Notice that $L_1\equiv L_1^Q$ in $(\Omega\backslash\Omega_{\m F,Q_0})\cup\Omega_{\m F,Q}$, and observe that for $B_s=B(x_Q,s)$ the ball in Lemma \ref{lm.carlesonball} ii) for which $\ell(Q)\lesssim s\leq\ell(Q)$ and (\ref{eq.carlesonball2}) is satisfied, we have that
\begin{multline*}
B_s\cap\Omega=B_s\cap\big[(\Omega\backslash\Omega_{\m F, Q_0})\cup\Omega_{\m F,Q_0}\big]=\big(B_s\cap(\Omega\backslash\Omega_{\m F, Q_0})\big)\cup\big(B_s\cap\Omega_{\m F,Q_0}\big)\\=\big(B_s\cap(\Omega\backslash\Omega_{\m F, Q_0})\big)\cup\big(B_s\cap\Omega_{\m F,Q}\big)=B_s\cap\big(\Omega_{\m F,Q}\cup\Omega\backslash\Omega_{\m F,Q_0}\big),
\end{multline*}
and hence $L_1\equiv L_1^Q$ in $B_s\cap\Omega$. Therefore, using Lemma \ref{lm.poissonequiv} we see that there exists $\tilde\ep$, $1\lesssim\tilde\ep\leq1$ so that
\[
k_1^{X^Q}(y)=k_{L_1}^{X^Q}(y)\approx k_{L_1^Q}^{X^Q}(y),\qquad\text{for }\sigma-\text{almost every } y\in\tilde\ep\Delta_Q=\Delta(x_Q,\tilde\ep\ell(Q)).
\]
We use this result, (\ref{eq.poissonrh1}), and the doubling property Lemma \ref{lm.doublinglemma} to deduce the estimate
\begin{multline}\label{eq.calc1step1ref}\notag
\Big(\frac1{\sigma(\tilde\ep\Delta_Q)}\int_{\tilde\ep\Delta_Q}(k_1^{X^Q})^p\,d\sigma\Big)^{\frac1p}\leq\Big(\frac1{\sigma(\tilde\ep\Delta_Q)}\int_Q(k_{L_1^Q}^{X^Q})^p\,d\sigma\Big)^{\frac1p}\\  \lesssim\frac1{\sigma(\tilde\ep\Delta_Q)}\int_Qk_{L_1^Q}^{X^Q}\,d\sigma\leq\frac1{\sigma(\tilde\ep\Delta_Q}\int_{\Delta(x_Q,A_0\ell(Q))}k_{L_1^Q}^{X^Q}\,d\sigma\\  \lesssim\frac1{\sigma(\tilde\ep\Delta_Q)}\int_{\tilde\ep\Delta_Q}k_{L_1^Q}^{X^Q}\,d\sigma \approx\frac1{\sigma(\tilde\ep\Delta_Q)}\int_{\tilde\ep\Delta_Q}k_1^{X^Q}\,d\sigma .
\end{multline}
We can now change poles from $X^Q$ to $X_0$ via Lemma \ref{lm.poissonpole} to obtain
\begin{equation}\label{eq.calc2step1ref}\notag
\Big(\frac1{\sigma(\tilde\ep\Delta_Q)}\int_{\tilde\ep\Delta_Q}(k_1^{X_0})^p\,d\sigma\Big)^{\frac1p}\lesssim\frac1{\sigma(\tilde\ep\Delta_Q)}\int_{\tilde\ep\Delta_Q}k_1^{X_0}\,d\sigma,\qquad\text{for each } Q\in\bb D_{Q_0}.
\end{equation}
Finally, we use Lemma \ref{lm.technical} to furnish

\begin{conclusion}[Step 1] We have that $\omega_1^{X_0}\in A_{\infty}^{\dyadic}(Q_0)$ uniformly in $\tau$ and $Q_0$. Hence we deduce that $\m P_{\m F}\omega_1^{X_0}\in A_{\infty}^{\dyadic}(Q_0)$ (uniformly in $\tau$ and $Q_0$) by Lemma \ref{lm.projdyadic}.
\end{conclusion}

\subsection{Step 2: Hide the ``bad'' Carleson regions.} We define the operator $L_2$ such that the disagreement with $L_1$ lives roughly inside the Carleson regions corresponding to the family $\m F$. More precisely, set
\[
L_2=\left\{\begin{matrix}L_{\tau},&\text{ in }&R_{Q_0}\backslash\Omega_{\m F},\\L_1,&\text{ in }&\Omega\backslash(R_{Q_0}\backslash\Omega_{\m F}).\end{matrix}\right.
\]
Note carefully that $R_{Q_0}\backslash\Omega_{\m F}\subseteq R_{Q_0}\backslash\Omega_{\m F,Q_0}$, but the opposite containment does not hold in general.    We write $\omega_1=\omega_{L_1}^{X_0}$ and $\omega_2=\omega_{L_2}^{X_0}$ for the corresponding harmonic measures (on $\Gamma$) for $L_1$ and $L_2$ with fixed pole at $X_0$. We also let $\omega_{\star,1}=\omega_{\star,1}^{X_0}$ and $\omega_{\star,2}=\omega_{\star,2}^{X_0}$ denote the harmonic measures of $L_1$ and $L_2$ on $\partial\Omega_{\m F}$, the boundary of the dyadically-generated sawtooth domain $\Omega_{\m F}$ (see the beginning of Section \ref{sec.djk}). Note that
\[
\Omega_{\m F}\subset\Omega\backslash(R_{Q_0}\backslash\Omega_{\m F}),
\]
whence $L_1\equiv L_2$ in $\Omega_{\m F}$ and consequently we have that  $\omega_{\star,1}\equiv\omega_{\star,2}$.

Next, we apply Lemma \ref{lm.djkg} to both $L_1$ and $L_2$ to deduce that for all $Q\in\bb D_{Q_0}$ and $F\subset Q$, the estimate
\[
\Big(\frac{\m P_{\m F}\omega_i(F)}{\m P_{\m F}\omega_i(Q)}\Big)^{\theta_i}\lesssim\frac{\m P_{\m F}\mu_i(F)}{\m P_{\m F}\mu_i(Q)}\lesssim\frac{\m P_{\m F}\omega_i(F)}{\m P_{\m F}\omega_i(Q)} 
\]
holds for $i=1,2$, and $\mu_i$ is defined in (\ref{eq.mug}). Observe that $\m P_{\m F}\mu_1\equiv\m P_{\m F}\mu_2$ since $\omega_{\star,1}\equiv\omega_{\star,2}$. Since $A_{\infty}^{dyadic}(Q_0)$ defines an equivalence relationship among dyadically doubling measures (which the projection measures $\m P_{\m F}\mu_i$ are; see Lemma 9.51 of \cite{dm} or Lemma B.2 of \cite{hm2}), and since we showed in Step 1 that $\m P_{\m F}\omega_1\in A_{\infty}^{dyadic}(Q_0)$, we obtain in this step that $\m P_{\m F}\omega_2\in A_{\infty}^{dyadic}(Q_0)$. For definiteness, we have
 
\begin{conclusion}[Step 2] There exist $\theta,\theta'>0$ (independent of $\tau$ and $Q_0$) such that for all $Q\in\bb D_{Q_0}$ and all Borel sets $F\subseteq Q$, we have the estimate
\[
\frac1C\Big(\frac{\sigma(F)}{\sigma(Q)}\Big)^{\theta}\leq\frac{\m P_{\m F}\omega_2^{X_0}(F)}{\m P_{\m F}\omega_2^{X_0}(Q)}\leq C\Big(\frac{\sigma(F)}{\sigma(Q)}\Big)^{\theta'}
\]
with $C$ uniform in $\tau$ and $Q_0$.
\end{conclusion}

\subsection{Step 3: Extend outside the Carleson region of $Q_0$} Observe that
\[
\Omega=\Omega_{\m F,Q_0}\cup(R_{Q_0}\backslash\Omega_{\m F})\cup(R_{Q_0}\cap\Omega_{\m F}\backslash\Omega_{\m F,Q_0})\cup(\Omega\backslash R_{Q_0}).
\]
We have successfully changed the operator from $L_0$ to $L_{\tau}$ on $\Omega_{\m F,Q_0}\cup(R_{Q_0}\backslash\Omega_{\m F})$ in the last two steps; it remains to change it in the set $\hat\Omega:=(R_{Q_0}\cap\Omega_{\m F}\backslash\Omega_{\m F,Q_0})\cup(\Omega\backslash R_{Q_0})$. Thus we define
\[
L_3=\left\{\begin{matrix}L_\tau,&\text{ in }&\hat\Omega,\\L_2,&\text{ in }&\Omega\backslash\hat\Omega.\end{matrix}\right.
\]
Hence, note that $L_3=L_\tau$ in $\Omega$, and $L_3\equiv L_2$ in $\Omega\backslash\hat\Omega=R_{Q_0}\backslash(\Omega_{\m F}\backslash\Omega_{\m F,Q_0})$. We will show that (\ref{eq.propextra}) holds, so fix $\ep\in(0,1)$ and take $E\subset Q_0$ with $\frac{\sigma(E)}{\sigma(Q_0)}\geq\ep$. In the case that $\m F=\{Q_0\}$, we have the trivial estimate
\[
\frac{\m P_{\m F}\omega_3^{X_0}(E)}{\m P_{\m F}\omega_3^{X_0}(Q_0)}=\frac{\frac{\sigma(E)}{\sigma(Q_0)}\omega_3^{X_0}(Q_0)}{\frac{\sigma(Q_0)}{\sigma(Q_0)}\omega_3^{X_0}(Q_0)}=\frac{\sigma(E)}{\sigma(Q_0)}\geq\ep.
\]
We thus suppose that $\m F\subseteq\bb D_{Q_0}\backslash\{Q_0\}$. For $t\ll1$, recall that we define $\Sigma_{t}=\Sigma_{Q_0,t}$ in Lemma \ref{lm.poissoncube}. Define $Q_{t}=Q_0\backslash\medcup_{Q'\in\m I_{t}}Q'$, where
\[
\m I_{t}=\{Q'\in\bb D_{Q_0}:t\ell(Q_0)<\ell(Q')\leq2t\ell(Q_0), Q'\cap\Sigma_{t}\neq\varnothing\}.
\]
It is easy to see that $\Sigma_{t}\subset\medcup_{Q'\in\m I_{t}}Q'\subset\Sigma_{Ct}$ for $C$ a uniform constant. Then, for all $t=t(\ep)$ small enough, we have that  
\[
\sigma(Q_0\backslash Q_{t})\leq\sigma(\Sigma_{Ct})\lesssim A_0(Ct)^{\zeta}\sigma(Q_0)\leq\frac{\ep}2\sigma(Q_0),
\]
where we have used Lemma \ref{lm.dyadiccubes} \ref{item.bdrythin}. Letting $F=E\cap\ Q_{t}$, it follows that
\[
\ep\sigma(Q_0)\leq\sigma(E)\leq\sigma(F)+\frac{\ep}2\sigma(Q_0),
\]
and therefore $\sigma(F)/\sigma(Q_0)\geq\ep/2$. Using the conclusion of Step 2, we see that
\[
\frac{\m P_{\m F}\omega_2^{X_0}(F)}{\m P_{\m F}\omega_2^{X_0}(Q_0)}\gtrsim \Big(\frac{\sigma(F)}{\sigma(Q_0)}\Big)^{\theta}\geq\Big(\frac{\ep}2\Big)^{\theta}.
\]
Now we claim that $\m P_{\m F}\omega_3^{X_0}(F)\geq c_{\ep}\m P_{\m F}\omega_2^{X_0}(F)$. The point of our argument is that the region  of discrepancy between $A_2$ and $A_3$ is uniformly far away (depending on $t=t(\ep)$) from most of $Q_0$, allowing us to compare the Poisson kernels of $\omega_2$ and $\omega_3$ in the set $F\subseteq E$, which has been chosen so that it retains most of $\sigma(E)$ while staying far from the region of discrepancy between $A_2$ and $A_3$.   Since $L_2\equiv L_3$ in $R_{Q_0}\backslash(\Omega_{\m F}\backslash\Omega_{\m F,Q_0})$, then we have by Lemma \ref{lm.poissoncube} that 
\[
k_2^{X_0}(y)\approx_{t} k_3^{X_0}(y)\quad\text{for }\sigma-\text{almost every }y\in Q_{t}\subset Q_0\backslash\Sigma_{t},
\]
where the implicit constants depends on $t$ and hence on $\ep$. It is then the case that $\omega_2^{X_0}(F\backslash(\cup_{Q_j\in\m F}Q_j)\approx\omega_3^{X_0}(F\backslash(\cup_{Q_j\in\m F}Q_j))$, and thus we observe the estimate
\begin{multline}\label{eq.calc5}
\m P_{\m F}\omega_3^{X_0}(F)=\omega_3^{X_0}(F\backslash(\cup_{Q_j\in\m F}Q_j))+\sum_{Q_j\in\m F}\frac{\sigma(F\cap Q_j)}{\sigma(Q_j)}\omega_3^{X_0}(Q_j)\\ \geq c_{\ep}\omega_2^{X_0}(F\backslash(\cup_{Q_j\in\m F}Q_j))+\sum_{Q_j\in\m F}\frac{\sigma(F\cap Q_j)}{\sigma(Q_j)}\omega_3^{X_0}(Q_j).
\end{multline}
It remains to estimate the last term. We need only consider the cubes in $\m F$ that meet $F$. Let $Q_j\in\m F$ be such a cube.   If $Q_j\subset Q_t$, then again by Lemma \ref{lm.poissoncube} we have that $\omega_3^{X_0}(Q_j)\geq c_{\ep}\omega_2^{X_0}(Q_j)$. Otherwise, $Q_j\cap(Q_0\backslash Q_{t})\neq\varnothing$, whence there exists $Q'\in\m I_{t}$ such that $Q'\subsetneq Q_j$ (since $Q_j\cap Q_t\neq\varnothing$). Accordingly, $\ell(Q_j)>t\ell(Q_0)$. Now let $\tilde Q\in\bb D_{Q_j}$ be a dyadic descendant of $Q_j$ which contains $x_{Q_j}$ and verifying $\ell(\tilde Q)=2^{-M}\ell(Q_j)$ with $M=2(1+\log_2(A_0a_0^{-1}))\approx1$, so that $\ell(\tilde Q)\approx\ell(Q_j)$. Let us see that $\tilde Q\subset Q_0\backslash\Sigma_{a_0t/2}$. Indeed, choose $x^*\in\Gamma\backslash Q_j$ and $y^*\in\tilde Q$ so that $\dist(\Gamma\backslash Q_j,\tilde Q)=|x^*-y^*|$, and reckon that
\begin{multline*}
\dist(\Gamma\backslash Q_0,\tilde Q)\geq\dist(\Gamma\backslash Q_j,\tilde Q)=|x^*-y^*|\geq|x^*-x_{Q_j}|-|x_{Q_j}-y^*|\\ \geq a_0\ell(Q_j)-\diam\tilde Q\geq[a_0-A_02^{-M}]\ell(Q_j)>\tfrac{a_0}2t\ell(Q_0),
\end{multline*}
which does give our claim. We may then apply Lemma \ref{lm.poissoncube} one last time to see that $k_2^{X_0}(y)\approx_{t}k_3^{X_0}(y)$ for $\sigma-$almost every $y\in\tilde Q$, and therefore
\[
\omega_3^{X_0}(Q_j)\geq\omega_3^{X_0}(\tilde Q)\approx_{\ep}\omega_2^{X_0}(\tilde Q)\gtrsim\omega_2^{X_0}(Q_j),
\]
where we have used the doubling property of the harmonic measure on the dyadic cubes. We now plug this result back into (\ref{eq.calc5}) to obtain that
\begin{multline}\notag
\m P_{\m F}\omega_3^{X_0}(F)\geq c_{\ep}\omega_2^{X_0}(F\backslash(\cup_{Q_j\in\m F}Q_j))+c_{\ep}\sum_{Q_j\in\m F}\frac{\sigma(F\cap Q_j)}{\sigma(Q_j)}\omega_2^{X_0}(Q_j)=c_{\ep}\m P_{\m F}\omega_2^{X_0}(F).
\end{multline}

We have arrived at

\begin{conclusion}[Step 3] There exists $\xi>0$ for which the following statement holds: given $\ep\in(0,1)$, there is $C_{\ep}<\infty$ such that for every $Q_0\in\bb D$, if $\m F=\{Q_j\}_j\subset\bb D_{Q_0}$ is a disjoint family satisfying $\Vert\f m_{\m F}\Vert_{\m C(Q_0)}<\xi$, then for any Borel set $F\subset Q_0$,
\[
\frac{\sigma(F)}{\sigma(Q_0)}\geq\ep~\implies~\frac{\m P_{\m F}\omega_{L_\tau}^{X_0}(F)}{\m P_{\m F}\omega_{L_\tau}^{X_0}(Q_0)}\geq\frac1{C_{\ep}},
\]
and $C_{\ep}$ is uniform in $\tau$.
\end{conclusion}

\subsection{Step 4: Fix the pole} The conclusion of Step 3 above almost looks like what is needed; but we ought to improve it so that its conclusion holds for any cube $Q\in\bb D_{Q_0}$ while keeping the pole $X_0$ fixed. Nevertheless, this is not  difficult; the following result is immediate from the method proof in \cite{hm1}; we omit the details of the proof carried out in \cite{chm} in a very similar setting (see our Remark \ref{rm.closetohm}).

\begin{proposition}[Proposition 4.25 of \cite{chm}]\label{prop.step4prop} There exists $\xi>0$ for which the following statement holds: given $\ep\in(0,1)$, there is $C_{\ep}<\infty$ such that for every $Q_0\in\bb D_{Q_0}$ and for all $Q\in\bb D_{Q_0}$, if $\m F=\{Q_j\}_j\subset\bb D_{Q}$ is a disjoint family satisfying $\Vert\f m_{\m F}\Vert_{\m C(Q)}<\xi$, then for any Borel set $F\subset Q$, we have that
\[
\frac{\sigma(F)}{\sigma(Q)}\geq\ep~\implies~\frac{\m P_{\m F}\omega_{L_\tau}^{X_0}(F)}{\m P_{\m F}\omega_{L_\tau}^{X_0}(Q)}\geq\frac1{C_{\ep}},
\]
with $C_{\ep}$ uniform in $\tau$. Consequently, $\omega^{X_0}_{\tau}\in A_{\infty}^{\dyadic}(Q_0)$ uniformly in $\tau$ and $Q_0$. In particular, there exists $1<q<\infty$ such that $k_\tau^{X_0}\in RH_q^{\dyadic}(Q_0)$ uniformly in $Q_0\in\bb D$ and $\tau>0$; and therefore $k_\tau\in RH_q$ with $RH_q$ characteristic independent of $\tau$. 
\end{proposition}

\subsection{Step 5: Pass to the limit in $\tau$} Proposition \ref{prop.step4prop} is the desired conclusion for each operator $L_{\tau}$, $\tau>0$, with $RH_q$ characteristic independent of $\tau$. We now ought to pass to the limit as $\tau\ra0$ and argue that these quantitative absolute continuity properties are preserved. The required technology is the following result.

\begin{lemma}[Limiting lemma]\label{lm.opconv} Let $A_0$, $A$ be two matrices satisfying (\ref{eq.elliptic}), and write $L_0=-\dv A_0\nabla$, $L=-\dv A\nabla$. For each $\tau\geq0$ small enough, define $A_{\tau}=A_0$ in $\Gamma_{\tau}$ (the $\tau-$tube around $\Gamma$, see Definition \ref{def.tube}) and $A_{\tau}=A$ otherwise in $\bb R^n$. Accordingly, define the operator $L_{\tau}:=-\dv A_{\tau}\nabla$. Let $\{\omega_0^X\},\{\omega^X\}, \{\omega_{\tau}^X\}$ be the families of harmonic measures associated to the operators $L_0, L, L_{\tau}$ respectively. Assume that there exists $q\in(1,\infty)$ such that $\frac{d\omega_{L_{\tau}}}{d\sigma}\in RH_q$ with the $RH_q$ characteristic  uniformly bounded in $\tau$. Then $\omega_L\ll\sigma$ and $\frac{d\omega_L}{d\sigma}\in RH_q$.
\end{lemma}

\noindent\emph{Proof.} Fix the surface ball $\Delta_0\subset\Gamma$,   let $X_0$ be a Corkscrew point for the ball $\Delta_0$, and suppose that $\tau<\delta(X_0)/4$. We first show that $\omega_{\tau}^{X_0}\rightharpoonup\omega^{X_0}$ on $\Delta_0$ as $\tau\searrow0$. Define the functionals $\Phi$ and $\Phi_{\tau}$ on $C_c(\Delta_0)$ by
\[
\Phi(f)=\int_{\Delta_0}f\,d\omega^{X_0},\qquad \Phi_{\tau}(f)=\int_{\Delta_0}f\,d\omega_{\tau}^{X_0},\qquad f\in C_c(\Gamma).
\]
Let $u$ (respectively, $u_{\tau}$)  be the  unique solution  to the Dirichlet problem  $Lu=0$ in $\Omega$, $u|_{\Gamma}=f$ (respectively,  $Lu_{\tau}=0$ in $\Omega$, $u_{\tau}|_{\Gamma}=f$). By (\ref{eq.solnharmonicmeas}), we have that $\Phi(f)=u(X_0)$, $\Phi_{\tau}(f)=u_{\tau}(X_0)$. In this setting, using an elementary approximation argument, the Cauchy-Schwartz inequality, and Lemma \ref{lm.greenfn} iv), it is easy to show that the identity (\ref{eq.diffid})  holds for \emph{all}  $X\in\Omega\backslash{\Gamma_{2\tau}}$, since the pole  $X_0$  lies far away from the support of $A_{\tau}-A$.  Thus, according to Lemma \ref{lm.diffid}, we may write
\begin{multline}\label{eq.diffid2}\notag
|\Phi(f)-\Phi_{\tau}(f)|=|u(X_0)-u_{\tau}(X_0)|\leq\Big|\dint_{\Gamma_{\tau}}(A-A_{\tau})^T(Y)\nabla_Yg^T_{\tau}(Y,X_0)\nabla u(Y)\,dY\Big|\\ \leq2(C_A+C_{A_0})\Big(\dint_{\Gamma_{\tau}}|\nabla_Yg^T_{\tau}(Y,X_0)|^2\,dm(Y)\Big)^{\frac12}\Big(\dint_{\Gamma_{\tau}}|\nabla u(Y)|^2\,dm(Y)\Big)^{\frac12}\\ \leq2(C_A+C_{A_0})\Big(\dint_{\Omega\backslash B(X_0,\delta(X_0)/2)}|\nabla_Yg^T_{\tau}(Y,X_0)|^2\,dm(Y)\Big)^{\frac12}\Big(\dint_{\Gamma_{\tau}}|\nabla u(Y)|^2\,dm(Y)\Big)^{\frac12}\\ \leq 2(C_A+C_{A_0})C^{\frac12}(\delta(X_0)/2)^{\frac{1-d}2}\Big(\dint_{\Gamma_{\tau}}|\nabla u(Y)|^2\,dm(Y)\Big)^{\frac12}\\ \lra0\text{ as }\tau\ra0,
\end{multline}
where we have employed Lemma \ref{lm.greenfn} iv),  the fact that $u\in W$, and the absolute continuity of the integral. We have shown that
\[
\Phi_{\tau}(f)\lra\Phi(f),\qquad\text{for all }f\in C_c(\Delta_0),
\]
which gives the claimed weak convergence.

Now suppose that $f\in C_c(\Delta_0)$ and $\Vert f\Vert_{L^{q'}(\Delta_0,\sigma)}=1$, where $q'$ is the H\"older conjugate of $q$. Since $\omega_{\tau}^{X_0}\in RH_q(\Delta_0)$ uniformly in $\tau$, then we have that (\ref{eq.poissonlpest}) holds with $\Delta$ replaced by $\Delta_0$ and $X_{\Delta}$ replaced by $X_0$. Consequently, we see that
\begin{multline*}
|\Phi(f)|=\Big|\int_{\Delta_0}f\,d\omega^{X_0}\Big|=\Big|\lim_{\tau\ra0}\int_{\Delta_0}f\,d\omega_{\tau}^{X_0}\Big|\leq\sup_{\tau>0}\Big|\int_{\Delta_0}fk_{\tau}^{X_0}\,d\sigma\Big|\\ \leq\sup_{\tau>0}\Big(\Vert k_{\tau}^{X_0}\Vert_{L^q(\Delta_0,\sigma)}\Vert f\Vert_{L^{q'}(\Delta_0,\sigma)}\Big)\leq C_0\sigma(\Delta_0)^{1-q}.
\end{multline*}
It follows that $\Phi$ is a bounded linear functional on $L^{q'}(\Delta_0,\sigma)$. Hence we must have that $\omega^{X_0}\ll\sigma$ and $k^{X_0}=\frac{d\omega^{X_0}}{d\sigma}\in L^q(\Delta_0,\sigma)$ satisfies the estimate (\ref{eq.poissonlpest}), so that $k^{X_0}\in RH_q(\Delta_0)$. Since $\Delta_0\subset\Gamma$ was arbitrary, we finally conclude that $\frac{d\omega}{d\sigma}\in RH_q$.\hfill{$\square$}

Using the previous lemma in conjunction with Proposition \ref{prop.step4prop}  yields the desired conclusion for the operator $L$ and Theorem \ref{thm.main} is shown.\hfill{$\square$}

\section{Proof of Theorem \ref{thm.main2}}\label{sec.smallthm}

In this section, we give the proof of Theorem \ref{thm.main2}; in fact, the proof of this theorem is very similar to but simpler than the proof of Theorem \ref{thm.main} in the previous section, as we will not need to use the projection lemma, Lemma \ref{lm.djkg}, nor the extrapolation theorem, Theorem \ref{thm.extrapolation}. As such, we mainly describe the set-up here and point out the differences to the proof of Theorem \ref{thm.main2}.

Let $A_0$, $A$ be two matrices that satisfy (\ref{eq.elliptic}), write $\m A_0=w^{-1}A_0$, $\m A=w^{-1}A$, and suppose that $d\lambda(X)=\frac{\f a(X)^2}{\delta(X)^{n-d}}dX$ is a (continuous) Carleson measure with $\vertiii{\lambda}_{\m C}\leq\ep_0$, where $\f a$ is defined in (\ref{eq.disagreement}) and $\ep_0$ is small and to be chosen later. As in Lemma \ref{lm.todyadic}, the natural discretization of the Carleson measure $\lambda$ is the collection $\f m=\{\alpha_Q\}_{Q\in\bb D}$ with $\alpha_Q$ as defined in (\ref{eq.discretecarleson}). By Lemma \ref{lm.todyadic}, we see that $\Vert\f m\Vert_{\m C}\lesssim\ep_0$. Let $L_0=-\dv A_0\nabla$,  $L=-\dv A\nabla$, and let $\omega_0$, $\omega$ be the harmonic measures of $L_0$, $L$ respectively. 

Fix $Q_0\in\bb D$, and recall that $R_{Q_0}\subset B(x_{Q_0},7\sqrt nA_0\ell(Q_0))$ (see (\ref{eq.sawtoothbdd})). Let $\hat B_0:=B(x_{Q_0},30c^{-1}\sqrt n A_2\ell(Q_0))$, and fix $X_0\in\Omega$ as a Corkscrew point (with Corkscrew constant $c$) for the surface ball $\hat B_0\cap\Gamma$. Then $X_0\in\Omega\backslash B(x_{Q_0},30\sqrt nA_2\ell(Q_0))\subset\Omega\backslash R_{Q_0}^{**}\subset\Omega\backslash R_{Q_0}$. Moreover, according to Corollary \ref{cor.dd}, we have that $\omega^{X_0}$ is dyadically doubling in $Q_0$, while we also have that $\delta(X_0)\approx\dist(X_0,Q_0)\approx\ell(Q_0)$. 
 
Since $\operatorname{(D)}_{p'}$ is solvable for $L_0$, then by Theorem \ref{thm.equivalence}, we have that $\omega_{L_0}^{X_0}=\omega_0^{X_0}\in RH_p(Q_0)$ (with reverse H\"older characteristic independent of $Q_0$).

{\bf Step 0.}  Owing to our assumptions and Theorem \ref{thm.main}, we a priori have that $\omega_L\ll\sigma$. However, we still ought to make a qualitative reduction as in Step 0 of the proof of Theorem \ref{thm.main}, to guarantee that $k\in L^p_{\loc}(\Gamma,\sigma)$ for the fixed $p$ in our hypothesis. Accordingly, we define $A_{\tau}$ as $A_{\tau}=A_0$ in the $\tau-$tube $\Gamma_{\tau}$ (see Definition \ref{def.tube}), and $A_{\tau}=A$ in $\bb R^n\backslash\Gamma_{\tau}$. We work with $L_{\tau}$ in place of $L$ in the steps below, and ultimately pass to the limit as $\tau\ra0$ using Lemma \ref{lm.opconv}, very similarly as in Step 5 of the proof of Theorem \ref{thm.main}. We omit further details.

{\bf Step 1.} Introduce the operator $L_1$ defined as $L_1=L$ in $\hat B_0$, and $L_1=L_0$ in $\Omega\backslash\hat B_0$. We write $\omega_1$ for the harmonic measure associated to the operator $L_1$, and $g_1$ for the Green function associated to $L_1$. We have that $k_0^{X_0}\in RH_p(Q_0)$, and in particular by Theorem \ref{thm.equivalence}, Harnack Chains and the Harnack Inequality, we have that  
\begin{equation}\label{eq.poissonrh02}
\int_{Q_0}\big(k_0^{X_0}\big)^p\,d\sigma\lesssim\sigma(Q_0)^{1-p}.
\end{equation}

Our immediate goal in Step 1 is to show that (\ref{eq.poissonrh02}) remains true when $k_0^{ X_0}$ is replaced by $k_1^{X_0}$, the Poisson kernel for the operator $L_1$ defined above. The proof follows essentially the same as that of Step 1 in the previous section, where $\m F=\varnothing$ in our situation; thus we omit it and present only the 
\begin{conclusion}[Step 1] The estimate (\ref{eq.poissonrh02}) holds with $k_0^{X_0}$ replaced by $k^{X_0}_1$.
\end{conclusion}

We remark that the smallness of $\ep_0$, used in the analogue of (\ref{eq.hide}), necessarily depends on the implicit constant in Theorem \ref{thm.djk} applied to the operator $L_0$; hence, $\ep_0$ depends on the $RH_p$ constant of $\omega_0$.

{\bf Step 2.} Now let $L_2=L$ in $\Omega_{\hat B_0}$ and $L_2=L_1$ in $\hat B_0$. Note that $L_2\equiv L$ in $\Omega$. According to our choice of $\hat B_0$ and Lemma \ref{lm.poissonequiv}, we have that $k_2^{X_0}(y)\approx k_1^{X_0}(y)$ for $\sigma-$almost every $y\in B(x_{Q_0},7\sqrt nA_0\ell(Q_0))\supset Q_0$. As such, we reckon the estimate
\[
\int_{Q_0}(k_2^{X_0})^p\,d\sigma\approx\int_{Q_0}(k_1^{X_0})^p\,d\sigma\lesssim\sigma(Q_0)^{1-p}.
\]
Let $X_{Q_0}\in\Omega$ be a Corkscrew point for $Q_0$. Then by Lemma \ref{lm.poissonpole} and the doubling property of the harmonic measure, we have that $k_2^{X_0}\approx k_2^{X_{Q_0}}$ for $\sigma-$almost every $y\in Q_0$. Therefore, we conclude that
\[
\int_{Q_0}(k^{X_{Q_0}})^p\,d\sigma\lesssim\sigma(Q_0)^{1-p}.
\]
Since $Q_0\in\bb D$ was arbitrary, then the Poisson kernel $k$ for the operator $L$ satisfies the condition (\ref{eq.dyadicanalogue}) in Remark \ref{rm.equivrh}. By this same remark, we know that then $k$ satisfies the condition (\ref{eq.poissonlpest}), and finally by Theorem \ref{thm.equivalence} we obtain the desired conclusion.\hfill{$\square$}

\appendix

\section{Sets of locally finite perimeter and the relative isoperimetric inequality}

In this appendix, we state a minimal set of well-known results in geometric measure theory with the goal of verifying that the boundary of our sawtooth domain satisfies a relative isoperimetric inequality, which is used in the proof of Theorem \ref{thm.sawtooth}.  We follow the results in \cite{eg1} Chapter 5.
Let $U$ be an open subset of $\bb R^n$. By $C_c^1(U;\bb R^n)$ we denote the space of compactly supported $C^1$ vector fields $\phi:U\ra\bb R^n$.  

\begin{definition}[Sets of finite perimeter and related notions]\label{def.bv} A function $f\in L^1(U)$ has \emph{bounded variation} in $U$ if
\[
\sup\Big\{\dint_Uf\dv\varphi\,dX~:~\phi\in C_c^1(U;\bb R^n), |\varphi|\leq1\Big\}<\infty.
\]
A function $f\in L^1_{\loc}(U)$ has \emph{locally bounded variation} in $U$ if for each open set $V\subset\subset U$,
\[
\sup\Big\{\dint_Vf\dv\varphi\,dX~:~\phi\in C_c^1(V;\bb R^n), |\varphi|\leq1\Big\}<\infty.
\]
If $f$ has bounded variation (resp. locally bounded variation), we say that $f\in BV(U)$ (resp. $f\in BV_{\loc}(U)$). An $\n L^n-$measurable subset $E\subset\bb R^n$ has \emph{finite perimeter} (resp. \emph{locally finite perimeter}) in $U$ if ${\bf 1}_E\in BV(U)$ (resp. ${\bf 1}_E\in BV_{\loc}(U)$.
\end{definition}

\begin{theorem}[Structure Theorem for $BV_{\loc}$ functions; \cite{eg1} 5.1 Theorem 1]\label{thm.mu} Let\\ $f\in BV_{\loc}(U)$. Then there exists a Radon measure $\mu$ on $U$ and a $\mu-$measurable vector field $\sigma:U\ra\bb R^n$ such that
\begin{enumerate}[(i)]
	\item $|\sigma(x)|=1$ for $\mu-$almost every $x$, and
	\item $\dint_Uf\dv\varphi\,dX=-\dint_U\varphi\cdot\sigma\,d\mu$
\end{enumerate}
for all $\varphi\in C_c^1(U;\bb R^n)$.
\end{theorem}

If $E$ is a set of locally finite perimeter in $U$, we write $\Vert\partial E\Vert$ for  $\mu$, and we write $\nu_E:=-\sigma$, where $\mu$ and $\sigma$ are given in Theorem \ref{thm.mu}.

\begin{definition}[Measure theoretic boundary] Let $E$ be an $\n L^n-$measurable subset of $\bb R^n$, and let $x\in\bb R^n$. We say that $x\in\partial_\star E$, the \emph{measure theoretic boundary} of $E$, if
\[
\limsup_{r\ra0}\frac{\n L^n(B(x,r)\cap E)}{r^n}>0,\qquad\text{and}\qquad\limsup_{r\ra0}\frac{\n L^n(B(x,r)\backslash E)}{r^n}>0.
\]
\end{definition}

It is obvious that $\partial_\star E\subseteq\partial E$. We have the following criterion for a set to have locally finite perimeter.
\begin{theorem}[Criterior for finite perimeter, \cite{eg1} 5.11 Theorem 1]\label{thm.criterionfp} Let $E\subset\bb R^n$ be $\n L^n-$measurable. Then $E$ has locally finite perimeter if and only if for any compact set $K\subset\bb R^n$, $\n H^{n-1}(K\cap\partial_\star E)<\infty$.
\end{theorem}

\begin{definition}[Reduced boundary] Suppose that $E$ is a set of locally finite perimeter in $\bb R^n$, and let $x\in\bb R^n$. We say that $x\in\partial^\star E$, the \emph{reduced boundary} of $E$, if
\begin{enumerate}[(i)]
\item $\Vert\partial E\Vert(B(x,r))>0$ for all $r>0$.
\item $\lim_{r\ra0}\dint_{B(x,r)}\nu_E\,d\Vert\partial E\Vert=\nu_E(x)$, and
\item $|\nu_E(x)|=1$.
\end{enumerate}
\end{definition}

In \cite{eg1}, it is seen that for any $x\in\partial^*E$,
\[
\liminf_{r\ra0}\frac{\n L^n(B(x,r)\cap E)}{r^n}>0,\qquad\liminf_{r\ra0}\frac{\n L^n(B(x,r)\backslash E)}{r^n}>0,
\]
which implies that $\partial^\star E\subset\partial E$ for any set $E$ of locally finite perimeter.

\begin{theorem}[Structure theorem for sets of finite perimeter, \cite{eg1} 5.7 Theorem 2]\label{thm.structure} \phantom{a} Assume that $E$ has locally finite perimeter in $\bb R^n$.
\begin{enumerate}[(i)]
\item Then
\[
\partial^\star E=\medcup_{k=1}^\infty K_k\cup N,
\]
where $\Vert\partial E\Vert(N)=0$ and $K_k$ is a compact subset of a $C^1-$hypersurface $S_k$ ($k=1,2,\ldots$).
\item Furthermore, $\nu_E|_{S_k}$ is normal to $S_k (k=1,\ldots)$, and
\item $\Vert\partial E\Vert=\n H^{n-1}\mres\partial^\star E$.
\end{enumerate}
\end{theorem}

Finally, we state the relative isoperimetric inequality.

\begin{theorem}[Relative isoperimetric inequality, \cite{eg1} 5.6 Theorem 2]\label{thm.iso} Let $E$ be a set of locally finite perimeter in $\bb R^n$. Then for each (open) ball $B(x,r)\subset\bb R^n$, we have that
\begin{equation}\label{eq.iso}
\min\big\{\n L^n(\overline{B(x,r)}\cap E)~,~\n L^n(\overline{B(x,r)}\backslash E)\big\}^{1-\frac1n}\leq a_n^{-1}\Vert\partial E\Vert(B(x,r))
\end{equation}
where $a_n$ is a uniform constant depending only on the dimension $n$.
\end{theorem}


\hypersetup{linkcolor=toc}

\bibliography{refs} 
\bibliographystyle{alphaurl-max} 


\end{document}